\documentclass[12pt]{article}
\usepackage{multirow}
\usepackage{latexsym, amsbsy, amssymb}
\usepackage{amsmath}
\usepackage{relsize}
\usepackage[geometry]{ifsym}
\usepackage{graphicx}

\usepackage{xcolor}

\usepackage{natbib}

\usepackage{hyperref}

\usepackage{subfig}
\newcommand{\stens}[1]{\mbox{\rsfstena #1}}

\newtheorem{proposition}{{\bf Proposition}}[section]
\newtheorem{lemma}{{\bf Lemma}}[section]
\newtheorem{corollary}{{\bf Corollary}}[section]
\newtheorem{assumption}{{\bf Assumption}}[section]
\newtheorem{definition}{{\bf Definition}}[section]
\newtheorem{theorem}{{\bf Theorem}}[section]

\newtheorem{remark}{{\bf Remark}}[section]

\DeclareMathOperator*{\argmin}{arg\,min}

\newtheorem{example}{{\bf Example}}[section]

\def\cid{\stackrel{\mbox{\tiny $\cal D$}}\longrightarrow}
\def\cip{\stackrel{\mbox{\tiny $\cal P$}}\longrightarrow}
\def\argmin{\mathrm{argmin}}
\def\proof{\noindent {\sc Proof. }}

\def\cip{\stackrel{\mbox{\tiny $P$}}\rightarrow}

\usepackage{graphicx}
\usepackage{epsfig}
\usepackage{amsmath}

\begin{document}
\vspace{-18pt}

\def\R{\mathbb R}
\def\Xb{\boldsymbol X}
\def\trans{^{\mbox{\tiny{\sf  T}}}}
\def\Yb{\boldsymbol Y}
\def\eqdis{\stackrel{\mbox{\tiny $\cal D$}}{=}}
\def\half{^{\mbox{\tiny $\frac{1}{2}$}}}
\def\mhalf{^{\mbox{\tiny $-\frac{1}{2}$}}}
\def\inv{^{\mbox{\tiny $-1$}}}
\def\thetabfs{{\mbox{\tiny $\thetabf$}}}
\def\Xbfs{{\mbox{\tiny $\Xbf$}}}
\def\var{\mathrm{var}}
\def\cov{\mathrm{cov}}
\newcommand{\indep}{\;\, \rule[0em]{.03em}{.67em} \hspace{-.27em}
\rule[-.02em]{.7em}{.03em} \hspace{-.27em}
\rule[0em]{.03em}{.67em}\;\,}
\def\tr{\mathrm{tr}}
\def\vec{\mathrm{vec}}
\def\mat{\mathrm{mat}}
\def\sign{\mathrm{sign}}
\def\S{\mathbb S}
\def\eop{\hfill $\Box$ \\
}
\def\trans{^{\mbox{\tiny{\sf T}}}}

\def\Xs{{\mbox{\tiny $X$}}}
\def\cid{\stackrel{\mbox{\tiny $\cal D$}}\rightarrow}
\def\gpc{{\cal G}}

\def\dy{D_{{\mbox{\tiny $Y$}}}}
\def\du{D_{{\mbox{\tiny $U$}}}}
\def\kappax{\kappa_{\mbox{\tiny $X$}}}
\def\Xsmall{\mbox{\tiny $X$}}
\def\Ysmall{\mbox{\tiny $Y$}}
\def\yygx{{\mbox{\tiny $YY|X$}}}
\def\xxsub{{\mbox{\tiny $XX$}}}
\def\yysub{{\mbox{\tiny $YY$}}}
\def\xysub{{\mbox{\tiny $XY$}}}
\def\yxsub{{\mbox{\tiny $YX$}}}
\def\Psub{{\mbox{\tiny $P$}}}
\def\real{\mathbb R}
\def\nano{\scriptscriptstyle}
\def\trans{^{\nano \sf T}}
\def\sigvu{\Sigma_{\nano VU}}
\def\adj{^*}
\def\siguv{\Sigma_{\nano UV}}
\def\lambdamax{\lambda_{\max}}
\def\half{^{1/2}}
\def\key{K}
\def\omxi{\Omega_{\nano X^i}}
\def\lxi{L_2(P_{\nano X^i})}
\def\lxo{L_2(P_{\nano X^1})}
\def\lxp{L_2(P_{\nano X^p})}
\def\lu{L_2(P_{\nano U})}
\def\pu{P_{\nano U}}
\def\pv{P_{\nano V}}
\def\lu{L_2 (\pu)}
\def\lv{L_2 (\pv)}
\def\corr{\mathrm{cor}}
\def\lxi{L_2(P_{\nano X^i})}
\def\lxj{L_2(P_{\nano X^j})}
\def\iff{\, \Leftrightarrow \, }
\def\sign{\mathrm{sign}}
\def\equal{&=\,}
\def\xk{X^k}
\def\lpxi{L_2 (P_{\nano X^i})}
\def\pxi{P_{\nano X^i}}
\def\px{P_{\nano X}}
\newfont{\rsfsten}{rsfs10 scaled 1100}
\newfont{\rsfstena}{rsfs10 scaled 800}
\newfont{\rsfstenb}{rsfs10 scaled 800}
\def\ox{\Omega_{\nano X}}
\def\oy{\Omega_{\nano Y}}
\def\a{\mbox{\rsfsten A\,}}
\def\as{{\mbox{\rsfstena A}}}
\def\prob{{\cal P}}
\def\lpx{L_2(P_{\nano X})}
\def\h{\mbox{\rsfsten H\, }}
\def\hs{\mbox{\rsfstena H}}
\def\ep{E_{\nano P}}
\def\graph{{\cal G}}
\def\xmij{X^{-(i,j)}}
\def\moment{\mbox{\rsfsten M}}
\def\agm{\mathrm{AGM}}
\def\xa{X^{\nano A}}
\def\prob{\mbox{\rsfsten P}}
\def\k{\mbox{\rsfsten K\,}}
\def\hs{{\mbox{\rsfstena H}}}
\def\hx{\h_{\nano X}}
\def\hy{\h_{\nano Y}}
\def\hz{\h_{\nano Z}}
\def\xb{X^{\nano B}}
\def\pxb{P_{\nano \xb}}
\def\hb{\h_{\nano B}}
\def\ha{\h_{\nano A}}
\def\px{P_{\nano X}}
\def\proba{\prob_{\nano A}}
\def\qxa{Q_{\nano \xa}}
\def\implies{$\Rightarrow$}
\def\moment{\mbox{\rsfsten M}}
\def\sigmaij{\Sigma_{\nano X^i X^j}}
\def\xa{X^{\nano A}}
\def\pxa{P_{\nano \xa}}
\def\sadditive{\mbox{\rsfstenb A}}
\def\vbar{\rule[-.21em]{0.03em}{0.96em}}
\def\twohbar{\rule[-.23em]{.35em}{.03em}\hspace{-.35em}\rule[.75em]{.35em}{.03em}}
\def\twovbar{\vbar\,\vbar}
\def\hiba{\rule[.67em]{.45em}{.03em}}
\def\leftdouble{\twohbar\hspace{-.35em}\twovbar\hspace{.2em}}
\def\rightdouble{\hspace{.2em}\twohbar\hspace{-.15em}\twovbar}

\def\cran{\overline{\mathrm{ran}\!}\, \, }

\newcommand\hi[1]{^{\mbox{\raisebox{1.4pt}{$\scriptscriptstyle #1$}}}}

\def\mp{m_{\nano P}}

\vspace{.3in}
\def\inner{\mbox{\rsfsten U}}
\def\ox{\Omega_{\nano X}}
\def\lpxj{L_2(P_{\nano X^j})}
\def\lpxj{L_2 (P_{\nano X^j})}
\def\sigmaxx{\Sigma_{\nano XX}}
\def\bdd{{\mbox{\rsfsten B}}}
\def\lp{L_2(P)}
\newcommand\lt[1]{L_2(P_{\nano #1})}
\def\piz{\Pi_{\nano Z}}
\def\lpz{L_2 (P_{\nano Z})}
\def\f{\mbox{\rsfsten F}}
\def\fx{\f_{\nano X}}
\def\fy{\f_{\nano Y}}
\def\oz{\Omega_{\nano Z}}
\def\pz{P_{\nano Z}}
\newcommand{\sten}[1]{\mbox{\rsfsten #1}\,}
\newcommand{\sub}[2]{#1_{\nano #2}}
\def\additive{\sten A}
\def\field{{\cal F}}
\def\of{{\scriptstyle \circ}}
\def\aindep{\indep_{\nano A}}
\def\dindep{\perp_{\nano d}}

\newcommand\lo[1]{_{\nano  #1}}

\newcommand\projminus[2]{P_{{\stens A}_{\nano #1}{\nano \ominus}{\stens A}_{\nano #2}}}

\noindent

\def\real{\mathbb R}
\def\tsum{
%\textstyle{
\sum}
% }

\parindent 0cm

\def\ail{\sten A}

\def\sfG{\mathsf{G}}
\def\sfE{\mathsf{E}}
\def\sfV{\mathsf{V}}

\def\HS{\mathrm{HS}}
\begin{center}
{\Large{\bf
Nonparametric and high-dimensional \\
functional graphical models }} \\
\vspace{.2in}
Eftychia Solea and Holger Dette
%\vspace{.2in}

{\em Crest-Ensai and Ruhr-Universit{\"a}t Bochum}
\end{center}

\def\nano{\scriptscriptstyle}
\def\of{\mbox{\raisebox{1pt}{$\nano{\circ}$}}}
\def\ka{\kappa}
\def\cspan{\overline{\mathrm{span}}}

\def \M {\mathfrak{M}}

 \abstract{We consider the problem of constructing nonparametric undirected graphical models for high-dimensional functional data.  Most existing statistical methods in this context assume either a Gaussian distribution on the vertices or linear conditional means.  In this article we provide a more flexible model which relaxes the linearity assumption by replacing it by an arbitrary additive form.  The use of  functional principal components offers an estimation strategy that uses a group lasso penalty to estimate the relevant edges of the graph.   We establish  statistical guarantees for the resulting estimators, which can be used to prove consistency if the dimension and the number of functional principal components  diverge to infinity with the  sample size.  We also investigate the empirical performance of our method through simulation studies and a real data application. \\
 

}

\section{Introduction} \label{sec1}
\def\theequation{1.\arabic{equation}}
\setcounter{equation}{0}

In recent years, there has been a large amount of work on estimating undirected graphical models that describe the conditional dependencies among the components of a $p$-dimensional random vector $X=(X \hi 1, \ldots, X \hi p ) \trans$.  Let $\sf V=\{ 1, \ldots, p \}$, and $\sf E$ denote a subset of $\{ (i,j) \in \sf V \times \sf V: \, i \ne j \}$, which satisfies $(i,j) \in \sf E$ if and only if $(j,i) \in \sf E$. The pair $\sf G = ( \sf V, \sf E)$ constitutes an undirected graph, with $\sf V$ representing the set of vertices and $\sf E$  the set of edges. The vector $X$ follows a graphical model  if
\begin{align}\label{eq: graphical model}
(i,j) \notin  {\sf E}   \quad \iff \quad X \hi i \indep X \hi j | X \hi {-(i,j)},
\end{align}
where $X \hi {-(i,j)}$ represents the vector $X$ with its $i$th and $j$th components removed, and for random elements $A$, $B$, and $C$,  $A \indep B | C$ means that  $A$ and $B$ are conditionally independent given  $C$.  The goal is to estimate the edge set $\sf E$ based on a random sample from $X$.

If  $X$ is assumed to follow a $p$-dimensional Gaussian distribution with expectation $0 \in \real \hi d$ and positive definite covariance matrix  $\Sigma \in \real \hi {p \times p}$, the model is called  {\it Gaussian Graphical Model (GGM)} and has become very popular.  For a Gaussian random vector $X=(X \hi 1, \ldots, X \hi p ) \trans$, the structure of the precision matrix $\Theta=\Sigma \hi {-1}$  characterizes the conditional independence relationships among the variables $X \hi 1, \ldots, X \hi p$ \citep{lauritzen1996graphical}. Specifically,

\begin{align}\label{gaussianrelation}
 X \hi i \indep X \hi j | X \hi {-(i,j)}  \quad \iff \quad \theta  \lo {ij}=0,
\end{align}
where $\theta \lo {ij}$ is the $(i, j)$th entry of the precision matrix $\Theta$.  Because of the relation \eqref{gaussianrelation}, the estimation of the edge set $\sf E$ reduces to estimating the sparsity pattern of the precision matrix $\Theta$. Hence, there exists a large amount of literature, which has  its focus  on estimating high-dimensional Gaussian graphical models. 
 For example, \cite{meinshausen2006high} introduced a neighbourhood-based approach by solving $p$ lasso linear regression problems for each node of the graph.  \cite{yuan2007} and \cite{friedman2008sparse} considered a penalized maximum likelihood approach with the lasso penalty imposed on the off-diagonal entries of the precision matrix $\Theta$.   Based on a relation between partial correlation and regression coefficient,  \cite{peng2009partial} proposed to estimate a sparse GGM by imposing the lasso penalty on the partial correlations.  Other developments on GGM include the SCAD and the adaptive lasso penalty
 \citep{lam2009sparsistency},  the Dantzig selector  \citep{cai2011constrained}  and  hard-thresholding \citep{bickel2008covariance}.

Despite of its simplicity, the Gaussian assumption can be very restrictive in practice and statistical inference based on the Gaussian distribution might be misleading if this assumption is violated. Therefore, more recent work has its focus on  considering graphical models under less restrictive assumptions. For example, \cite{liu2009nonparanormal, liu2012high} and \cite{xue2012regularized} relaxed the marginal Gaussian assumption on the vertices of the graph using copula transformations, and  \cite{voorman2013graph} allowed the conditional means of the variables to take an additive form.   \cite{liaci} and \cite{lee2016additive} proposed a non-Gaussian graphical model based on  additive conditional independence (ACI), a three-way statistical relation  that captures the spirit of conditional independence.
% but without resorting to high-dimensional kernels.

Most of the literature discussed so far has its focus on graphical models for finite dimensional data.  However, many recent applications involve functional data, such as electroencephalogram (EEG) and functional magnetic resonance imaging (fMRI) data, where each sampling unit is modelled as a realization of a stochastic process varying over a time interval.  In this paper, we are interested in estimating a nonparametric and high-dimensional undirected graphical model for multivariate functional data.  

In contrast to the finite dimensional case, less literature can be found  on  graphical models for multivariate functional data.  
 \cite{qiao2018functional} proposed the Functional Gaussian Graphical Model (FGGM)  assuming that   $X$ is a multivariate Gaussian random process. 
Roughly speaking, they used  a truncated Karhunen-Lo{\`e}ve expansion, say   of order $m_{n}$,  to  reduce the infinite dimensional problem  to 
 a  $pm_{n}$-dimensional problem for the principal component scores.  The conditional independencies of the graph define  a  block sparsity  structure, such that 
the properties of the precision matrix of the scores can be used  to identify the edge set using a group lasso penalty. They called this method functional glasso, or simply fglasso, and 
the authors   showed that,  when $m \lo n$ approaches infinity, consistent estimation of the edge set is possible. 
%The FGGM can also be represented as a multivariate linear regression model with respect to the principal scores (see Example \ref{ex1} below  for details). 
\cite{zhu2016bayesian}  proposed a Bayesian framework under the Gaussian assumption on the random functions for the analysis of functional graphical models, while  \cite{li2018nonparametric} relaxed this assumption by extending the concept of ACI to the functional setting. 

  \begin{figure}[h]
  \begin{center}
\includegraphics[scale=0.45]{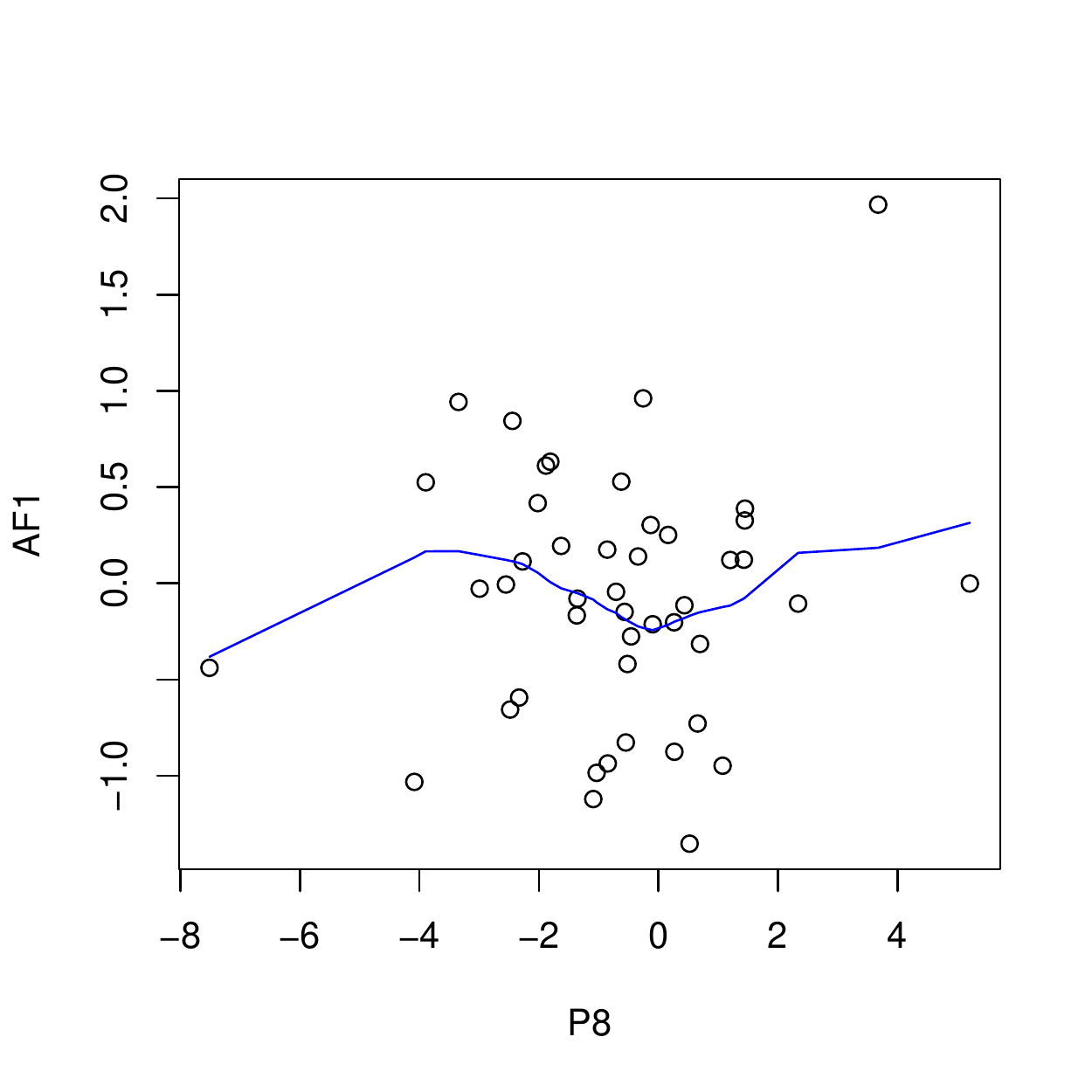} 
 \includegraphics[scale=0.45]{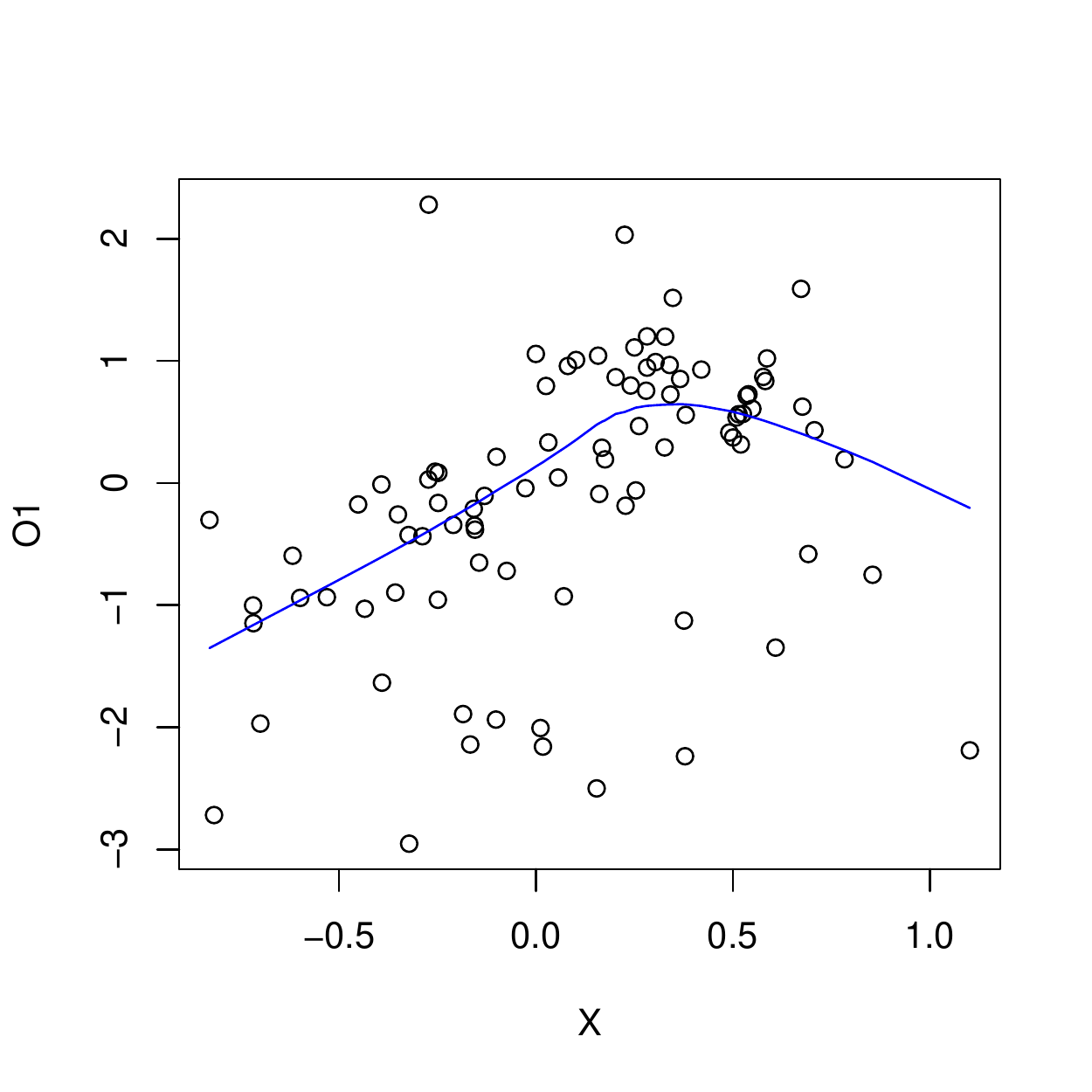}
 \caption{\it Pairwise scatterplots  for the control group  between channels AF1 and P8 (left) and channels O1 and X (right).  \label{fig1}}
   \end{center}
\end{figure}

{
 In this paper, we introduce an alternative approach to relax the Gaussian assumption in the functional graphical model. Our research is motivated by the fact that in many applications the relation between the functional principal scores is rarely linear as implied by the assumption of the FGGM. To illustrate this observation, we consider an electroencephalography (EEG) dataset that consists of two groups of subjects: 77 subjects in the alcoholic group, and 45 in the control group \citep{zhang1995event, ingber1997statistical}.   For each subject, an EEG activity was recorded at 256 time points over a one second time interval using 64 electrodes placed on the subject's scalp.  The goal is to construct a functional graphical model to characterize brain network connectivity for the two groups of subjects based on the functional data collected by the electrodes.  
In  Figure \ref{fig1},  we display the pairwise scatterplots between channels using the first two principal components for the random functions of the control group. Clearly, this figure
 indicates that the conditional relationships among the scores corresponding to different vertices of the graph are nonlinear. 
 Therefore, the Gaussian assumption of the FGGM is difficult to justify for the analysis of this type of data. 
 
As an alternative, we propose a new nonparametric functional graphical model that allows the conditional relationships among the principal scores to take an additive structure.
 Our approach uses  the traditional probabilistic concept of conditional independence,  and then applies the additive structure to the scores of   the Karhunen-Lo{\'e}ve expansion of   each random function. We  approximate each nonparametric additive component by a linear combination of B-splines basis functions. This enables us to  estimate the edge set of the graph  by imposing the group lasso penalty on a matrix formed by the coefficients in the spline approximation. 
We derive statistical guarantees for  the resulting estimates, which can be used to prove consistency if the dimension $p$ and the number of scores   diverge to infinity with the sample size. This provides a useful methodology for general nonparametric analysis of high-dimensional functional graphical models.  Our approach differs from related work  on nonparametric functional graphical models by   \cite{leefaro} and  \cite{li2018nonparametric}, who used additive conditional independence and reproducing kernel Hilbert spaces for  neighbourhood selection via a functional additive regression operator.

The remainder of the article is organized as follows. Section \ref{methodology} describes the methodology and proposes the nonparametric functional graphical model. Section \ref{estimation} presents the estimation procedure. In Section \ref{theory} we study the theoretical  properties of the resulting estimator. In Section \ref{simulations} we conduct simulation studies to evaluate the finite sample properties of the proposed methodology, and in Section \ref{sec52} we apply the new model to the motivating EGG dataset. We conclude with some final remarks in Section \ref{conclusion}, while all proofs of the theoretical results are deferred to the Appendix.

\section{Additive functional graphical models }\label{methodology}
\def\theequation{2.\arabic{equation}}
\setcounter{equation}{0}

We first provide a formal definition of an {\it additive function-on-function regression model} which will be used to define the functional graphical model considered in this paper.  We begin by introducing some basic concepts from functional data analysis.

Throughout this paper $\mathcal{L}  \hi 2([0,1])$ denotes the space of all square-integrable functions defined on the interval $[0,1] \subset \R$. We denote by $\langle f ,  g \rangle=\int \lo {[0,1]} f(t) g(t) dt$ the common inner product in $\mathcal{L}  \hi 2([0,1])$ and by $\| f \|=\langle f, f \rangle \hi {1/2}$ the corresponding norm.  
 Let $X=(X \hi 1, \ldots, X \hi p ) \trans$ denote a $p$-dimensional random element with mean $0$  whose $i$th component $X \hi i $ is an element of $\mathcal{L}  \hi 2([0,1])$ such that $E\| X \hi i \| \hi 2< \infty$. 
For each $X^i$, we define the corresponding  covariance operator  
\begin{align}\label{covop}
\Sigma \lo {X \hi i X \hi i}  (f) (t) = \int \lo T f(s) \sigma \lo {X \hi i X \hi i} (s, t) d s, \quad f \in \mathcal{L}  \hi 2([0,1]),
\end{align}
where ${\mathcal{\sigma}   \lo {X \hi i X \hi i}} (s,t)= \cov(X  \hi i (s), X \hi i (t))=E(X  \hi i (s) X \hi i (t))$ is  the covariance function of the random element $X \hi i$.  The operator $\Sigma \lo {X \hi i X \hi i}$ is a compact  Hilbert-Schmidt operator \citep[see, for example,][]{hsing2015theoretical}, and there exists a spectral decomposition of the covariance function of the form
\begin{align}\label{trueeigenpairs}
\sigma \lo {X \hi i X \hi i} (s, t) =  \tsum \lo {r=1}  \hi {\infty} \lambda \hi i \lo r \phi \hi i \lo r (s) \phi \hi i \lo r(t),
\end{align}
where $\lambda \hi i \lo 1 \geq  \lambda \hi i \lo 2 \geq \ldots$ are the eigenvalues and $\{\phi \hi i \lo k \} \lo {k \in \mathbb{N}}$ are orthonormal eigenfunctions satisfying
\begin{align*}
\int \lo T  \sigma \lo {X \hi i X \hi i } (s, t) \phi \lo r \hi i (s) d s= \lambda \lo r \hi i \phi \lo r \hi i (t).
\end{align*}   
Consequently, each $X \hi i \in \mathcal{L}  \hi 2([0,1])$ can be represented by its Karhunen-Lo{\'e}ve expansion
 \begin{align} \label{karhunen}
 X \hi i = \tsum \lo {r=1}  \hi {\infty} \sqrt{\lambda \hi i \lo r} \xi \lo {r} \hi i  \, \phi \lo r \hi i \quad i=1, \ldots, p,
\end{align}
where the random variables  $\xi \lo {r} \hi i=\langle X \hi i, \phi \lo r \hi i \rangle / \sqrt{\lambda \hi i \lo r}$ are called the {\it {functional principal component scores}} and satisfy  $E(\xi \lo {r} \hi i)=0$,  $\var( \xi \lo {r} \hi i)= 1$, $E(\xi \hi i \lo {q} \xi \hi i \lo {r})=0$  for $q \ne r$.  

We next give our formal definition of the functional graphical model.
Suppose $\sf G = ( \sf V, \sf E)$ is an undirected graph,  where $\sf V$ denotes the finite set $\{ 1, \ldots, p \}$, and $\sf E$ denotes a subset of $\{ (i,j) \in \sf V \times \sf V: \, i \ne j \}$, which satisfies $(i,j) \in \sf E$ if and only if $(j,i) \in \sf E$.
\begin{definition}
A vector of random functions  $X=(X \hi 1, \ldots X \hi p) \hi {\top} \in  \mathcal{L}  \hi 2([0,1]) \times \ldots, \times \mathcal{L}  \hi 2([0,1])$ is said to follow a functional graphical  model  with respect to an undirected graph $\sf G = ( \sf V, \sf E)$ if and only if
\begin{align*}
X \hi i \indep   X \hi j | X \hi {-(i,j)}, \ \ \ \forall \, (i, j) \notin \sf E.
\end{align*}
\end{definition}

 {
 \begin{example} \label{ex1} 
 {\rm 
 \cite{qiao2018functional}  assumed that $X=(X \hi {1} , \ldots , X \hi {p}) \hi {\top }$ is a Gaussian process on  $\mathcal{L}  \hi 2([0,1]) \times  \ldots \times \mathcal{L}  \hi 2([0,1])$
 and define a Functional Gaussian Graphical model  (FGGM) by  the condition
\begin{align} \label{fggm}
(i,j) \notin {\sf E} \quad \iff \quad \mbox{cov} [ X \hi i (s), X \hi j (t)| X \hi {-(i,j)} ] =0  \; \; \;  \forall \ s, t \in [0,1].
\end{align}
They  proposed to  approximate each $X \hi i$ by the first $m \lo n$ coefficients from the Karhunen-Lo{\'e}ve  expansion  \eqref{karhunen}. Thus 
 for each $X \hi i $, one obtains  a $p m \lo n$-dimensional Gaussian random vector $\xi \hi {\top}=((\xi \hi { 1}) \hi {\top}, \ldots, (\xi \hi { p}) \hi {\top})$ of scores, where $\xi \hi {i}=(\xi \hi { i} \lo 1, \ldots, \xi \hi { i} \lo {m \lo n})^{\top}$ is the vector of the first $m \lo n$ functional principal  {component scores in the 
 Karhunen-Lo{\`e}ve expansion  \eqref{karhunen} of each $X \hi i$.  
% 
% infinity with increasing sample size.
% Note that the FGGM is the extension of the Gaussian graphical model of \cite{yuan2007} to the functional setting.  \\
 Using the  $m \lo n$-truncation 
 \cite{qiao2018functional}  also 
 showed that the FGGM can be represented as a conditional multivariate linear regression model with respect to the scores.  Indeed, each $\xi \hi { i } \lo {q}$ can be expressed as
\begin{align} \label{flm}
\xi \hi { i} \lo {q}  =\tsum \lo{j \ne i} \hi {p}  \tsum \lo {r=1} \hi {m \lo n}   B  \hi {ij}  \lo {qr}    \xi \hi j \lo r + \epsilon \hi i \lo q, \quad i \in {\sf V}, q=1, \ldots, m \lo n,
\end{align}
such that $ ({\epsilon \hi i \lo q}) \lo {1 \leq q \leq m \lo n}$ is uncorrelated with $ (\xi \hi j \lo r)  \lo {1 \leq r \leq m \lo n}, i \ne j$ if and only if 
 \begin{align*}
B \hi {ij} \lo n=( B  \hi {ij}  \lo {qr} ) \lo {1 \leq q ,r \leq m \lo n}=- (\Theta \hi {ii} \lo n) {\hi {-1}} \Theta \hi {ij} \lo n, \quad (i, j) \in {\sf V} \times {\sf V}, i \ne j,
\end{align*}
 where $ \Theta  \lo {  n} \hi {ij} \in \R \hi {m \lo n \times m \lo n}$ is the $(i, j)$th element of the block precision matrix $\Theta  \lo {  n}=(\Theta \hi {ij}  \lo {  n}) \lo {1 \leq i, j \leq p} \in \R \hi {p m \lo n \times p m \lo n}$ of  the $pm_{n}$-dimensional vector $\xi$.  
Hence, under the Gaussian assumption the conditional relationships between nodes $i$ and $j$ are linear, and the network structure  of the FGGM can also be recovered by the sparsity structure of the regression coefficient matrix $B  \hi {ij} \lo n $.  
They used  group-lasso penalized maximum likelihood estimation to address the  blockwise sparsity of the precision matrix
 and showed that  the  precision matrix  is a consistent estimate of the set  $\sf E $,  when $p$ and  $m \lo n$ approach infinity with increasing sample size.
 Note that the FGGM is the extension of the Gaussian graphical model of \cite{yuan2007} to the functional setting.
 }
}
\end{example}
}

We use a generalization of  the representation \eqref{flm} to  give a formal definition of the additive function-on-function model for multivariate functional data.

%$****************************definition *******************************************$

\begin{definition} Consider a vector  of random functions $X=(X \hi 1, \ldots X \hi p) \hi {\top} \in  \mathcal{L}  \hi 2([0,1]) \times \ldots, \times \mathcal{L}  \hi 2([0,1])$ and  
suppose that each $X \hi i$ has a Karhunen-Lo{\'e}ve  expansion of the form \eqref{karhunen}.
The  vector $X$ follows the function-on-function additive  model if for each pair $(i, j) \in \sf V \times \sf V$ there exists a sequence of smooth functions $f  \hi {ij}=\{ f  \hi {ij} \lo {qr}:q, r \in \mathbb{N} \}$ defined on $\R$ with $E  [f \hi {ij} \lo {qr}(\xi \hi j \lo {r} )]=0, q, r \in \mathbb{N}$, such that
  \begin{align}\label{afm}
E[\xi \hi { i} \lo {q}|\{ \xi \hi j \lo {r}, j \ne i\}] =\tsum \lo{j \ne i}  \hi {p}  \tsum \lo {r=1}  \hi {\infty}   f  \hi {ij}  \lo {qr}    (\xi \hi j \lo {r})
\end{align}
 \end{definition}
  Similar to the functional additive regression model of  \cite{han2018smooth}, model \eqref{afm} relaxes the linearity assumption in FGGM by imposing an additive structure on the scores in the  Karhunen-Lo{\'e}ve expansion, giving rise to a more flexible model than the FGGM.  By definition, the scores $\xi \hi i \lo r$ are uncorrelated, but we also require them to be independent in the following discussion as also postulated in  \cite{han2018smooth}.  Furthermore, we assume that they take values in a closed and bounded interval $[-1, 1]$. For example, this can be achieved  by taking a monotone transformation $\Psi: \R \to [-1, 1]$   \citep[see][]{zhu2014structured,wong2019partially}. 

% \subsection{Additive functional graphical model}

We now define a new nonparametric functional graphical model which we call the Additive Functional Graphical Model (AFGM).  
{\begin{definition}
Suppose  $X=(X \hi 1, \ldots X \hi p) \hi {\top} \in  \mathcal{L}  \hi 2([0,1]) \times \ldots, \times \mathcal{L}  \hi 2([0,1])$ is associated with a functional graphical model $\sf G = ( \sf V, \sf E)$. If $X$ is additionally a function-on-function additive model of the form \eqref{afm},  then we say that $X$ follows an additive functional graphical model, and write this statement as $X \sim \mathrm{AFGM}(\sf G)$.
\end{definition}}

The definition implies that the independence structure of $X$ can be recovered by the sparse structure of the additive components $f \lo {qr} \hi {ij}$ in the representation \eqref{afm}.  
 Since each random function is infinite-dimensional, some type of regularization is needed by truncating the Karhunen-Lo{\'e}ve  expansion  \eqref{karhunen} at a finite number of principal components, say $m \lo n$, where $(m \lo n) \lo {n \in \mathbb{N}}$ is a sequence converging to infinity with increasing sample size.   
Thus, we obtain  a truncated version of model \eqref{afm}, that is  
\begin{align}\label{scorestrunc}
E[ \xi ^i \lo q  |  \{ \xi \hi j \lo {r}, j \ne i; r=1, \ldots , m_n \}] 
=\tsum \lo{j \ne i}  \hi {p}  \tsum \lo {r=1}  \hi {m \lo n}   f  \hi {ij}  \lo {qr}    (\xi \hi j \lo {r}) , \quad q=1, \ldots, m \lo n, i \in {\sf V}.
\end{align}
%\begin{align}\label{scorestrunc}
%\xi \hi { i} \lo {q} =\tsum \lo{j \ne i}  \hi {p}  \tsum \lo {r=1}  \hi {m \lo n}   f  \hi {ij}  \lo {qr}    (\xi \hi j \lo {r}) + \epsilon \lo q \hi i, \quad q=1, \ldots, m \lo n, i \in {\sf V},
%\end{align}
Then, our goal is to  estimate the ``truncated'' edge set
\begin{align}\label{truncatededge}
{\sf {E} \lo n} =\{(i, j) \in {\sf {V}} \times {\sf {V}}: i \ne j, f \lo {qr} \hi {ij}  \ne 0 \,  \mathrm{~for \, some~} \, q, r=1, \ldots, {m \lo n} \}
\end{align}
 Note that we aim to recover the edge set when each $X \hi i$ is approximated by  a finite sum of $m \lo n$ terms rather than an  infinite sum. Our theoretical results  in Section \ref{theory} show that the edge set can be identified   with probability converging to $1$ as    $m \lo n \to \infty$, $p \to \infty$ and $n \to \infty$.

\begin{remark} \label{rem1} ~~

{\rm \begin{itemize} 
\item[(1)] Note that it is not necessary to fix the sign of the eigenfunctions in the definition of scores $\xi \hi i \lo r= \langle X \hi i, \phi \hi i \lo r \rangle / \sqrt{\lambda \lo r}$ used in the         representation \eqref{afm} or \eqref{scorestrunc},  because a sign change can always be compensated by choosing the function $f \hi {ij} \lo {qr}(-x)$ instead of $f \hi {ij} \lo {qr}(x)$.
\item[(2)]
Model \eqref{scorestrunc} can be regarded as the nonparametric and additive version of the FGGM \eqref{flm}, and the generalization of the model of \cite{voorman2013graph} to the functional setting where they propose a semi-parametric method  which allows the conditional means of the random variables to take on an arbitrary additive structure. 
 
 \end{itemize} }
 \end{remark}

 %%%%%%%%%%%%%%%%%%%%%%%%%%%%%%%%%%%%
 % Estimation
 %%%%%%%%%%%%%%%%%%%%%%%%%%%%
 
 \section{Estimation and computation}\label{estimation}
 \label{sec3}
\def\theequation{3.\arabic{equation}}
\setcounter{equation}{0}

  In this section, we develop an estimation procedure for fully observed functional data to estimate the scores $\xi \hi i$ for each $X \hi i$, which is used afterwards  for the estimation of the edge set $\sf E \lo n$.

% \subsection{Empirical  KarhunenÃ¢ÂÂLo{\`e}ve expansion}

To be precise, let   $X \lo 1, \ldots, X \lo n$ be an independent sample from  $X$, such that for each $u = 1, \ldots, n$,   $X \lo u= (X \lo u \hi 1 , \ldots, X \lo u \hi p )^{\top}$ is a vector in $ \mathcal{L}  \hi 2([0,1]) \times \ldots \times  \mathcal{L}  \hi 2([0,1])$. 
Then, for each $i = 1, \ldots, p$,  the covariance operator $ {\Sigma} \lo {X \hi i X \hi i}$ can be estimated by
\begin{align*}
\hat  \Sigma \lo {X \hi i X \hi i}  (f) (t) = \int \lo {[0,1]} f(s) \hat \sigma \lo {X \hi i X \hi i} (s, t) d s, \quad f \in \mathcal{L}  \hi 2([0,1]),
\end{align*}
where
\begin{align*}
 \hat \sigma \lo {X \hi i X \hi i } (s, t) = \frac{1}{n} \tsum \lo {u=1} \hi n    X \hi i \lo u (s)    X \hi i \lo u (t) 
\end{align*}
is the common estimator of the covariance function
(note that  the $X \hi i$ are centered).  Let $\hat \lambda \lo r \hi i$ and $\hat \phi \lo r \hi i$ be the sample eigenvalues and eigenfunctions obtained by solving the equation
\begin{align*}
\int \lo 0 \hi 1 \hat \sigma \lo {X \hi i X \hi i } (s, t) \phi \lo r \hi i (s) d s= \lambda \lo r \hi i \phi \lo r \hi i (t), \quad r=1, \ldots, m \lo n,
\end{align*}
subject to the constraints $\langle  \phi \lo q \hi i ,  \phi \lo r \hi i  \rangle=0$, for $q \ne r, q, r=1, \ldots, m \lo n$ and $\|\phi \lo r \hi i \|=1$.  Then, the estimated scores $\hat \xi \hi i \lo {ur}$ are
given by  
  \begin{align*}
\hat \xi \lo {ur} \hi i = (\hat \lambda \lo r \hi i) \hi {-1/2} \langle  X \lo u \hi i, \hat \phi \lo r \hi i \rangle, \quad u = 1, \ldots, n, r=1, \ldots, m \lo n, i \in {\sf V}.
\end{align*}
%%%%%%%%%%%%%%%%%%%%%%%%%%
%  Algorithm
%%%%%%%%%%%%%%%%%%%%%%%%%%
For each $i \in \sf V$, let $\hat \xi \hi i \lo u=(\hat \xi \hi i \lo {u1}, \ldots, \hat \xi \hi i \lo {u m \lo n})^{\top}$ be the $m \lo n$-dimensional vector of the estimated scaled scores corresponding to the observation $X_u$, $u = 1, \ldots, n$.  Following  \cite{huang2010variable} we use B-spline functions to approximate the additive components $f \hi {ij} \lo {qr}$ in model \eqref{scorestrunc}.   To be precise, let $-1 = \tau \lo 0 < \tau \lo 1 < \ldots < \tau \lo {L \lo n}< \tau \lo {L \lo n+1}=1$ be an equidistant partition of the interval $[-1, 1]$ into $L \lo n +1$ subintervals $I \lo {b}=[\tau \lo b, \tau \lo {b+1}), b=0, \ldots, L \lo n-1$, and $I \lo {L \lo n}=[\tau \lo {L \lo n}, \tau \lo {L \lo n+1} ]$.  

  For the number of knots we make the following assumption.
Define $S \lo { \ell L \lo n}$ as the space of polynomial splines of degree $\ell \geq 1$ consisting of functions $s$ satisfying: 
 (i) the restriction of $s$ to the interval $I \lo {b}$ is a polynomial of degree $\ell$ for $1 \leq b \leq L \lo n;$ 
 (ii) for $\ell \geq 2$ and $1 \leq \ell'  \leq \ell-2$, $s$ is a $\ell'$ times continuously differentiable on the  interval $[-1,1]$.  
Then, there exists a basis of normalized B-splines functions $(h  \lo {k}) \lo { 1 \leq k \leq  k \lo n }$ for the space $S \lo {\ell L \lo n}$, where $k \lo n=L  \lo n+ \ell+1$, such that every  function $s \in S \lo {\ell L \lo n}$ 
can be represented as 
\begin{align*}
s (x)=\tsum \lo {k=1} \hi {k \lo n} \beta  \lo {k} h \lo {k}(x)   
 \end{align*} 
  (see \cite{schumaker2007spline}).
Under some smoothness conditions, the  additive functions $f \hi {ij} \lo {qr}$ can be represented by linear combinations of B-splines functions
\begin{align}\label{frepresentation}
f  \lo {qr} \hi {ij} (x) = \tsum \lo {k=1} \hi {\infty} \beta \hi {ij} \lo {qrk}   h  \lo {k}(x) \,  \quad q, r=1, 2, \ldots, %, (i,j) \in {\sf V} \times {\sf V}, i \ne j.
 \end{align} 
where  the sum of squared coefficients is summable, that is
\begin{align}\label{splinesasu1}
\sum \lo {k=1} \hi {\infty} ({\beta  \hi {ij} \lo {qrk}}) \hi 2 <\infty.  
\end{align}
By truncation of this series, we obtain the following approximation
\begin{align}\label{fapprox}
f  \lo {qr} \hi {ij} (x) \approx \tsum \lo {k=1} \hi {k \lo n} h  \lo {k}(x) \, \beta \hi {ij} \lo {qrk}, \quad q, r=1, 2, \ldots, %, (i,j) \in {\sf V} \times {\sf V}, i \ne j.
 \end{align}
 where the sequence $(k \lo n) \lo {n \in \mathbb{N}}$ diverges to infinity as $n \to \infty$ (note that this can always be achieved by increasing the number of knots in the partition).   
Hence, the corresponding function $f \lo {qr} \hi {ij}$ will be zero approximately if and only if  $\| \beta \hi {ij} \lo {qr} \| \hi 2 \lo 2=0$, where $\| \cdot \| \lo 2$ denotes the Euclidean norm of the $k \lo n$-dimensional vector  $\beta  \hi {ij} \lo {qr} =(\beta \hi {ij} \lo {qr1}, \ldots, \beta \hi {ij} \lo {qr k \lo n})^\top $, $q, r=1, \ldots, m \lo n$.   Thus, to encourage sparsity we propose to minimize the criterion
\begin{align*}
 {PL} \lo i (\beta, \hat \xi )=\frac{1}{2n} \tsum \lo {q=1} \hi {m \lo n} \tsum \lo {u=1} \hi { n}  \Big(\hat \xi \hi i \lo {ur} - \tsum \lo {j \ne i} \hi {p} \tsum \lo {r=1}  \hi {m \lo n} h \trans   (\hat \xi \hi j  \lo {ur}) \beta \lo {qr}  \hi {ij} \Big ) \hi 2 + {\lambda \lo {n}} \tsum \lo {j \ne i} \hi {p} \Big \{ \tsum \lo {q=1} \hi {m \lo {n}} \tsum \lo {r=1}  \hi {m \lo {n}} \| \beta \hi {ij} \lo {qr} \| \lo 2 \hi 2  \Big  \} \hi {1/2},
 \end{align*}
 subject to the constraint 
 \begin{align}\label{splinesconstraint}
\tsum \lo {u=1} \hi {n} h \trans   (\hat \xi \hi j  \lo {ur}) \beta \lo {qr}  \hi {ij} =0, \quad q, r =1, \ldots, m \lo n, j \in \sf V,
  \end{align}
 where $h \trans  (x) =(h  \lo {1} (x), \ldots, h  \lo { k \lo n}(x)) $ is the $ k \lo n$-dimensional vector of the B-splines basis functions and $\lambda \lo {n}$ is a tuning parameter.  The group lasso penalty $ \tsum \lo {j \ne i} \hi {p} \{ \tsum \lo {q=1} \hi {m \lo {n}} \tsum \lo {r=1}  \hi {m \lo {n}} \| \beta \hi {ij} \lo {qr} \| \lo 2 \hi 2  \} \hi {1/2}$ enforces all  regression coefficients $\beta  \hi {ij} \lo {qr1}, \ldots, \beta  \hi {ij} \lo {qr k \lo n}$ to either be all 0 or all nonzero, $q, r=1, \ldots, m \lo n$.
 Note that the centering constraint \eqref{splinesconstraint} accounts for the fact that the function $f \hi {ij} \lo {qr}$ in model \eqref{scorestrunc} satisfies $ E(f \hi {ij} \lo {qr}(\xi \hi j \lo r ))=0,$ for $q, r=1, \ldots, m \lo n$.% (see the approximation in equation \eqref{fapprox}).  

  This problem can be converted to an unconstrained optimisation problem by centering the basis functions. More precisely, defining
  {
 \begin{align}\label{centersplines}
\tilde h  \lo {nk}  (\hat \xi \hi j  \lo {ur})= h  \lo k (\hat \xi \hi j  \lo {ur}) - \frac{1}{n} \tsum \lo {u=1}  \hi {n}   h  \lo k  (\hat \xi \hi j  \lo {ur}), \; k=1, \ldots, k \lo n, r=1, \ldots, m \lo n, j \in \sf V,
  \end{align}}
 we consider the unconstrained optimization problem 
  \begin{align}\label{objective1}
 \widehat{{PL}} \lo i (\beta, \hat \xi )=\frac{1}{2n} \tsum \lo {q=1} \hi {m \lo n} \tsum \lo {u=1} \hi { n}  \Big(\hat \xi \hi i \lo {ur} - \tsum \lo {j \ne i} \hi {p} \tsum \lo {r=1}  \hi {m \lo n}  \tilde h  \lo {n} \trans  (\hat \xi \hi j  \lo {ur}) \beta \lo {qr}  \hi {ij} \Big) \hi 2 +   {\lambda \lo {n}} \tsum \lo {j \ne i} \hi {p} 
 \Big \{ \tsum \lo {q=1} \hi {m \lo {n}} \tsum \lo {r=1}  \hi {m \lo {n}} \| \beta \hi {ij} \lo {qr} \| \hi 2 \lo 2 \Big \} \hi {1/2},
 \end{align}
where
 \begin{align}\label{centersplines1}
  \tilde h  \lo {n}   (\hat \xi \hi j  \lo {ur})=(\tilde h \lo {n1}   (\hat \xi \hi j  \lo {ur}), \ldots, \tilde h   \lo {n k \lo n}   (\hat \xi \hi j  \lo {ur}))  \trans,
    \end{align}
 is the $ k \lo n$-dimensional vector  of the centered B-splines evaluated at the estimated scores.

Now  let $\hat \xi \hi i =(\hat \xi \hi i \lo {ur}) \lo {1 \leq u \leq n, 1 \leq r \leq m \lo n}$ be the $n \times m \lo n$  matrix of the estimated scores, and define
 \begin{align}\label{centersplinesmatrix}
 \mathbf{\tilde H \trans \lo n } (\hat \xi \hi {-i})=(\tilde H \lo n   (\hat \xi \hi {1}), \ldots, \tilde H \lo n  (\hat \xi \hi {i-1}), \tilde H \lo n   (\hat \xi \hi {i+1}), \ldots, \tilde H \lo n   (\hat \xi \hi {p})) \in \R \hi {n \times (p-1) k \lo n m \lo n }
     \end{align}
   as the vector of matrices $ \tilde H \lo n (\hat \xi \hi j)= ( \tilde h \lo {n} \trans (\hat \xi \hi j  \lo {ur})) \lo {1 \leq u \leq n, 1 \leq r \leq m \lo n} \in \R \hi {n \times k \lo n m \lo n}$.  Similarly, let 
    \begin{align*}
   B \hi i=(B   \hi {ij}, j \in \sf {V} \setminus \{i\} ) \in \R \hi {(p-1) k \lo n m \lo n   \times m \lo n},  
        \end{align*}
   be the vector of matrices $B \hi {ij}=({\beta \hi {ij} \lo {qr}} ) \lo {1 \leq q \leq m \lo n, 1 \leq r \leq m \lo n} \in \R \hi {k \lo n m \lo n \times  m \lo n}, j \ne i$. 
Then, following some algebraic manipulations, the objective function in \eqref{objective1} can be rewritten as
   \begin{align}\label{finalpenfunction}
 {\widehat{PL}} \lo i (B, \hat \xi )=\frac{1}{2n} \| \hat \xi \hi i -  \mathbf{\tilde H \lo n} \trans   (\hat \xi \hi {-i}) B  \hi {i}  \| \lo {F} \hi 2 +  \lambda \lo {n} \tsum \lo {j \ne i} \hi p \| B \hi {ij} \| \lo F,
 \end{align}
where $\| \cdot \| \lo {F}$ denotes the Frobenius norm. Finally, we define $\hat B \hi i \lo n$ as the solution of
  \begin{align*}
 \hat B \hi i \lo n =\mathrm{argmin} \{ \widehat {PL} \lo i (B  , \hat \xi  ): B  \in \R \hi {(p-1)k {\lo n m \lo n} \times m \lo n} \},
  \end{align*}
and propose to estimate the set $\sf E \lo n$ in \eqref{truncatededge} by 
\begin{align*}
\hat {\sf E} \lo n=\{ (i,j) \in {\sf V} \times {\sf V}: i \ne j, \| \hat B \hi {ij} \lo n \| \lo F >0 \; \mathrm{or} \; \| \hat B \hi {ji} \lo n \| \lo F >0 \}.
\end{align*}
We summarize the algorithm below

\begin{itemize}
\item[(1)]  Implement FPCA to obtain the estimated scores $ \hat \xi \hi i \lo  {ur}$ of each observation $X \hi {i} \lo u$ and then transform the scores into the range $[-1, 1]$ using a monotone transformation.  Choose $m \lo {n}$ such that   90\% of the total variation is explained. 
\item[(2)]    For a given $\lambda \lo n$ and for each $i \in \sf V$ solve the optimisation problem \eqref{finalpenfunction} using, for example, distance convex programming techniques, to find a sparse estimate of $B  {\hi {i}}$.
\item[(3)]   
Declare that there is an edge between node $i$ and node $j$ if and only if either $\| \hat  B  {\hi {ij} \lo n}\| \lo {F} \hi 2$ or $\| \hat  B   {\hi {ji} \lo n}\| \lo F \hi 2$ are not zero.
\end{itemize}

\section{Statistical guarantees} \label{theory}
\label{sec4}
\def\theequation{4.\arabic{equation}}
\setcounter{equation}{0}

In this section we study the theoretical properties of the proposed estimator of the graph structure of the AFGM where we allow  the number of nodes $p$ to diverge to infinity with increasing sample size.  A particular technical challenge in deriving the asymptotic theory consists in the fact that the additive structure is applied to the unobserved variables $\xi \hi i  \lo {ur}$, and the estimator $\hat B  \hi i \lo n$ obtained from minimizing  \eqref{finalpenfunction} is based on the estimated scores. Thus, the  error in these estimated coefficients must be taken into account for the analysis of the procedure.

We begin by introducing some notation. For any two positive sequences of real numbers $(a \lo n) \lo {n \in \mathbb{N}}$ and $(b \lo n) \lo {n \in \mathbb{N}}$, we write $a \lo n \lesssim  b \lo n$ if $a \lo n  \leq K \lo 1 b \lo n$ for some  constant $0 < K \lo 1 < \infty$ which does not depend on $n$.  We use the notation $a \lo n \asymp b \lo n$ representing the property $A \leq \mathrm{inf} \lo n | \frac{a \lo n}{ b \lo n} | \leq \mathrm{sup} \lo n | \frac{a \lo n}{ b \lo n} | \leq B$, for  positive constants $A$ and $B$.  Moreover,  given a matrix $A=(a \lo {ij}) \lo {1 \leq i \leq M \lo 1,1 \leq j \leq M \lo 2} \in \R \hi { M \lo 1 \times  M \lo 2}$, we use $\|A \| \lo F$ for  the Frobenious norm and $\| A \|  \lo {2}$ for  the operator norm. Finally, for any two symmetric matrices $A$ and $B$, we use the notation $  A \preceq B$ to denote the property that the matrix $B-A$ is nonnegative definite. 

Let 
 \begin{align}\label{boptimuntr}
 B  \hi {*i} \lo {m \lo n} =(B \hi {*ij}  \lo {m \lo n}, j \in \sf {V} \setminus \{i\}),
 \end{align}
with $B \hi {*ij}  \lo {m \lo n} =\{\beta \hi {*ij} \lo {qrk}: 1 \leq q, r \leq m \lo n, k \in \mathbb{N} \}$ be the true population matrix of parameters for the optimal prediction, defined by 
\begin{align*}
B \hi {*i}  \lo {m \lo n}  =\argmin_{\beta \hi {ij} \lo {qrk}}    \sum \lo {q=1} \hi {m \lo n}   E \Big  [ \xi \hi i \lo {q} - \sum \lo {j \ne i} \hi {p} \sum \lo {r=1} \hi {m \lo n} \sum \lo {k=1} \hi {\infty} ( h \lo {k}( \xi \hi j \lo {r}) -{E( h \lo {k}( \xi \hi j \lo {r}))} ) \, \beta \hi {ij} \lo {qrk}  \Big ] \hi 2,
 \end{align*}
 where $h \lo {k} (\cdot) \lo {k \geq 1}$ are the  B-splines functions used in the representation \eqref{frepresentation}. We define the truncated neighbourhood $\sf N \hi i \lo n$ of each node $i \in \sf V$ by
  \begin{align*}
{\sf N \hi i \lo n}=\{ j \in {\sf V} \setminus \{i \}:  \| B {\hi {*ij} \lo {m \lo n}}   \| \lo F > 0 \} 
\end{align*}
(note that $\| B {\hi {*ij} \lo {m \lo n}}   \| \lo F < \infty$ by assumption \eqref{splinesasu1}). Using this  representation and observing the expansion \eqref{frepresentation} of $f  \lo {qr} \hi {ij}$, the edge set $\sf E \lo n$ defined in \eqref{truncatededge} can be rewritten as
\begin{align}\label{truegraph2}
{\sf E \lo n}= \{ (i, j) \in {\sf V} \times {\sf V}: i \ne j, i \in {\sf N \hi j \lo n} \; \mathrm{or}   \;  j  \in { \sf N \hi i \lo n} \}.  
  \end{align}
Let 
$$
f  \lo {qr} \hi {ij} ( \xi \hi j \lo {r}) = \tsum \lo {k=1} \hi {\infty} \beta \hi {*ij} \lo {qrk} h \lo {k} ( \xi \hi j \lo {r}) = \tsum \lo {k=1} \hi {\infty}  \beta \hi {*ij} \lo {qrk} \tilde  h \lo {k} ( \xi \hi j \lo {r}) 
$$
(for the second equality we use the fact that  $E[ f  \lo {qr} \hi {*ij} ( \xi \hi j \lo {r})] =0$), where the functions $ \tilde{h} \lo k$ are defined by
 \begin{align}\label{centertruesplines}
  \tilde{h} \lo k ( \xi \hi j  \lo {r})= h \lo {k}( \xi \hi j \lo {r}) -{E( h \lo {k}( \xi \hi j \lo {r}))}. 
 \end{align}
We obtain  from \eqref{scorestrunc} the representation
% an anlogue  of the representation \eqref{scorestrunc}, that is 
\begin{align}\label{optpredictor}
 \xi \hi {i} \lo {q} =  \tsum \lo { j \in {\sf N \hi i \lo n}} \tsum \lo {r =1} \hi {m \lo n} f  \lo {qr} \hi {ij} ( \xi \hi j \lo {r}) + \epsilon \lo q \hi {i}, \quad q=1, \ldots, m \lo n, i= 1, \ldots, p, 
\end{align}
where 
$\epsilon \lo q = \xi \hi {i} \lo {q} -  E[ \xi\hi {i}
\lo q  |  \{ \xi \hi j \lo {r}, j \ne i; r=1, \ldots , m_n \} ] $. 
Thus, the best predictor of $\xi \hi i \lo q $ is an additive function of the scores in the set of neighbours $\sf N \hi i \lo n$ of the node $i$ only.  
  
Let $  \tilde{h}  ( \xi \hi j  \lo {r})=( \tilde{h} \lo {1}   ( \xi \hi j  \lo {r}), \ldots,  \tilde{h}   \lo { k \lo n}   ( \xi \hi j  \lo {r})) \trans  \in \R \hi {k \lo n}$ be the  vector of the centered $k \lo n$ B-splines evaluated at the unobserved scaled scores $\xi \hi j \lo {r}$, $r=1, \ldots, m \lo n$, $j \in \sf V$, {  and define the $ 1 \times k \lo n m \lo n $  and 
$ 1 \times n \hi i k \lo n m \lo n $ vectors
    \begin{align} \nonumber
   \tilde{H} (\xi \hi j)  &=(  \tilde{h}\trans    ( \xi \hi j  \lo {r}) ) \lo {1 \leq r \leq m \lo n}  %\in \R \hi{k \lo n m \lo n}
   , \\
\mathbf{ \tilde{H}\trans (\xi \hi { \sf N \hi i \lo n} ) }& = ( \tilde{H} (\xi \hi j)  , j \in \sf N  \hi i \lo n)
%\in \R \hi{n \hi i k \lo n m \lo n}
,  
\label{hr1}
    \end{align}
   where $1 \leq  n \hi i \leq p $ is the cardinality of the set  $\sf N \hi i \lo n$.  }
   Finally, we introduce the matrices  
{{ \begin{align} \label{truemat1}
 \Sigma \hi {*} \lo {\sf N \lo n \hi i \sf N \lo n \hi i} = E \Big ( \mathbf {\tilde{H} (\xi \hi { \sf N \hi i \lo n}} ) \mathbf{ \tilde{H}  \trans (\xi \hi { \sf N \hi i \lo n} ) } \Big  ) \in \R \hi {n \hi i k \lo n m \lo n  \times n \hi i k \lo n m \lo n}
\end{align}
 and 
    \begin{align}     \label{truemat2}
   \Sigma \hi * \lo {{\xi \hi j}  \sf N \lo n \hi i  }
 = E \Big (\tilde{H} \trans ( \xi \hi { j}) \mathbf{ \tilde{H} 
  \trans (\xi \hi { \sf N \hi i \lo n} )  }  \Big )   ~ \in \R \hi { k \lo n m \lo n  \times n \hi i k \lo n m \lo n}.
\end{align}}}
            Let 
   \begin{align*}
 B  \hi {*i} \lo { n } =(B \hi {*ij} \lo {m_nk_n} , j \in \sf V \setminus \{i \})  ~ \in \R \hi  {(p-1) k \lo n m \lo n \times m \lo n}
   \end{align*}
 denote the truncated  version of the population matrix defined in \eqref{boptimuntr}, where
    \begin{align}\label{optimbeta}
 B \hi {*ij}  \lo {m \lo n k \lo n}= \{ \beta \lo {qrk} \hi {*ij}: 1 \leq q, r \leq m \lo n, 1 \leq k \leq k \lo n \},
    \end{align}
 and define
\begin{align}
   \label{hd2} 
  B \hi {* \sf N  \hi i \lo n} \lo n=(B \hi {*ij} \lo {m \lo n k \lo n}, j\in \sf N  \hi i \lo n) \in \mathbb{R} \hi {n \hi i k \lo n m \lo n \times m \lo n}.
    \end{align}
 Recalling that $\hat B  \hi {i} \lo n=(\hat B \hi {ij}, j \in \sf V \setminus \{ i \} ) \in \R \hi {(p-1)k \lo n m \lo n \times m \lo n}$ is the solution to the minimization problem \eqref{finalpenfunction}.   The estimated neighbourhood for each node $i \in \sf V$ is  defined
 \begin{align*}
   {\sf  \hat  N \hi i \lo n}=\{ j \in {\sf V} \setminus \{i \}:  \| \hat B {\hi {ij} \lo { n}}   \| \lo F > 0 \},
   \end{align*}
which yields an alternative representation of the estimated edge set
\begin{align}\label{estimgraph2}
\hat {\sf E \lo n}=\{ (i, j) \in \sf V \times \sf V: i \in \sf  \hat  N \hi j \lo n \text{ or }   j \in \sf  \hat  N \hi i \lo n\}.
\end{align}
 
 For the statement of our theoretical results we require several assumptions. Assumption \ref{as1} is a similar assumption as made by \cite{qiao2018functional} and refers to the eigensystem of the covariance operator defined in \eqref{covop}.

 \begin{assumption}. \label{as1} 
{\rm ~ \\ (i) There exist  positive constants $d \lo 0, d \lo 1$ and $d \lo 2$ such that
\begin{align*}
d \lo 0  r \hi {-\beta} \leq \lambda  \hi i \lo r \leq  d \lo 1 r \hi {-\beta},  \quad   \lambda  \hi i \lo r -  \lambda  \hi i \lo {r+1} \geq d \lo 2 \inv  r \hi {-1-\beta}  \quad \text{for } r \geq 1,
\end{align*}
and for some $\beta >1$}.
{\rm ~ \\ (ii) The number of principal component scores $m \lo n$ satisfies $m \lo n \asymp n \hi \alpha$ for some constant $\alpha \in [0, \frac{1}{2+3\beta})$.}
{\rm ~ \\ (iii) The eigenfunctions $\phi \hi i \lo r$   of the covariance operator defined in \eqref{trueeigenpairs} 
  are continuous and satisfy 
  \begin{align*}
  \max \lo {j \in \sf V} \sup \lo {s \in [0,1]} \sup \lo {r \in \mathbb{N}} |\phi \hi j \lo r (s)| \leq C < \infty.
  \end{align*}
}
\end{assumption}
The next two conditions refer to the smoothness of the  functions $f \hi {ij} \lo {qr}$ in model \eqref{afm}.  To be precise, let $\kappa$ be a nonnegative integer and let $\rho \in (0, 1]$. We define $\sten {F} \lo  {{\kappa,\rho}}$ as 
 the H{\"o}lder space of functions $f:[0,1] \to \R $ whose $\kappa$th derivative exists and satisfies a Lipschitz condition of order $\rho$, and additionally  
 satisfy the condition 
 \begin{align}\label{fbounded}
  \| f  \| \lo {\infty}=\mathrm{sup} \lo {x \in [0,1]}  | f (x)| \leq F
  \end{align}
   for some $F>0$. 
\begin{assumption} Let $d=\kappa+ \rho >0.5$ and assume \label{as3} $f \hi {ij} \lo {qr} \in \sten {F}_{{\kappa,\rho}}$ and $E [f \hi {ij} \lo {qr}(\xi \hi j \lo {ur})]=0$, for all   $q, r=1, \ldots, m \lo n$ and $(i, j) \in \sf V \times \sf V$.
\end{assumption}
{{\begin{assumption}  \label{as4} The joint density function, say $p \hi j$, of  the random vector $\xi \hi j=(\xi \hi j \lo {1}, \ldots, \xi \hi j \lo {m \lo n}) \trans$ is bounded away from zero and infinity on $[0,1] \hi {m \lo n}$ for every $j=1, \ldots, p.$
\end{assumption}}}
In order  to derive graph estimation consistency,  we make the following assumption about the errors $\epsilon^i_1,\ldots,\epsilon^i_{m_n}$ in model \eqref{scorestrunc}. A similar condition was also postulated by  \cite{voorman2013graph} for joint additive models in the multivariate setting.
\begin{assumption}\label{aserrors}  There exists a constant $C>0$ such that 
$P(|\epsilon \hi i \lo {q}| > x) \leq 2 \exp (- C x \hi 2)$ for all $x \geq 0$ and $q=1, \ldots, m \lo n, i \in \sf V$.
\end{assumption}

\begin{assumption} \label{assnu}
$k \lo n=O(n \hi {\nu})$ for some $\nu >0$, where $\nu \leq \frac{\alpha(2+3\beta)}{2d-4}$ if  $d \geq 2$.
\end{assumption}

% We consider the following sparsity assumption.

\begin{assumption}\label{sparsity1} (Sparsity) There exists a constant $\theta>0$ such that for all $i \in \sf V$
\begin{align*}
\tsum \lo {j \in  \sf N \hi i \lo n} \|  B \hi {*ij}  \lo {m \lo n k \lo n}  \| \lo F < \theta.
\end{align*}
\end{assumption}

\begin{assumption} \label{boundedeigen} (Bounded eigenspectrum) \label{as5} The minimum eigenvalue $\Lambda \lo {\min} (\Sigma \hi * \lo { \sf N \lo n \hi i \sf N \lo n \hi i})$  of the matrix $\Sigma \hi * \lo { \sf N \hi i \sf N \hi i}$ defined in \eqref{truemat1} satisfies
\begin{align*}
\Lambda \lo {min}(\Sigma \hi * \lo { \sf N \lo n \hi i \sf N \lo n \hi i}) >C \lo {\min} %\leq \Lambda \lo {\max} (\Sigma \hi * \lo { \sf N \lo n \hi i \sf N \lo n \hi i}) \leq C \lo {\max} < \infty.
\end{align*}
for some constant $C \lo {\min} >0$.
\end{assumption}
{{\begin{assumption} \label{irrepr} (Irrepresentable condition) There exists a constant $0 < \eta \leq 1$  such that 
\begin{align}\label{irreprequation}
\max \lo { j \notin \sf N \hi i \lo n } \| \Sigma \hi * \lo {\xi \hi j \sf N \hi i \lo n} (\Sigma \hi * \lo { \sf N \hi i \lo n \sf N \hi i \lo n}) \hi {-1} \| \lo F \leq \frac{ 1-\eta}{\sqrt{n \hi i}}.
\end{align}
\end{assumption}}}
Assumption \ref{boundedeigen} states that the minimum eigenvalue of the population matrix $\Sigma \lo { \sf N \lo n \hi i \sf N \lo n \hi i}$ is bounded away from 0.   Assumption \ref{irrepr} is the classical irrepresentable condition which is necessary and sufficient to show  model selection consistency of the group lasso \citep{ yuan2006model, bach2008consistency}.  According to \cite{meinshausen2009lasso} if the irrepresentable condition is relaxed, the lasso selects the correct non-zero coefficients but it may select some additional zero components. \cite{ravikumar2009sparse} and  \cite{obozinski2011support} considered similar assumptions for sparse additive and for  high-dimensional models, respectively. We now state our main  theoretical result for  the  estimator $ { \hat {\sf N }\hi i \lo n}$  of the neighbourhood  corresponding to the node $i \in \sf V$.
\def\cip{\stackrel{\nano P}{\rightarrow}}
%%%%%%%%%%%%
%main theorem
%%%%%%%%%%%%%%%%%
\begin{theorem}\label{theorem2section4}  Suppose that Assumptions \ref{as1} - \ref{irrepr} are satisfied and the regularization parameter $\lambda \lo n$ satisfies for all $i$
 %\begin{align}\label{as1theorem2}
%n \hi {1-\alpha(2+3\beta)}  > C [  (n \hi i) \hi {5/2} m \lo n \hi 2 k \lo %n \hi 4  \log (n \hi i m \lo n k \lo n)],
 % \end{align}
%for some  large constant $C>0$, and the regularization parameter $\lambda \lo n$ satisfies for all $i$
   \begin{align}\label{as2theorem2}
 \frac{n \hi i m \lo n \hi {3/2}}{ k \lo n \hi {d}  \tsum \lo {j \in  \sf N \hi i \lo n} \|  B \hi {*ij}  \lo {m \lo n k \lo n}  \| \lo F}   \lesssim  \lambda \lo n \lesssim   (n \hi i) \hi {-3/2} (b \hi {*i} \lo n) \hi 3 (\tsum \lo {j \in  \sf N \hi i \lo n} \|  B \hi {*ij}  \lo {m \lo n k \lo n}  \| \lo F) \hi {-2},
  \end{align}
 where $b \hi {*i} \lo n= \min \lo {j \in \sf N \hi i \lo n} \| B \hi {* ij }  \lo {m \lo n k \lo n}  \| \lo F$.  Then, 
   \begin{equation*}
\begin{split}
 P \Big( \hat {\sf N} \hi i \lo n \ne {\sf N }\hi i \lo n  \Big) \lesssim   \exp \Big ( - C \lo 1 \frac{n \hi {1- \alpha (2+3\beta)} (\lambda \lo n \tsum \lo {j \in  \sf N \hi i \lo n } \|  B \hi {*ij} \lo {m \lo n k \lo n} \| \lo F) \hi 2  }{n \hi i m \lo n \hi 2  k \lo n \hi {4} }+2 \log (p m \lo n k \lo n) \Big ),
\end{split}
  \end{equation*}
where $C \lo 1>0$.
\end{theorem}
The proof of Theorem \ref{theorem2section4} is complicated and given in the Appendix.
%and is divided into two parts.  First, in Section \ref{sec82} we prove that the sample versions of the matrices \eqref{truemat1} and \eqref{truemat2} satisfy the Assumptions \ref{boundedeigen} and \ref{irrepr} with high probability. Then, in Section \ref{sec83} we prove a conditional result (stated as Proposition \ref{sampletheorem} in the Appendix), for the sample design matrices.  
A major difficulty consists in the fact that
% In addition to this, since 
the objective function \eqref{finalpenfunction} is based on the estimated scores  and 
%the extra difficulty in developing the proof  is to 
one has to establish concentration bounds in the estimation of the sample design matrix $ \Sigma \hi n \lo {\sf N \hi i \lo n \sf N \hi i \lo n}$ using the estimated scores, rather than the true scores (stated as Theorem \ref{theorem1} in the Appendix).

Recalling the representations \eqref{truegraph2} and \eqref{estimgraph2} for the edge set $\sf E \lo n$ and its estimate $\hat {\sf E} \lo n$ respectively.   Using the union bound of probability and Theorem \ref{theorem2section4} we obtain the following result. 

{ \begin{corollary} 
If the assumptions of Theorem \ref{theorem2section4}
are satisfied, we have for a positive constant $C_1 >0$
$$
P(\hat {\sf E} \lo n \ne {\sf E} \lo n) 
\lesssim  
\exp \Big ( - C \lo 1 \frac{n \hi {1- \alpha (2+3\beta)} (\lambda \lo n 
\min_{i=1}^p \tsum \lo {j \in  \sf N \hi i \lo n } \|  B \hi {*ij} \lo {m \lo n k \lo n} \| \lo F) \hi 2  }{p   m \lo n \hi 2  k \lo n \hi {4} }+2 \log (p   m \lo n k \lo n) \Big )
. 
$$
\end{corollary}
}

\begin{remark}
{\rm ~~

\noindent
{(a) 
Under the assumptions of Theorem \ref{theorem2section4}
it follows that 
 \begin{align*}
 P(\hat {\sf E} \lo n \ne {\sf E} \lo n) \to 0,
 \end{align*}
if $n \to \infty,p \to \infty$ and
 \begin{align}\label{corconditions}
\frac{n \hi {1- \alpha (2+3\beta)} (\lambda \lo n \min_{i=1}^p b \lo n \hi {*i})}{ n \hi i m \lo n \hi 2 k \lo n \hi 4} \to \infty ~ \text{ and } ~~
\frac{n \hi i m \lo n \hi 2 k \lo n \hi 4 \log (p m \lo n k \lo n)}{\lambda \lo n \hi 2} = o \big (  \min_{i=1}^p b \lo n \hi {*i} \big ).
 \end{align}
 For example, if $m \lo n=O(n \hi \alpha)$ with $\alpha \in [0, \frac{1}{2+3 \beta})$, $k \lo n=O(n \hi \nu)$
 and $\max \lo {i \in \sf V} n \hi i=O(n \hi {\theta})$ with  $0 \leq \theta < 1$, then, \eqref{corconditions} reduces to
\begin{align*}
\frac{\log(p m \lo n k \lo n)}{n \hi {1-(\alpha(4+3\beta)+\theta+4 \nu)} \lambda \lo n}=o \big (\min_{i=1}^p b \lo n \hi {*i} \big).
\end{align*}
(b)
In the case of scalar data ($\alpha=0$), the conditions of  Theorem \ref{theorem2section4} will be implied by
\begin{align*}
\frac{n \hi i}{k \lo n \hi d} \lesssim \lambda \lo n b \lo n \hi {*i}, \quad \lambda \lo n (n \hi i) \hi {3/2} \lesssim (b \lo n \hi {*i}) \hi {3}
~\mbox{ and }~~
\sqrt{\frac{n \hi i k \lo n \hi 4 \log({p } k \lo n)}{n\lambda \lo n \hi 2}}=o(b \lo n \hi {*i})~,
\end{align*}
which are similar to the assumptions made in \cite{ravikumar2009sparse} and  
\cite{voorman2013graph}  for the analysis of scalar data by sparse additive models.
}}
\end{remark}

%%%%%%%Numerical examples%%%%%%%%%%%%%%

\section{Finite sample properties} \label{simulations}
\def\theequation{5.\arabic{equation}}
\setcounter{equation}{0}

\subsection{Simulated data}  \label{sec51}
In this section we investigate the finite sample performance of the proposed model (AFGM) by means of a simulation study. We also compare the new methodology with  the functional additive precision operator (FAPO) of \cite{li2018nonparametric} and the FGGM of  \cite{qiao2018functional}, where  we
 consider two scenarios: nonlinear dependence and linear dependence.

Given an edge set $\sf E$ of a directed acyclic graph, 
we generate functional data by the model
\begin{align}\label{simulatefunctions}
X \lo u \hi i (t \lo s )=\tsum \lo {(i, j ) \in \sf E} \tsum \lo {q=1}  \hi {5} \tsum \lo {r=1} \hi {5} f \hi {ij} \lo {qr}( \xi \hi {j} \lo {ur}) \phi  \lo q (t \lo s) \, + \epsilon \hi i \lo {us} , \quad u=1, \ldots, n, 
%s=1, \ldots, 100~,
 \end{align} 
 where $ \phi \hi i \lo 1 , \ldots, \phi \hi i \lo 5  $ are the first $5$ functions of the orthonormal Fourier basis, and  the errors $\epsilon \hi i \lo {us}$ form an i.i.d. sample from a $\mathcal N(0, 0.5 \hi 2)$ distribution.  
 In all simulation experiments in this section, this data is smoothed to obtain continuous  functions $X \lo u \hi i $ using $10$ B-spline basis functions of order $4$; that is, piecewise polynomials of degree $3$.
 As a consequence the scores satisfy a  structural equation of the form
\begin{align}\label{simulatescores}
\xi \lo {uq} \hi i=\tsum \lo {(i, j ) \in \sf E}  \tsum \lo {r=1} \hi {5} f \hi {ij} \lo {qr}( \xi \hi {j} \lo {ur}) + \tilde \epsilon \hi i \lo {uq}, \quad u=1, \ldots, n, q=1, \ldots, 5
 \end{align} 
 \citep[see][]{pearl2002causality},
where the errors $ \tilde \epsilon \hi i \lo {uq}$ form an i.i.d. sample a centred  normal distribution.
For simplicity we assume  $f \hi {ij} \lo {qr}( x)=f (x)$ for all $q, r=1, \ldots $ and for all $(i, j) \in E$. In all examples,  we center $ f ( \xi \hi {j} \lo {ur})$ to have 0 mean, and we generated $n=100$ functions observed at 100 equally spaced time points $0=t \lo 1, \ldots, t \lo {100}=1$. 

We consider directed acyclic graphs with $p=100$ nodes so that  1\% of pairs of vertices are randomly selected as edges.  Then, we moralized the directed graph in order to obtain the undirected graph.     We choose $m \lo n=5$  functional principal components scores  so that at least 90\% of the total variation is explained.  Furthermore, we approximate each additive function using B-splines of order $4$. For simplicity we choose the same spline functions for all $j = 1,\ldots,p$ and for all $r=1, \ldots, m \lo n$. For the choice of $k \lo n$,  we follow \cite{meier2009high} and take $k \lo n = 4+\lceil \sqrt{n} \rceil$. 
 
 For each scenario, we produce the average ROC curves (over 50 replications) for a range of 50 tuning parameters  for the 3 functional graphical models estimators.  To draw the curves, we compute for different regularization parameters $\lambda$ the  positive rate (sensitivity) and false positive rate (1-specificity) which are defined as
 \begin{align*}
 \mathrm{TP}= \frac{\tsum \lo{1 \leq j < i \leq p} I \{(i, j) \in \sf E \lo n, (i, j ) \in \hat E \lo n \}}{\tsum \lo{1 \leq j < i \leq p} I \{(i, j) \in \sf E \lo n \}},  \,
  \mathrm{FP}= \frac{\tsum \lo{1 \leq j < i \leq p} I \{(i, j) \notin \sf E \lo n, (i, j ) \in \hat E \lo n \}}{\tsum \lo{1 \leq j < i \leq p} I \{(i, j) \notin \sf E \lo n \}}. 
 \end{align*}

  \subsubsection{Scenario 1: nonlinear models}
 We use the following nonlinear models, where the linearity assumption \eqref{flm} does not hold.
 %  \paragraph{ Example 1:} 
   We first consider the following model, Model I, used in  \cite{zhu2014structured}
 \begin{align*}
\begin{split}
\mbox{  Model I: } \qquad
 f ( x) &=1.4+ 3x-\frac{1}{2}+\sin (2 \pi (x-\frac{1}{2}))+ 8 (x-\frac{1}{3})\hi 2-\frac{8}{9}. \\
 \end{split}
   \end{align*}
 For the choice of scores in \eqref{simulatescores}, we simulate $\xi \hi {i} \lo {ur}$ independently from the uniform distribution $U [-1,1]$ for all $r=1, \ldots, m \lo n, i \in {\sf V}, u=1, \ldots, n$.  Furthermore, the errors $ \epsilon \hi i \lo {uq}$ in \eqref{simulatescores} form an i.i.d. sample from the normal distribution with variance $0.1$, that is $\mathcal N(0, 0.5 \hi 2)$.
 
%  \paragraph{ Example 2:}
  The second example was considered in \cite{meier2009high}
 \begin{align*}
\begin{split}
\mbox{  Model II: } \qquad
 f ( x) &=-\sin (2x) + x \hi 2-25/12+x+\exp (-x)-2/5 \cdot \mathrm {sinh} (5/2). \\
 \end{split}
   \end{align*}
 The scores $\xi \hi {i} \lo {ur}$ were simulated independently from the uniform distribution $U [-2.5,2.5]$ for all $r=1, \ldots, m \lo n, i \in {\sf V}, u=1, \ldots, n$, and  the errors $ \epsilon \hi i \lo {uq}$ in \eqref{simulatescores} were simulated from the normal distribution  $\mathcal N(0, 1)$.

 The left and middle panel of Figure 1 show the averaged ROC curves over 40 replications corresponding to  the two models.   In first two lines Table 1, we report the means and standard deviations (in parentheses) of the associated area-under-curve (AUC) values.  An AUC close to 1 means a better performance for the estimator.
 We observe from the plots  in Figure 1 and from Table 1, that for the AFGM estimator the areas under the ROC are substantially larger than for the FGGM, indicating that our new method AFGM dominates the FGGM. Similarly, the AFGM performs better than FAPO, indicating the benefit of a sparse and high-dimensional scheme.
  \begin{center}
 \includegraphics[scale=0.3]{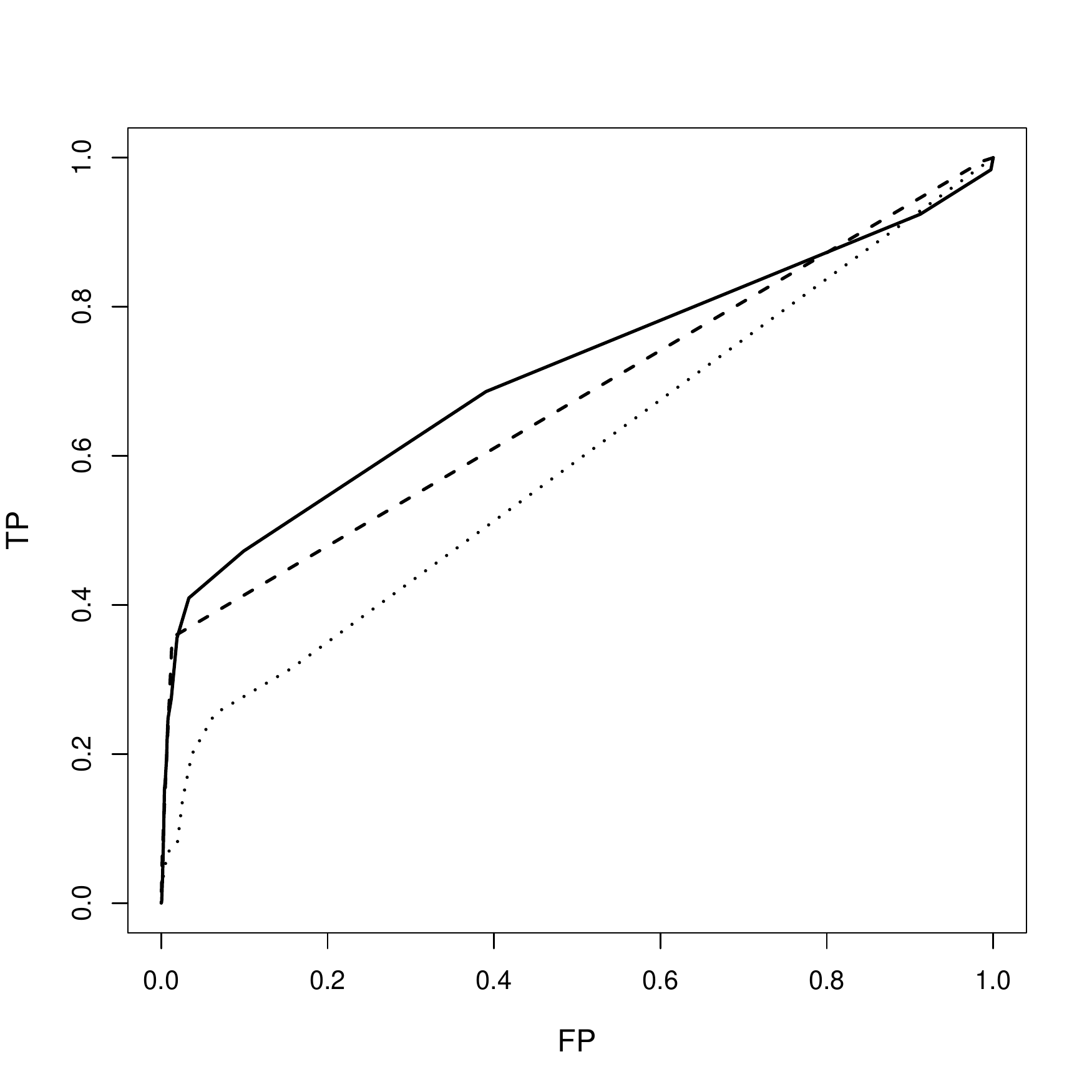} 
  \includegraphics[scale=0.3]{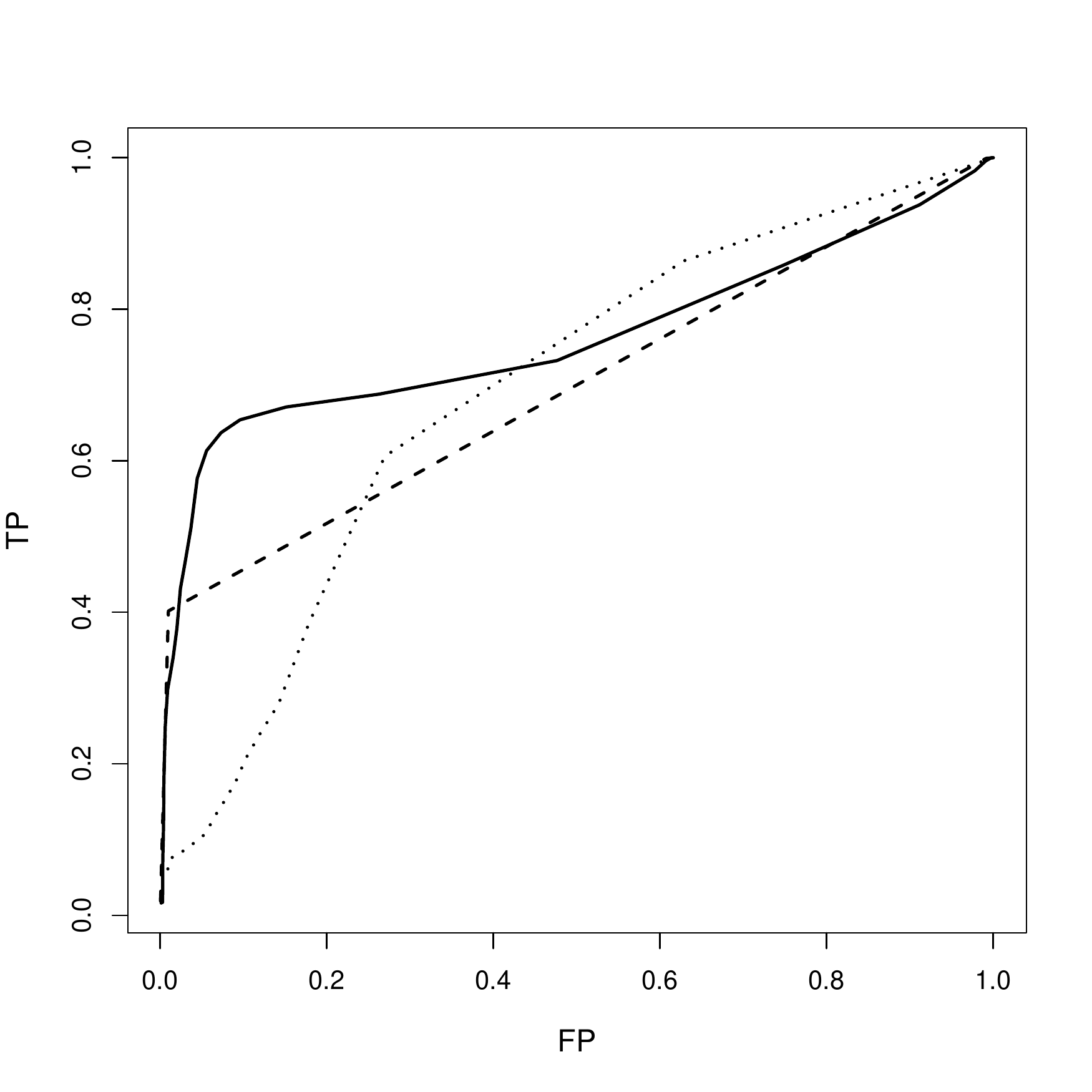} 
  \includegraphics[scale=0.3]{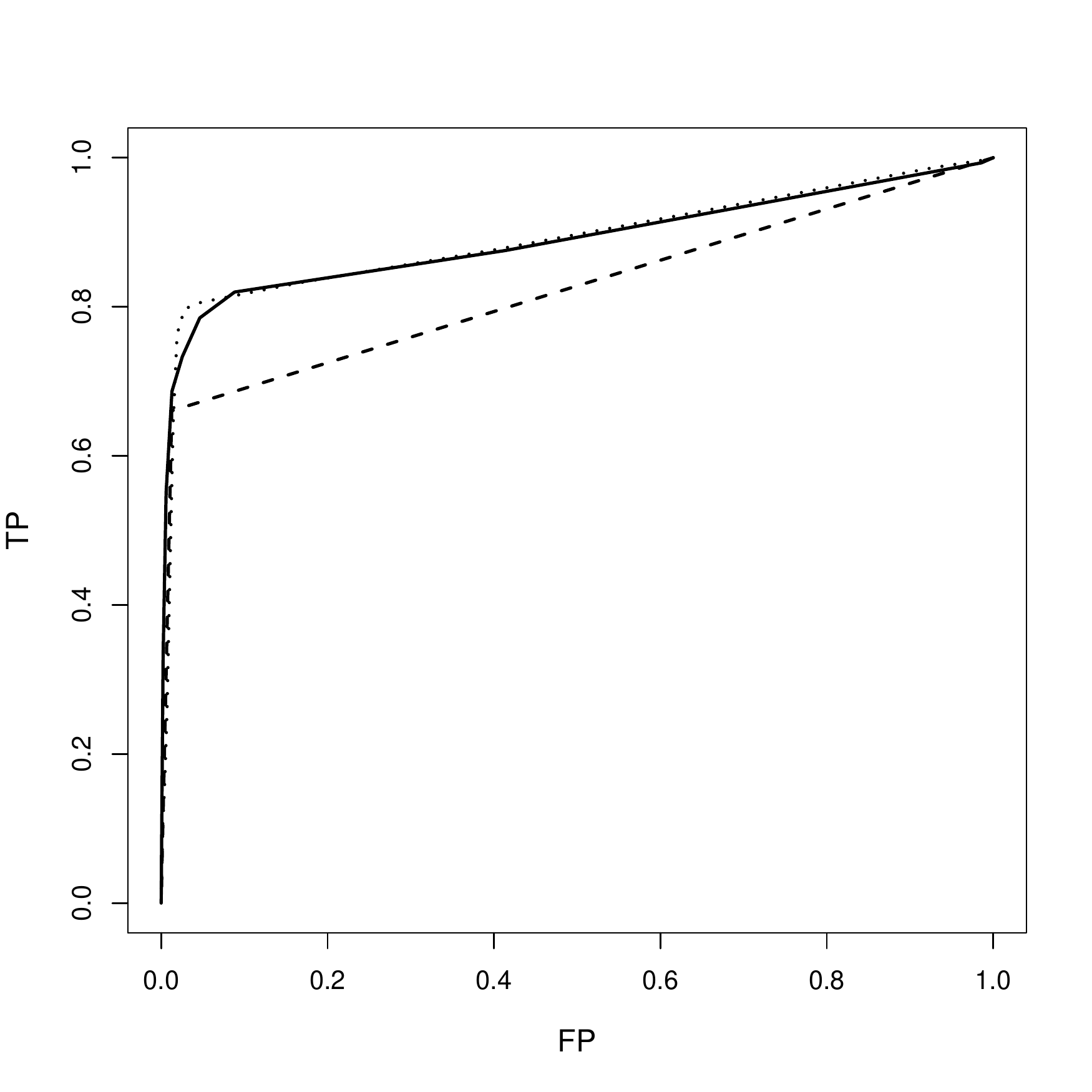}
   \begin{center}
{\small{\parbox{5in}{\raisebox{9pt}{Figure 1.} \ \, \parbox{4in}{ROC curves ((AFGM ($-$),  FAPO ($- - -$), FGGM ($\cdot \cdot \cdot$)) \\ for Model I (left) and Model II (middle) and Model III (right).}}}}
\end{center}
\end{center}

\bigskip

{\small
\begin{center}
\begin{tabular}{|c|c |c c c |}\hline\hline

\multirow{2}{*}{$p$}&\multirow{2}{*}{\scriptsize Models}		& \multicolumn{3}{|c|}{\scriptsize Methods}\\ \cline{3-5}
&	 			&  {\scriptsize  AFGM}	 &	  {\scriptsize  FAPO}	          & 	   {\scriptsize  FGGM}		         \\ \hline
\multirow{3}{*}{$100$} & I    & 0.73 (0.02)   & 0.67 (0.02)   &0.59 (0.01)     \\
 &{II}                &  0.76  (0.01) & 0.70 (0.02)     &  0.69 (0.02)         \\
& {III}
                    & 0.89 (0.01)   &  0.82 (0.01) &   0.90 (0.01)  \\
\hline\hline     \end{tabular}
\medskip \\
{ {Table 1.  Means and standard errors (in parentheses) for AUC \\for models I and II.}} 
\end{center}
}

\medskip

  \subsubsection{Scenario 2: linear model.}
 Next, we consider a model, where the linearity assumption is satisfied, to see how much efficiency might be lost by employing a nonparametric model under the Gaussian assumption. The model is generated by  \eqref{simulatescores} and \eqref{simulatefunctions} with
   \begin{align*}
\mbox{  Model III: } \qquad
 f ( x) =x  
 \end{align*}
 The scores $\xi \hi {i} \lo {ur}$ were simulated independently from the standard Gaussian distribution.  To implement the AFGM we truncate  the scores such that they are located in the interval $[-1,1]$.  The right panel in Figure 1 presents the averaged ROC curves for Model III.
The lower part of Table 1 reports the means and standard deviations of AUC. We can see that under the linearity assumption the AFGM is comparable with the FGGM and both methods show an improvement compared to FAPO.
%    \begin{center}
%\includegraphics[scale=0.3]{ROC-FAGMvsFGLASSO-AND-LINEAR-Gaussian-scoresn-100-p=100.pdf}
% \begin{center}
%{\small{\parbox{5in}{\raisebox{8pt}{Figure 2.} \ \, \parbox{4in}{ROC curves ((AFGM ($-$),  FAPO ($- - -$), FGGM ($\cdot \cdot \cdot$)) \\ for Model III. }}}}
%\end{center}
%\end{center}
%\bigskip
%{\small
%\begin{center}
%{ {Table 2.  Means and standard errors (in parentheses) for AUC \\for model III.}}\\
%\medskip
%\begin{tabular}{|c |c c c |}\hline\hline
%\multirow{3}{*}{$p$}		& \multicolumn{3}{|c|}{\scriptsize Methods}\\ %\cline{2-4}
%&	 			  {\scriptsize  AFGM}	 &	  {\scriptsize  FAPO}	          & 	   {\scriptsize  FGGM}		         \\ \hline
%\multirow{1}{*}{$50$}      &  0.89 (0.007)    & 0.91 (0.009)   & 0.98 (0.02)    \\
%\hline
%\multirow{1}{*}{$100$} 
%                                           & 0.89 (0.01)   &  0.82 (0.01) &  % 0.90 (0.01)  \\
%\hline\hline     \end{tabular}
%}
 \begin{center}
 \includegraphics[scale=0.4]{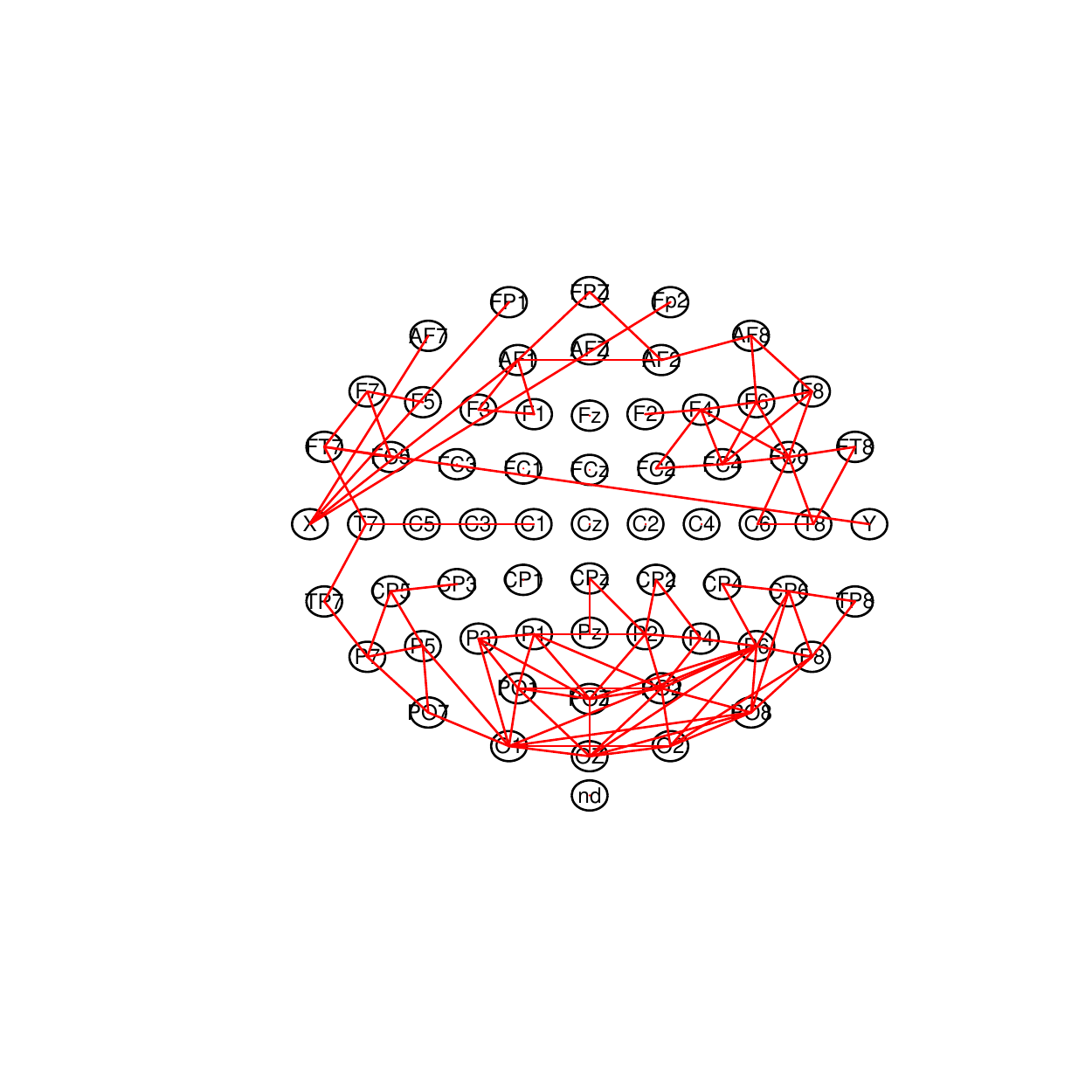} 
   \includegraphics[scale=0.4]{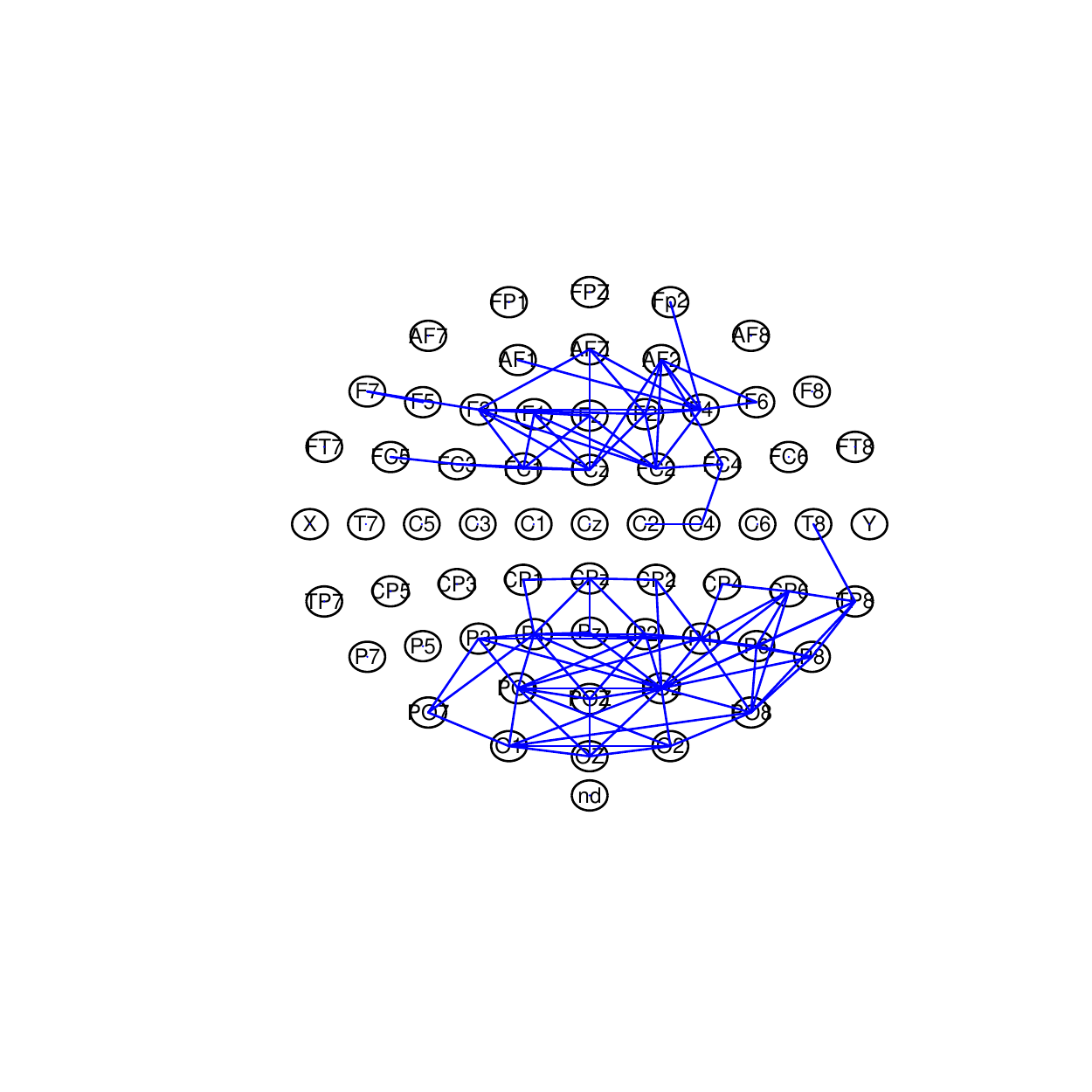}
      \includegraphics[scale=0.4]{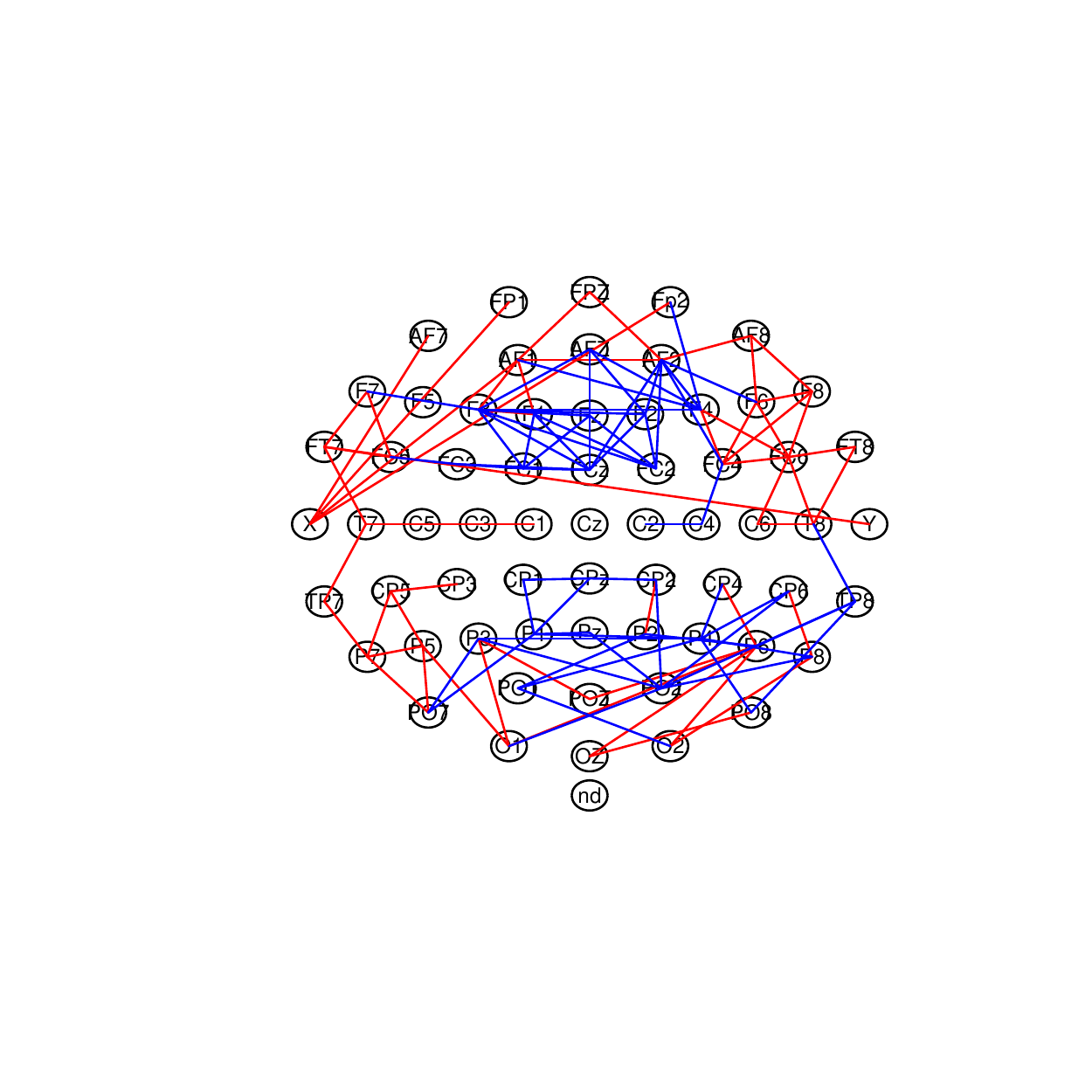}\\
      \includegraphics[scale=0.4]{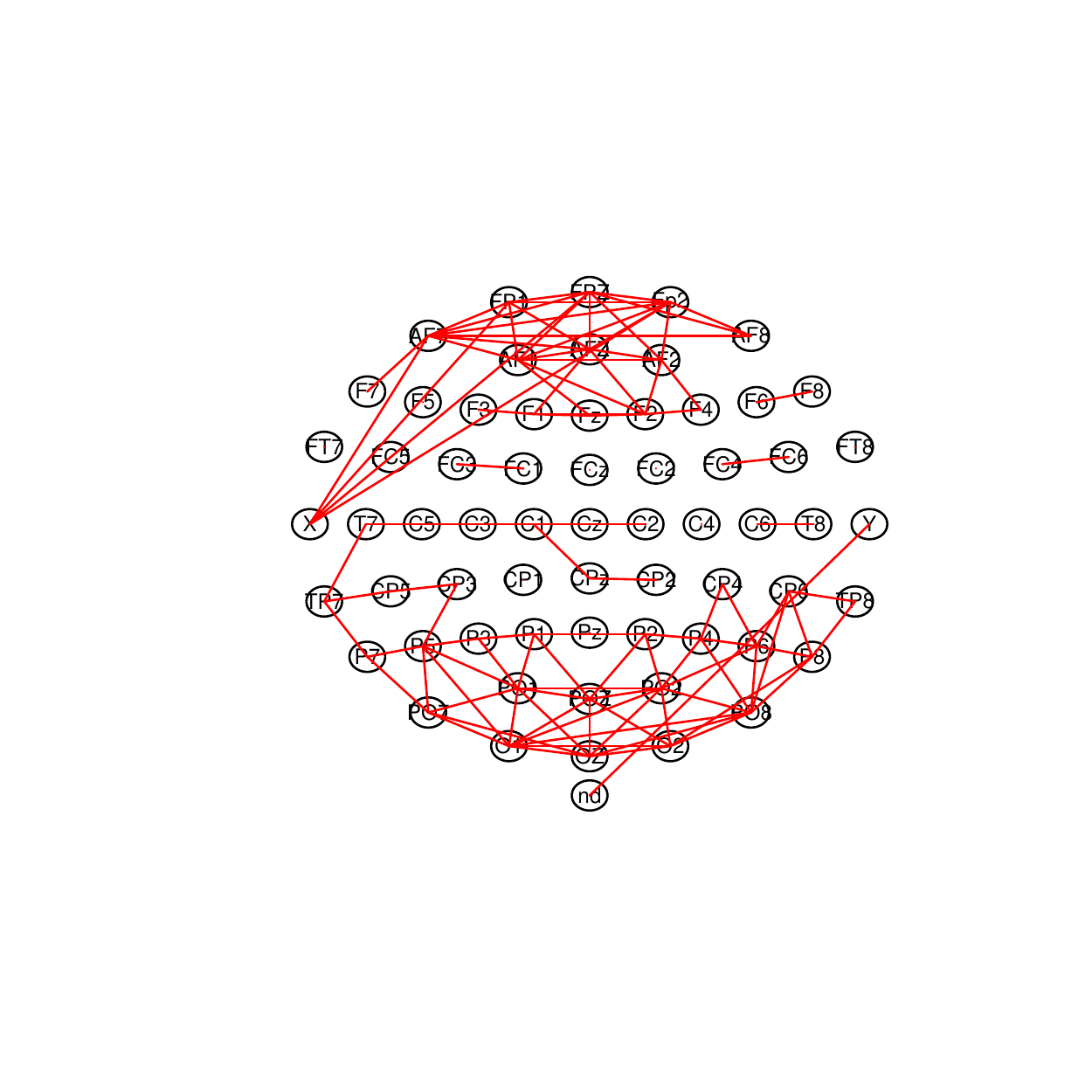} 
            \includegraphics[scale=0.4]{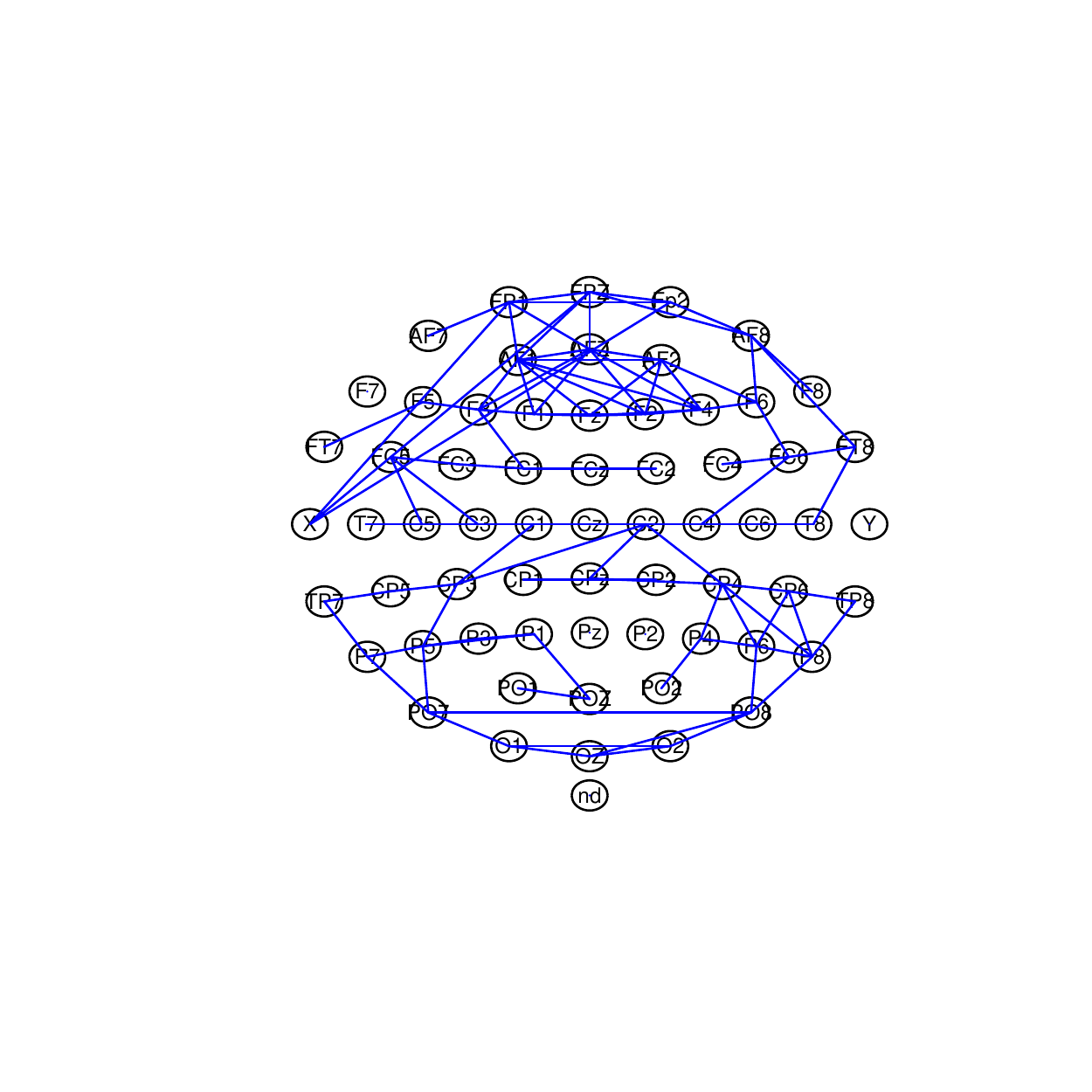} 
      \includegraphics[scale=0.4]{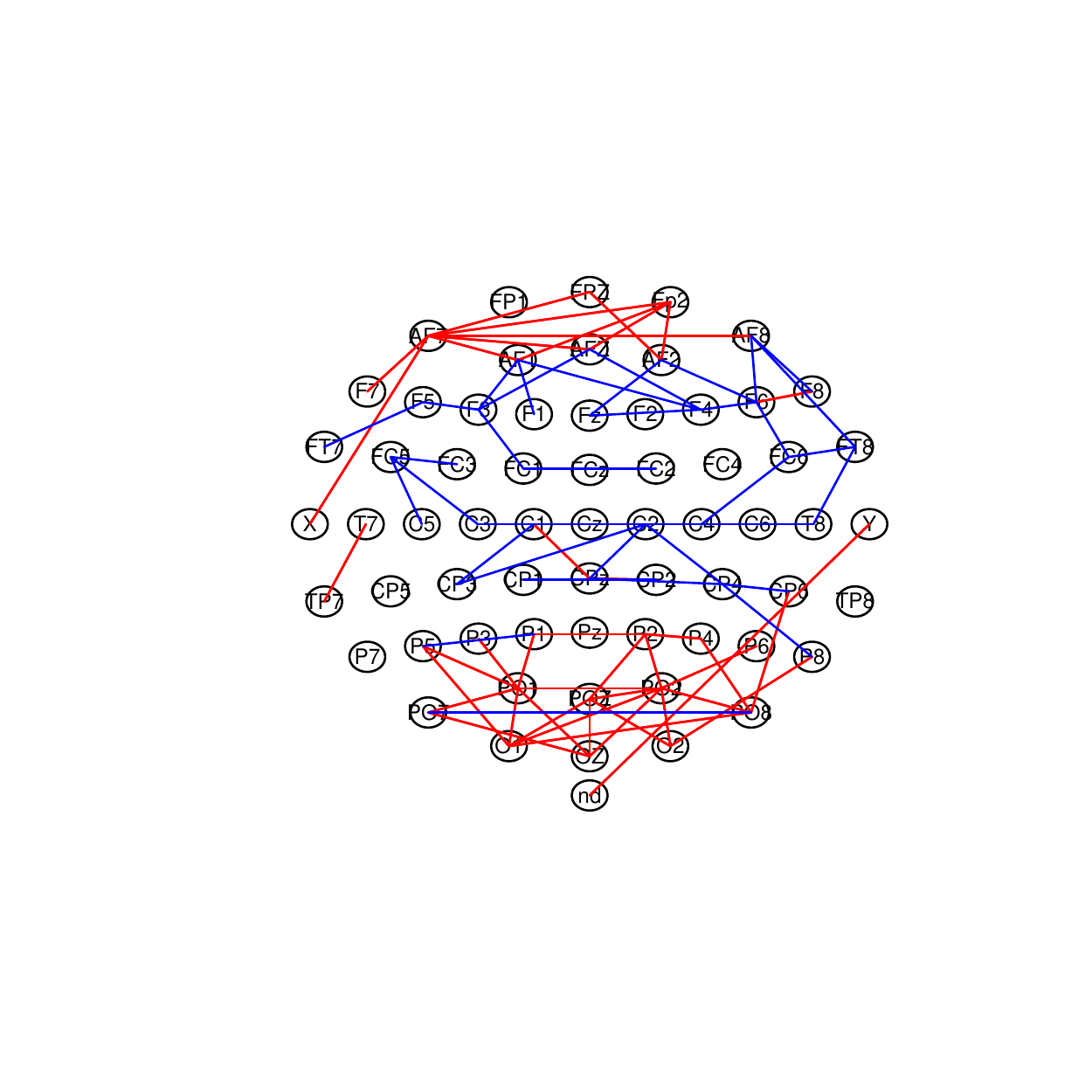} \\
 \includegraphics[scale=0.4]{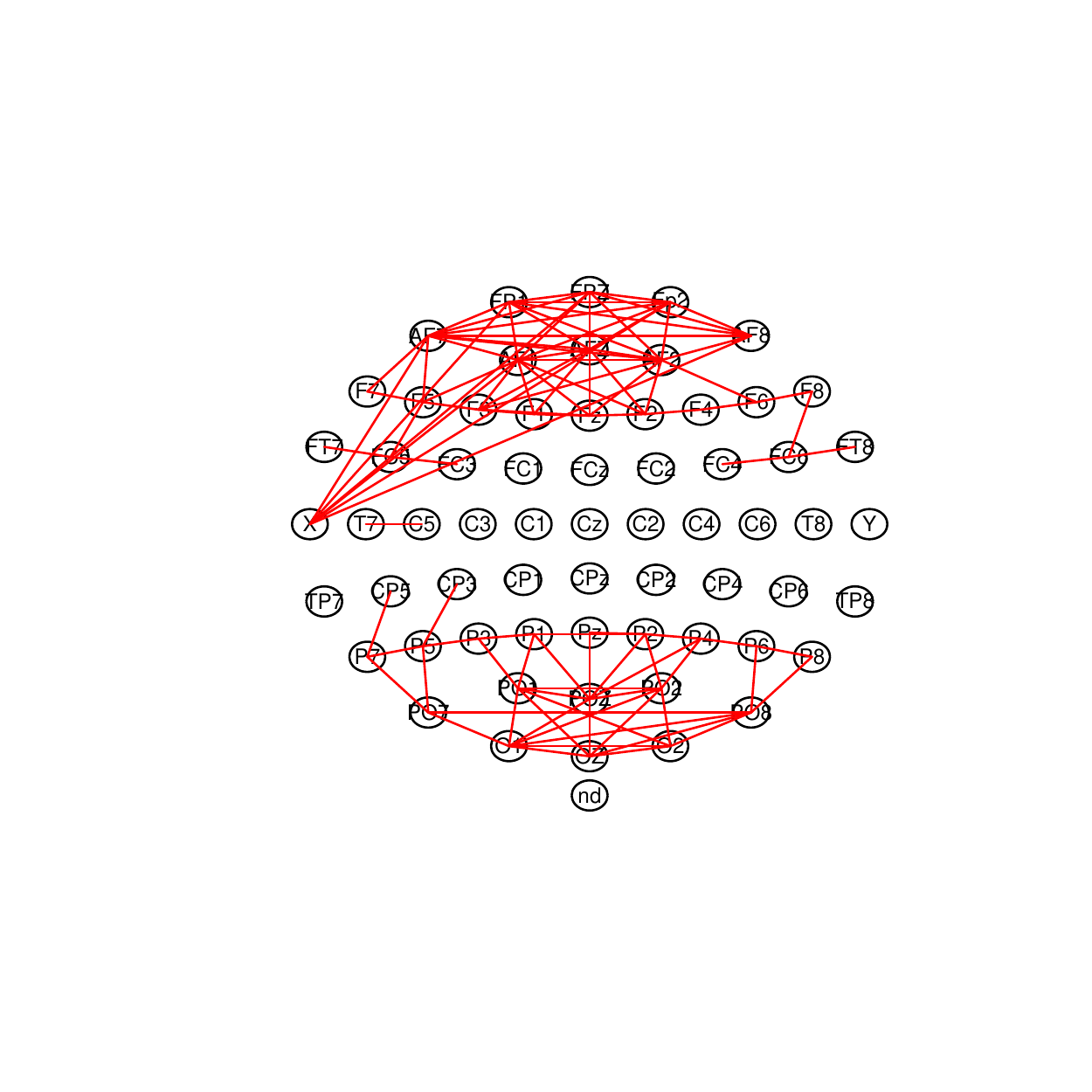}
 \includegraphics[scale=0.4]{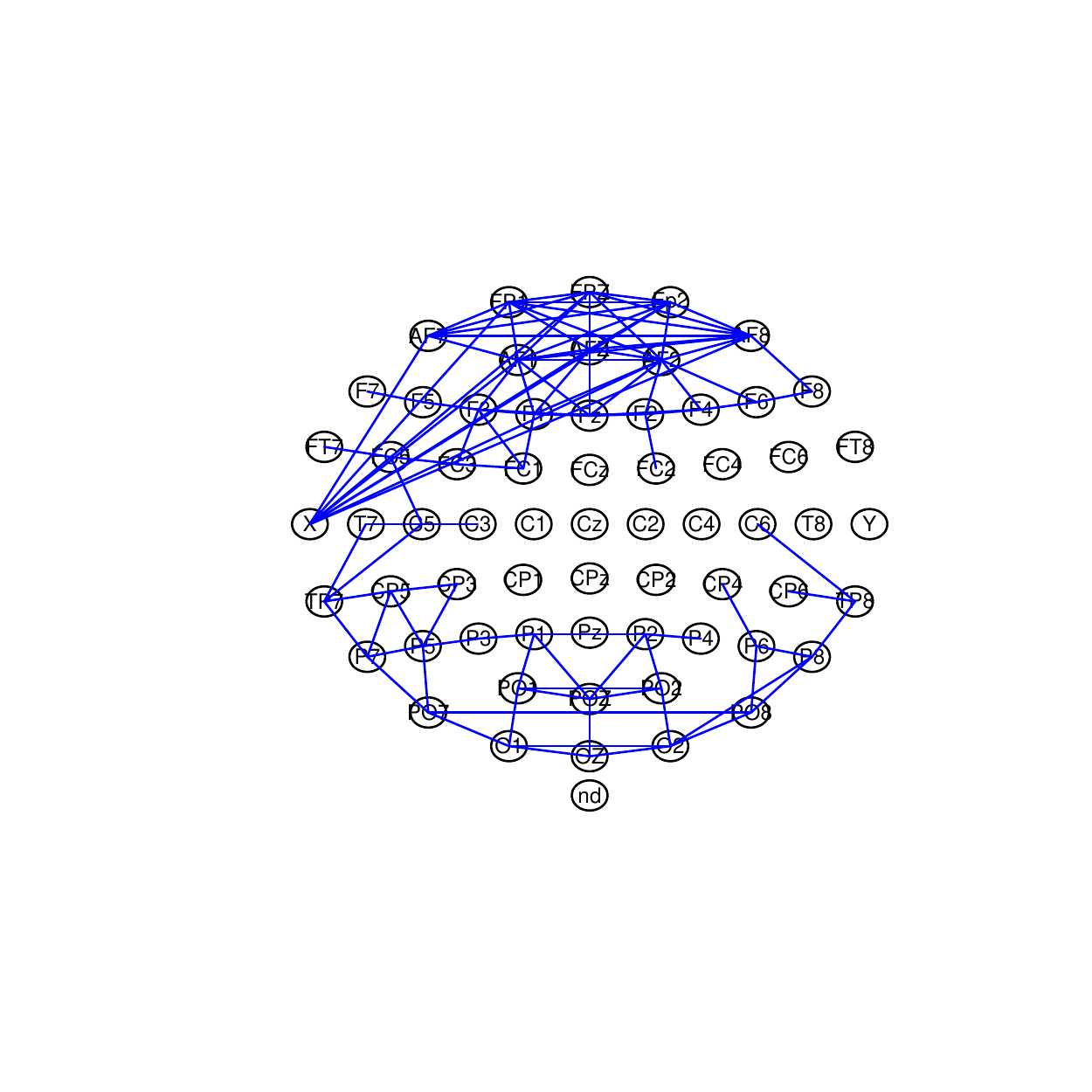}
 \includegraphics[scale=0.4]{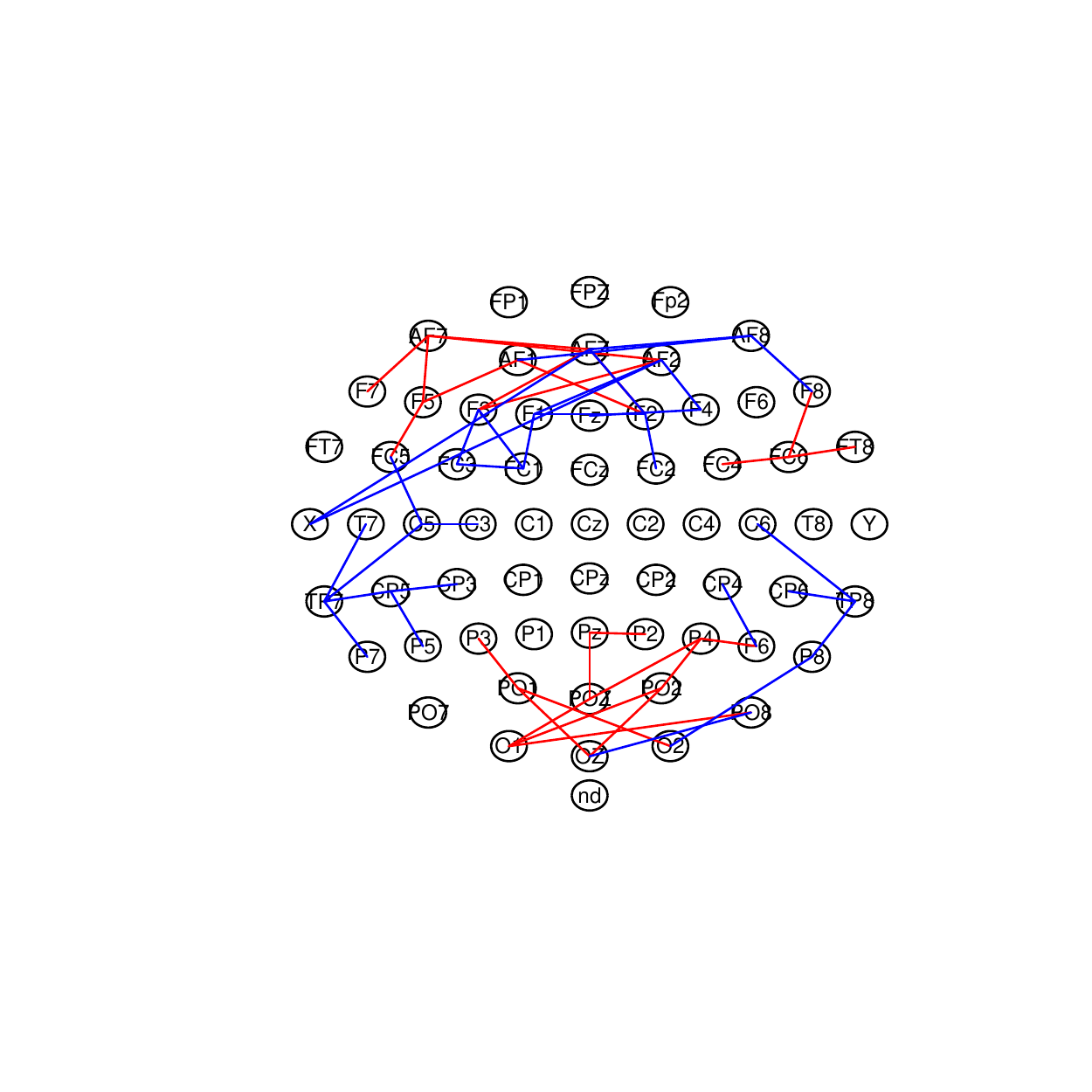}
 \begin{center}
{\small{\parbox{5in}{\raisebox{20pt}{Figure 2.} \ \, \parbox{4.5in}{Estimated brain networks by AFGM (upper panels), FAPO (middle panels) and FGGM (lower panels) for the alcoholic group (left), the control group (middle) and differential brain networks (right).}}}}
\end{center}
\end{center}

 \section{Real data application} \label{sec52}
\def\theequation{6.\arabic{equation}}
\setcounter{equation}{0}

In this section  we apply  the new method to the EEG data set available at  {UCI Machine Learning Repository}. The data involve 77 subjects in the alcoholic group and 45 subjects in the control group. Each subject was exposed to a stimulus while brain  activities were recorded from the 64 electrodes placed  on the subject's scalp,  over a one-second period in which  256 time points were sampled. See \citet{zhang1995event} and  \cite{ingber1997statistical} for more backgrounds of this data. The goal is to characterize functional connectivity among the 64 nodes  for the two groups, based on the functional data collected from the electrodes.

We choose $k \lo n = 4+\lceil \sqrt{n} \rceil$ B-spline functions of order $4$ and number of scores equal to $m \lo n=5$.  Since our goal is to capture outstanding differences in brain connectivity between the alcoholic and control groups, we take the tuning  constant $\lambda \lo n$ to be such that 5\% of the $64 \choose 2$ pairs of vertices are retained as edges.

Figure 2 shows the estimated brain networks constructed  by the three methods for the alcoholic group (left), control group (middle).  The right plots in Figure 2 represent the differential brain networks, where the red lines indicate the edges that are in the alcoholic network but not in the control network, and the blue lines indicate the edges that are in the control network but not in the alcoholic network.
\\
We observe that the brain networks have different patterns for the two groups.
 For example, we observe for all methods, that there is increased functional connectivity in the left frontal area for the alcoholic group relative to the control.

%\begin{center}
 %\includegraphics[scale=0.55]{dif-afgm-graph-eegK13M7.pdf} 
% \begin{center}
%{\small{\parbox{5in}{\raisebox{24pt}{Figure 4.} \ \, \parbox{4in}{Differential brain networks.  The red lines represent the edges that are in the alcoholic network but not in the control; the blue lines indicate the edges that are in the control network but not in the alcoholic.}}}}
%\end{center}
%\end{center}

\section{Conclusions}\label{conclusion}

In this paper, we utilise the idea of generalized additive models to develop a new nonparametric graphical model for multivariate functional data which does not require the assumption of a Gaussian distribution.  The conditional relationships among the principal scores in the Karhunen-Lo{\'e}ve expansion of a random function are allowed to take an arbitrary additive rather than a linear form as imposed by the assumption of Gaussianity.  The additive functions are then approximated by linear combinations of B-splines.  This approximation allows us to develop a group lasso algorithm to estimate the graph that encourages blockwise sparsity to  a matrix formed by the coefficients in the spline approximation.  We have established consistency of the procedure while both the number of principal components and the number of nodes diverge to infinity with increasing sample size.  By simulation study and an analysis of a data example we demonstrate the applicability of the new methodology.

The proposed model and methodology suggests many  directions for future research.  First, the asymptotic results here are developed under the framework where the order $m \lo n$ in  the expansion \eqref{scorestrunc} tends to infinity with the sample size, and it is of interest if similar statistical guarantees  can be obtained in  model \eqref{afm} with the  infinite representation. Second, the  theory is developed under the assumption that the  random functions are fully observed. Therefore, an interesting and important question for future research is the extension of the methodology  to smooth functions that are observed on a dense time grid  such that the covariance operators and the functions are consistently estimated.  Another important direction is the consideration of  the sparse setting, where the functions are observed on a relatively small number of time points and are contaminated with noise. 
In this case, alternative approaches such as \cite{yao2005functional}, \cite{xiao2018fast} and \cite{petrovich2018functional} might be useful and will be further investigated in the future.

\bigskip
\medskip

{\bf Acknowledgements} This work has been supported in part by the
Collaborative Research Center ``Statistical modeling of nonlinear
dynamic processes'' (SFB 823, Teilprojekt A1,C1) of the German Research Foundation
(DFG).  The authors would also like to thank Martina Stein who typed parts of this manuscript with considerable technical expertise.
  
\baselineskip=11pt
\bibliographystyle{asa}
\bibliography{lib_afgm.bib}

\begin{thebibliography}{44}
\newcommand{\enquote}[1]{``#1''}
\expandafter\ifx\csname natexlab\endcsname\relax\def\natexlab#1{#1}\fi

\bibitem[{Bach(2008)}]{bach2008consistency}
Bach, F.~R. (2008), \enquote{Consistency of the group lasso and multiple kernel
  learning,} \textit{Journal of Machine Learning Research}, 9, 1179--1225.

\bibitem[{Bickel and Levina(2008)}]{bickel2008covariance}
Bickel, P.~J. and Levina, E. (2008), \enquote{Covariance regularization by
  thresholding,} \textit{Annals of Statistics}, 36, 2577--2604.

\bibitem[{Bosq(2012)}]{bosq2012linear}
Bosq, D. (2012), \textit{Linear Processes in Function Spaces: Theory and
  Applications}, vol. 149, Springer Science \& Business Media.

\bibitem[{Boucheron et~al.(2013)Boucheron, Lugosi, and
  Massart}]{boucheron2013concentration}
Boucheron, S., Lugosi, G., and Massart, P. (2013), \textit{Concentration
  Inequalities: A Nonasymptotic Theory of Independence}, Oxford University
  Press.

\bibitem[{Cai et~al.(2011)Cai, Liu, and Luo}]{cai2011constrained}
Cai, T., Liu, W., and Luo, X. (2011), \enquote{A constrained $\ell_1$
  minimization approach to sparse precision matrix estimation,} \textit{Journal
  of the American Statistical Association}, 106, 594--607.

\bibitem[{De~Boor(1978)}]{de1978practical}
De~Boor, C. (1978), \textit{A Practical Guide to Splines}, vol.~27,
  springer-verlag New York.

\bibitem[{Friedman et~al.(2008)Friedman, Hastie, and
  Tibshirani}]{friedman2008sparse}
Friedman, J., Hastie, T., and Tibshirani, R. (2008), \enquote{Sparse inverse
  covariance estimation with the graphical lasso,} \textit{Biostatistics}, 9,
  432--441.

\bibitem[{Han et~al.(2018)Han, M{\"u}ller, Park, et~al.}]{han2018smooth}
Han, K., M{\"u}ller, H.-G., Park, B.~U., et~al. (2018), \enquote{Smooth
  backfitting for additive modeling with small errors-in-variables, with an
  application to additive functional regression for multiple predictor
  functions,} \textit{Bernoulli}, 24, 1233--1265.

\bibitem[{Hoeffding(1963)}]{hoeffding1963probability}
Hoeffding, W. (1963), \enquote{Probability inequalities for sums of bounded
  random variables,} \textit{Journal of the American Statistical Association},
  58, 13--30.

\bibitem[{Hsing and Eubank(2015)}]{hsing2015theoretical}
Hsing, T. and Eubank, R. (2015), \textit{Theoretical Foundations of Functional
  Data Analysis, with an Introduction to Linear Operators}, John Wiley \& Sons.

\bibitem[{Huang et~al.(2010)Huang, Horowitz, and Wei}]{huang2010variable}
Huang, J., Horowitz, J.~L., and Wei, F. (2010), \enquote{Variable selection in
  nonparametric additive models,} \textit{Annals of Statistics}, 38,
  2282--2313.

\bibitem[{Ingber(1997)}]{ingber1997statistical}
Ingber, L. (1997), \enquote{Statistical mechanics of neocortical interactions:
  Canonical momenta indicators of electroencephalography,} \textit{Physical
  Review E}, 55, 4578.

\bibitem[{Lam and Fan(2009)}]{lam2009sparsistency}
Lam, C. and Fan, J. (2009), \enquote{Sparsistency and rates of convergence in
  large covariance matrix estimation,} \textit{Annals of Statistics}, 37, 4254.

\bibitem[{Lauritzen(1996)}]{lauritzen1996graphical}
Lauritzen, S.~L. (1996), \textit{Graphical Models}, Oxford University Press.

\bibitem[{Lee et~al.(2016{\natexlab{a}})Lee, Li, and Zhao}]{lee2016additive}
Lee, K.-Y., Li, B., and Zhao, H. (2016{\natexlab{a}}), \enquote{On an additive
  partial correlation operator and nonparametric estimation of graphical
  models,} \textit{Biometrika}, 103, 513--530.

\bibitem[{Lee et~al.(2016{\natexlab{b}})Lee, Li, and Zhao}]{lee2016variable}
--- (2016{\natexlab{b}}), \enquote{Variable selection via additive conditional
  independence,} \textit{Journal of the Royal Statistical Society, Series B
  (Statistical Methodology)}, 78, 1037--1055.

\bibitem[{Lee et~al.(2020)Lee, Li, Li, and Zhao}]{leefaro}
Lee, K.-Y., Li, L., Li, B., and Zhao, H. (2020), \enquote{Nonparametric
  Functional Graphical Modeling through Functional Additive Regression
  Operator,} \textit{Working paper}.

\bibitem[{Li et~al.(2014)Li, Chun, and Zhao}]{liaci}
Li, B., Chun, H., and Zhao, H. (2014), \enquote{On an Additive Semigraphoid
  Model for Statistical Networks With Application to Pathway Analysis,}
  \textit{Journal of the American Statistical Association}, 109, 1188--1204.

\bibitem[{Li and Solea(2018)}]{li2018nonparametric}
Li, B. and Solea, E. (2018), \enquote{A nonparametric graphical model for
  functional data with application to brain networks based on fMRI,}
  \textit{Journal of the American Statistical Association}, 113, 1637--1655.

\bibitem[{Liu et~al.(2012)Liu, Han, Yuan, Lafferty, Wasserman,
  et~al.}]{liu2012high}
Liu, H., Han, F., Yuan, M., Lafferty, J., Wasserman, L., et~al. (2012),
  \enquote{High-dimensional semiparametric Gaussian copula graphical models,}
  \textit{The Annals of Statistics}, 40, 2293--2326.

\bibitem[{Liu et~al.(2009)Liu, Lafferty, and Wasserman}]{liu2009nonparanormal}
Liu, H., Lafferty, J., and Wasserman, L. (2009), \enquote{The nonparanormal:
  Semiparametric estimation of high dimensional undirected graphs,}
  \textit{Journal of Machine Learning Research}, 10, 2295--2328.

\bibitem[{Meier et~al.(2009)Meier, Van~de Geer, B{\"u}hlmann,
  et~al.}]{meier2009high}
Meier, L., Van~de Geer, S., B{\"u}hlmann, P., et~al. (2009),
  \enquote{High-dimensional additive modeling,} \textit{Annals of Statistics},
  37, 3779--3821.

\bibitem[{Meinshausen and B{\"u}hlmann(2006)}]{meinshausen2006high}
Meinshausen, N. and B{\"u}hlmann, P. (2006), \enquote{High-dimensional graphs
  and variable selection with the lasso,} \textit{Annals of Statistics},
  1436--1462.

\bibitem[{Meinshausen et~al.(2009)Meinshausen, Yu,
  et~al.}]{meinshausen2009lasso}
Meinshausen, N., Yu, B., et~al. (2009), \enquote{Lasso-type recovery of sparse
  representations for high-dimensional data,} \textit{Annals of Statistics},
  37, 246--270.

\bibitem[{Obozinski et~al.(2011)Obozinski, Wainwright, Jordan,
  et~al.}]{obozinski2011support}
Obozinski, G., Wainwright, M.~J., Jordan, M.~I., et~al. (2011),
  \enquote{Support union recovery in high-dimensional multivariate regression,}
  \textit{Annals of Statistics}, 39, 1--47.

\bibitem[{Pearl(2002)}]{pearl2002causality}
Pearl, J. (2002), \enquote{Causality: models, reasoning, and inference,}
  \textit{IIE Transactions}, 34, 583--589.

\bibitem[{Peng et~al.(2009)Peng, Wang, Zhou, and Zhu}]{peng2009partial}
Peng, J., Wang, P., Zhou, N., and Zhu, J. (2009), \enquote{Partial correlation
  estimation by joint sparse regression models,} \textit{Journal of the
  American Statistical Association}, 104, 735--746.

\bibitem[{Petrovich et~al.(2018)Petrovich, Reimherr, and
  Daymont}]{petrovich2018functional}
Petrovich, J., Reimherr, M., and Daymont, C. (2018), \enquote{Functional
  regression models with highly irregular designs,} \textit{arXiv preprint
  arXiv:1805.08518}.

\bibitem[{Qiao et~al.(2018)Qiao, Guo, and James}]{qiao2018functional}
Qiao, X., Guo, S., and James, G.~M. (2018), \enquote{Functional graphical
  models,} \textit{Journal of the American Statistical Association}, 114,
  1--12.

\bibitem[{Ravikumar et~al.(2009)Ravikumar, Lafferty, Liu, and
  Wasserman}]{ravikumar2009sparse}
Ravikumar, P., Lafferty, J., Liu, H., and Wasserman, L. (2009), \enquote{Sparse
  additive models,} \textit{Journal of the Royal Statistical Society: Series B
  (Statistical Methodology)}, 71, 1009--1030.

\bibitem[{Ravikumar et~al.(2011)Ravikumar, Wainwright, Raskutti, Yu,
  et~al.}]{ravikumar2011high}
Ravikumar, P., Wainwright, M.~J., Raskutti, G., Yu, B., et~al. (2011),
  \enquote{High-dimensional covariance estimation by minimizing
  $\ell_1$-penalized log-determinant divergence,} \textit{Electronic Journal of
  Statistics}, 5, 935--980.

\bibitem[{Schumaker(2007)}]{schumaker2007spline}
Schumaker, L. (2007), \textit{Spline Functions: Basic Theory}, Cambridge
  University Press.

\bibitem[{Stone et~al.(1985)}]{stone1985additive}
Stone, C.~J. et~al. (1985), \enquote{Additive regression and other
  nonparametric models,} \textit{Annals of Statistics}, 13, 689--705.

\bibitem[{Voorman et~al.(2013)Voorman, Shojaie, and Witten}]{voorman2013graph}
Voorman, A., Shojaie, A., and Witten, D. (2013), \enquote{Graph estimation with
  joint additive models,} \textit{Biometrika}, 101, 85--101.

\bibitem[{Wong et~al.(2019)Wong, Li, and Zhu}]{wong2019partially}
Wong, R.~K., Li, Y., and Zhu, Z. (2019), \enquote{Partially linear functional
  additive models for multivariate functional data,} \textit{Journal of the
  American Statistical Association}, 114, 406--418.

\bibitem[{Xiao et~al.(2018)Xiao, Li, Checkley, and Crainiceanu}]{xiao2018fast}
Xiao, L., Li, C., Checkley, W., and Crainiceanu, C. (2018), \enquote{Fast
  covariance estimation for sparse functional data,} \textit{Statistics and
  computing}, 28, 511--522.

\bibitem[{Xue et~al.(2012)Xue, Zou, et~al.}]{xue2012regularized}
Xue, L., Zou, H., et~al. (2012), \enquote{Regularized rank-based estimation of
  high-dimensional nonparanormal graphical models,} \textit{Annals of
  Statistics}, 40, 2541--2571.

\bibitem[{Yao et~al.(2005)Yao, M{\"u}ller, and Wang}]{yao2005functional}
Yao, F., M{\"u}ller, H.-G., and Wang, J.-L. (2005), \enquote{Functional data
  analysis for sparse longitudinal data,} \textit{Journal of the American
  Statistical Association}, 100, 577--590.

\bibitem[{Yuan and Lin(2006)}]{yuan2006model}
Yuan, M. and Lin, Y. (2006), \enquote{Model selection and estimation in
  regression with grouped variables,} \textit{Journal of the Royal Statistical
  Society: Series B (Statistical Methodology)}, 68, 49--67.

\bibitem[{Yuan and Lin(2007)}]{yuan2007}
--- (2007), \enquote{Model selection and estimation in Gaussian graphical
  model,} \textit{Biometrika}, 94, 19--35.

\bibitem[{Zhang et~al.(1995)Zhang, Begleiter, Porjesz, Wang, and
  Litke}]{zhang1995event}
Zhang, X.~L., Begleiter, H., Porjesz, B., Wang, W., and Litke, A. (1995),
  \enquote{Event related potentials during object recognition tasks,}
  \textit{Brain Research Bulletin}, 38, 531--538.

\bibitem[{Zhou et~al.(1998)Zhou, Shen, Wolfe, et~al.}]{zhou1998local}
Zhou, S., Shen, X., Wolfe, D., et~al. (1998), \enquote{Local asymptotics for
  regression splines and confidence regions,} \textit{Annals of Statistics},
  26, 1760--1782.

\bibitem[{Zhu et~al.(2016)Zhu, Strawn, and Dunson}]{zhu2016bayesian}
Zhu, H., Strawn, N., and Dunson, D.~B. (2016), \enquote{Bayesian graphical
  models for multivariate functional data,} \textit{Journal of Machine Learning
  Research}, 17, 7157--7183.

\bibitem[{Zhu et~al.(2014)Zhu, Yao, and Zhang}]{zhu2014structured}
Zhu, H., Yao, F., and Zhang, H.~H. (2014), \enquote{Structured functional
  additive regression in reproducing kernel Hilbert spaces,} \textit{Journal of
  the Royal Statistical Society: Series B (Statistical Methodology)}, 76,
  581--603.

\end{thebibliography}
\baselineskip=21pt

%%%%%%%%%%%%%%%%%%%%%%%%%%%%
%%%%%%%APPENDIX%%%%%%%%%%%%%%
%%%%%%%%%%%%%%%%%%%%%%%%%%%%

\newpage

\section{Appendix: Proofs}\label{appendix}
\def\theequation{8.\arabic{equation}}
\setcounter{equation}{0}

\subsection{Auxiliary results}  \label{sec81} 

{In this section we state some auxiliary results, which will be  used in the proof of the Theorem  \ref{theorem2section4}.  {The next Lemma provides a concentration inequality for the norm $\| \hat \Sigma \lo {X \hi i X \hi i}-\Sigma \lo {X \hi i X \hi i}\| \lo {HS}$. It can be proved  by similar arguments as given in the proof of Lemma 6 in \cite{qiao2018functional} observing the independence
of the random variables $\xi \hi i \lo {ur}$. The details are omitted for the sake of brevity.}

\begin{lemma}\label{qiaolemma}  Suppose that Assumption \ref{as1} is satisfied. Then, there exists a  constant $C \lo 1$ such that for all $0 < \epsilon \leq C \lo 1$ and for each $i=1, \ldots, p$
\begin{align*}
P \Big(\| \hat \Sigma \lo {X \hi i X \hi i}-\Sigma \lo {X \hi i X \hi i}\| \lo {HS}  \geq  \epsilon \Big) \lesssim  \exp (- C \lo 1 n \epsilon \hi 2).
\end{align*}
\end{lemma}
{{
Let $\xi \hi i=(\xi \hi i \lo {ur})^{1\leq r \leq m \lo n} \lo {1 \leq u \leq n} \in \mathbb{R}^{n \times m \lo n}$ be the matrix of unobserved  scores, and define 
$$
\tilde{h}  \lo n  ( \xi \hi j  \lo {ur})=( \tilde{h} \lo {n1}   ( \xi \hi j  \lo {ur}), \ldots,  \tilde{h}   \lo { nk \lo n}   ( \xi \hi j  \lo {ur}))\trans
$$
as the  vector of the centered $k \lo n$ B-splines functions evaluated at the score $\xi \hi {j} \lo {ur}$, where 
\begin{align}\label{centersamplesplines}
  \tilde{h} \lo {nk} ( \xi \hi j  \lo {ur})= h \lo {k}( \xi \hi j \lo {ur}) -{\frac{1}{n} \tsum \lo {u=1} \hi {n} h \lo {k}( \xi \hi j \lo {ur})}, \quad k=1, \ldots k \lo n.
 \end{align}
 Note that this definition corresponds to \eqref{centertruesplines}, where the expectation has been replaced by its empirical counterpart, and to \eqref{centersplines}, where the estimated scores $\hat \xi \hi j \lo {ur}$ have been replaced by the unobserved scores $\xi \hi j \lo {ur}$.
 Let 
    \begin{eqnarray} \label{hol1}
   \tilde{H} \lo n  (\xi \hi j)  &=&(  \tilde{h} \lo n \trans ( \xi \hi j  \lo {ur})) \lo {1 \leq u \leq n, 1 \leq r \leq m \lo n} \in \R \hi {n \times k \lo n m \lo n}, ~\\
  \mathbf{ \tilde{H} \lo n} \trans (\xi \hi { \sf N \hi i \lo n} )  &=& ( \tilde{H} \lo n (\xi \hi j)  , j \in \sf N  \hi i \lo n) \in \R \hi { n \times n \hi i k \lo n m \lo n}, \label{hol2}
    \end{eqnarray}
    and define
\begin{align}
 \label{samplemat1} 
 \Sigma  \hi n \lo {\sf N \lo n \hi i \sf N \lo n \hi i}
 &= \frac{1}{n} \mathbf{ \tilde{H} \lo n  (\xi \hi { \sf N \hi i \lo n} )} \mathbf{ \tilde{H} \lo n \trans (\xi \hi { \sf N \hi i \lo n} )}   ~ \in \R \hi {n \hi i k \lo n m \lo n  \times n \hi i k \lo n m \lo n},
    \end{align}
    which is the sample analog of the matrix $ \Sigma \hi {*} \lo {\sf N \lo n \hi i \sf N \lo n \hi i}$ defined in \eqref{truemat1}.
Similarly,  let $\hat \xi \hi i =(\hat \xi \hi i \lo {ur}) \lo {1 \leq u \leq n, 1 \leq r \leq m \lo n}$ be the $n \times m \lo n$  matrix of the estimated scores, and 
 \begin{align}\label{estimatedmat1}
\mathbf{ \tilde{H} \lo n} \trans (\hat \xi \hi { \sf N \hi i \lo n})  =(\tilde H \lo n   (\hat \xi \hi {j}), j \in \sf N \hi i \lo n) \in \R \hi { n \times n \hi i k \lo n m \lo n},
     \end{align}
   where 
    \begin{align}\label{estimatedmat2}
    \tilde H \lo n (\hat \xi \hi j)= ( \tilde h \lo {n} \trans (\hat \xi \hi j  \lo {ur})) \lo {1 \leq u \leq n, 1 \leq r \leq m \lo n} \in \R \hi {n \times k \lo n m \lo n},
         \end{align}
         and $\tilde h_n(\hat \xi^j_{ur})$ is defined in \eqref{centersplines1}.
 Then, 
 \begin{align}\label{estimatedmat3}
  \hat \Sigma \hi n \lo { \sf N \hi i \lo n \sf N \hi i \lo n}  = {1 \over n }
  \mathbf{{\tilde H} \lo n}   (\hat \xi \hi { \sf N \hi i \lo n} ) \mathbf{\tilde H \lo n} \trans (\hat \xi \hi { \sf N \hi i \lo n} )  
  \in \R \hi { n \hi i k \lo n m \lo n \times n \hi i k \lo n m \lo n}
  \end{align}
  is the estimated version of the sample design matrix $\Sigma \hi n \lo {\sf N \lo n \hi i \sf N \lo n \hi i}$  in \eqref{samplemat1}.
 The next result provides  tail bounds for all entries of the matrix  $ \hat \Sigma \hi n \lo {\sf N \lo n \hi i \sf N \lo n \hi i}- \Sigma \hi n \lo {\sf N \lo n \hi i \sf N \lo n \hi i}$.

%%%%%%%%%%%%%%%% THEOREM 1%%%%%%%%%%%%%%%%%%%%

\begin{theorem}\label{theorem1} Suppose that Assumption \ref{as1} holds.  Then, there exists a positive constant $C \lo 1$  such that for any $\delta >0$ satisfying $0 < \delta \leq C \lo 1$ and for all $(i , j ) \in \sf V \times \sf V$, $i \ne j $, $r , q=1, \ldots, m \lo n$ and $k , \ell=1, \ldots, k \lo n$, we have
\begin{align*}
P \Big  ( \Big  |\frac{1}{n} \sum \lo {u=1} \hi {n} \left( \tilde h \lo {n k } (\hat \xi \hi {i} \lo {u r })\tilde h \lo {n \ell} (\hat \xi \hi {j } \lo {u q})- \tilde h \lo {nk } ( \xi \hi {i} \lo {u r }) \tilde h \lo {n\ell} ( \xi \hi {j } \lo {u q}) \right) \Big | \geq \delta  \Big )  \lesssim \exp \left( - C \lo 1 {n  \hi {1-\alpha(2+3 \beta)}  k \lo n \hi {-2} \delta \hi 2} \right).
\end{align*}
 \end{theorem}

\proof First, we have
 \begin{align*}
 \big |  \sum \lo {u=1} \hi {n} \big (\tilde h \lo {n k } (\hat \xi \hi { i} \lo {u r })\tilde h \lo {n \ell} (\hat \xi \hi {j } \lo {u q})- \tilde{h} \lo {nk } ( \xi \hi {i} \lo {u r }) \tilde{h} \lo {n\ell} ( \xi \hi {j } \lo {u q}) \big ) \big | \leq T \lo 1 + T \lo 2,
\end{align*}
where the terms $ T \lo 1$ and $ T \lo 2$ are defined as
 \begin{align*}
 T \lo 1= \Big | \sum \lo {u=1} \hi {n}  \tilde{h} \lo {n \ell} ( \xi \hi {j } \lo {u q}) \big(\tilde h \lo {n k} (\hat \xi \hi {i } \lo {u r})- \tilde{h} \lo {nk} ( \xi \hi {i} \lo {u r}) \big) \Big |,    \quad
 T \lo 2  = \Big | \sum \lo {u=1} \hi {n}  \tilde h \lo {n k }(\hat \xi \hi { i} \lo {u r }) \left(\tilde h \lo {n \ell} (\hat \xi \hi {j } \lo {u q})- \tilde{h} \lo {n \ell} ( \xi \hi {j } \lo {u q}) \right)  \Big |.
 \end{align*}
Consequently, for any $\delta >0$, 
\begin{align} \label{theorem1eq1}
P \Big  ( \Big |  \sum \lo {u=1} \hi {n} \left(\tilde h \lo {n k } (\hat \xi \hi { i} \lo {u r })\tilde h \lo {n \ell} (\hat \xi \hi {j } \lo {u q})- \tilde{h} \lo {nk } ( \xi \hi {i} \lo {u r }) \tilde{h} \lo {n\ell} ( \xi \hi {j } \lo {u q}) \right) \Big | \geq 2 n \delta  \Big ) \leq P (T \lo 1 \geq n \delta) + P (T \lo 2 \geq n \delta),
\end{align}
and therefore it is sufficient to derive   inequalities for the two   probabilities on  the right-hand side of \eqref{theorem1eq1}.
\smallskip

\noindent
(a) We start with the probability $P (T \lo 1 \geq n \delta)$.  By the definition of $\tilde h \lo {n k} (\hat \xi \hi {i } \lo {u r})$  and $\tilde h \lo {n k} ( \xi \hi {i } \lo {u r})$ in \eqref{centersplines} and \eqref{centersamplesplines}, respectively,  and some elementary calculations, we obtain for any $\delta >0$
\begin{align*}
 P (T \lo 1 \geq n \delta) \leq P \Big (T \lo {11} \geq \frac{n \delta}{2} \Big)+P \Big (T \lo {12} \geq \frac{n \delta}{2} \Big),
\end{align*}
where
\begin{align*}
T \lo {11}&= \Big | \sum \lo {u=1} \hi {n} h \lo \ell ({\xi \hi j \lo {uq}}) (h \lo k (\hat \xi \hi i \lo {ur})-h \lo k ( \xi \hi i \lo {ur})) \Big |,\\
T \lo {12}&= \Big | n \inv \sum \lo {v=1} \hi {n} h \lo \ell ({\xi \hi j \lo {vq}}) \sum \lo {u=1} \hi n (h \lo k (\hat \xi \hi i \lo {ur})-h \lo k ( \xi \hi i \lo {ur})) \Big |.
\end{align*}
We now derive a concentration inequality for $T \lo {11}$.  By Cauchy-Schwarz inequality and using the fact that $| h \lo \ell ({\xi \hi j \lo {uq}})  | \leq 1$ we have
\begin{align*}
 P \Big (T \lo {11} \geq \frac{n \delta}{2} \Big) \leq P \Big ( \sum \lo {u=1} \hi {n} | h \lo k (\hat \xi \hi i \lo {ur})-h \lo k ( \xi \hi i \lo {ur})| \hi 2 \geq \frac{n \delta \hi 2}{4} \Big),
\end{align*}
and Taylor's expansion gives
\begin{align}\label{taylor}
 \left|  h \lo { k} (\hat \xi \hi {i } \lo {u r})- h \lo {k} ( \xi \hi {i} \lo {u r})  \right|= \Big| \frac{\partial  h \lo {k} (\xi \hi * ) } {\partial \xi \hi {i} \lo {u r}} \left( \hat \xi \hi {i } \lo {u r}-  \xi \hi {i} \lo {u r} \right)\Big|,
\end{align}
where $\xi \hi * $ lies in the line segment between $ \hat \xi \hi {i } \lo {u r}$ and  $\xi \hi {i} \lo {u r}$. From the derivative formula of the B-splines (see \cite{de1978practical}, Ch.10) there exists a constant $M>0$, independent of $n$, such that for all $k=1, \ldots, k \lo n$ and $x \in [-1,1]$,
\begin{align} \label{boundderivativesplines}
\Big | \frac{\partial  h \lo {k} (x ) } {\partial x} \Big | \leq M L \lo n,
\end{align}
where $L \lo n$ is the number of  knots.
As a result, using \eqref{taylor} and \eqref{boundderivativesplines} we obtain
\begin{align*}
P (T \lo {11} \geq n \delta) \leq P \Big (   \sum \lo {u=1} \hi {n} |h \lo { k} (\hat \xi \hi {i } \lo {u r}) - h \lo { k} ( \xi \hi {i } \lo {u r})| \hi 2  \geq \frac{n \delta \hi 2}{4} \Big ) \leq P \Big (  \sum \lo {u=1} \hi {n}  | \hat \xi \hi {i } \lo {u r}-  \xi \hi {i} \lo {u r}| \hi 2  \geq \frac{n \delta \hi 2}{4 M \hi 2 L \lo n \hi 2} \Big ).
\end{align*}
Recall that $\hat \xi \hi {i } \lo {u r}= (\hat \lambda \lo r \hi i) \hi {-1/2} \langle X \hi i \lo u, \hat \phi \lo r \hi i \rangle$ and $ \xi \hi {i } \lo {u r}= ( \lambda \lo r \hi i) \hi {-1/2} \langle X \hi i \lo u,  \phi \lo r \hi i \rangle$.  Then,
\begin{align*}
%\label{diffscores0}
 \hat \xi \hi {i } \lo {u r}-  \xi \hi {i} \lo {u r} =&   (\hat \lambda \lo r \hi i) \hi {-1/2} \langle X \hi i \lo u, \hat \phi \lo r \hi i \rangle - ( \lambda \lo r \hi i) \hi {-1/2} \langle X \hi i \lo u,  \phi \lo r \hi i \rangle   \nonumber \\
\leq &
   ((\hat \lambda \lo r \hi i) \hi {-1/2} -  (\lambda \lo r \hi i) \hi {-1/2} )   \langle X \hi i \lo u, \hat \phi \lo r \hi i \rangle  + ( \lambda \lo r \hi i) \hi {-1/2}  \langle X \hi i \lo u, \hat \phi \lo r \hi i- \phi \lo r \hi i \rangle  . 
\end{align*}
{
By the Cauchy-Schwarz inequality and using the fact that $ \| \hat \phi \lo r \hi i \| =1$ we obtain,
\begin{align*}
%\label{diffscores}
 | \hat \xi \hi {i } \lo {u r}-  \xi \hi {i} \lo {u r} | 
 &\leq  
 |(\hat \lambda \lo r \hi i) \hi {-1/2} -  (\lambda \lo r \hi i) \hi {-1/2}  | \|  X \hi i \lo u \|+  ( \lambda \lo r \hi i) \hi {-1/2}  \|  X \hi i \lo u \|  \| \hat \phi \lo r \hi i- \phi \lo r \hi i  \| \\
 & \leq  
 |(\hat \lambda \lo r \hi i) \hi {-1/2} -  (\lambda \lo r \hi i) \hi {-1/2}   | \|  X \hi i \lo u \|+  ( \lambda \lo r \hi i) \hi {-1/2} d \hi i \lo r  \|  X \hi i \lo u \| \| \hat \Sigma  \lo {X \hi i X \hi i}-\Sigma \lo {X \hi i X \hi i}\| \lo {HS}, \nonumber
\end{align*}
where we have used the inequality 
   $\| \hat \phi \hi {i} \lo r - \phi \hi {i} \lo r\| \leq d \hi i \lo r \| \hat \Sigma \lo {X \hi i X \hi i}-\Sigma \lo {X \hi i X \hi i}\| \lo {HS}
$ 
\citep[see Lemma 4.3 in][]{bosq2012linear} 
and assume w.l.o.g. that  $\hat \phi \hi {i} \lo r$ can be chosen to satisfy $\text{sgn} \langle \hat \phi \hi {i} \lo r, \phi \hi {i} \lo r \rangle=1$   (see the discussion in Remark \ref{rem1}).}
Using the inequality $(a + b) \hi 2 \leq 2 (a \hi 2+ b \hi 2), a, b \in \R$, this implies
\begin{equation}\label{t1}
\begin{split}
P \left(T \lo {11} \geq  {n\delta} \right)  & \leq    P \Big ( |(\hat \lambda \lo r \hi i) \hi {-1/2} -  (\lambda \lo r \hi i) \hi {-1/2}   | \hi 2 \sum \lo {u=1} \hi {n}   \|  X \hi i \lo u \| \hi 2 \geq  \frac{n \delta \hi 2}{16 M \hi 2 L \lo n \hi 2 } \Big) \\ 
& \quad + P \Big (( \lambda \lo r \hi i) \hi {-1} (d \hi i \lo r ) \hi 2 \| \hat \Sigma \lo {X \hi i X \hi i}-\Sigma \lo {X \hi i X \hi i}\| \hi 2 \lo {HS}  \sum \lo {u=1} \hi {n} \|  X \hi i \lo u \| \hi 2 \geq  \frac{n\delta \hi 2}{16 M \hi 2 L \lo n \hi 2 } \Big).
\end{split}
\end{equation}
We now consider the first term at the right-hand side of  \eqref{t1}. Observe that    $ \sum \lo {u=1} \hi {n}   \|  X \hi i \lo u \| \hi 2= \sum \lo {r=1} \hi {\infty} \lambda \lo {r} \hi i \sum \lo {u=1} \hi n (\xi \hi i \lo {ur}) \hi 2 $, and recall that $\xi \hi i \lo {ur} \in [-1,1]$ with $E(\xi \hi i \lo {ur})=1$, $E((\xi \hi i \lo {ur}) \hi 2 )=1$.  Thus $| (\xi \hi i \lo {ur}) \hi 2-1| \leq 2 $, which implies that  for each $r=1, \ldots, m \lo n$,  $ \sum \lo {u=1} \hi n ((\xi \hi i \lo {ur}) \hi 2 -1)$ is a sub-Gaussian random variable with parameter proxy $\sigma \hi 2=4 n$. Consequently, we obtain by Theorem 2.1 of \cite{boucheron2013concentration} 
  \begin{align*}
E \Big \{ \sum \lo {u=1} \hi n ((\xi \hi i \lo {um}) \hi 2 -1) \Big \}  \hi {2k} \leq k ! (16 n) \hi k, \quad k \geq 1.
\end{align*}
Using the convexity of the function $x \mapsto x \hi {2k}$ and Jensen's inequality, it follows
\begin{align*}
E \Big | \sum \lo {u=1} \hi {n}  ( \|  X \hi i \lo u \| \hi 2- E   \|  X \hi i \lo u \| \hi 2) \Big  | \hi {2k} &=E \Big \{ \sum \lo {r=1} \hi {\infty} \lambda \lo {r} \hi i \sum \lo {u=1} \hi n ((\xi \hi i \lo {um}) \hi 2 -1) \Big \} \hi {2k}\\
&\leq \sum \lo {r=1} \hi {\infty} \lambda \lo {r} \hi i E \Big \{ \sum \lo {u=1} \hi n ((\xi \hi i \lo {um}) \hi 2 -1) \Big \}  \hi {2k} (\sum \lo {r=1} \hi {\infty} \lambda \lo {r} \hi i) \hi {2k-1}\\
&\leq k! (16 \lambda \lo 0 \hi 2 n) \hi k,  \quad k \geq 1,
\end{align*}
where $\lambda \lo 0 =\sup \lo {i \leq p} \sum \lo {r=1} \hi {\infty} \lambda \lo {r} \hi i < \infty$ (due to Assumption \ref{as1}).  Hence, from Theorem 2.1 of \cite{boucheron2013concentration}, we obtain for all $\epsilon >0$,
\begin{align*}
P \Big (  \sum \lo {u=1} \hi {n}  ( \|  X \hi i \lo u \| \hi 2- E   \|  X \hi i \lo u \| \hi 2) \geq \epsilon \Big ) & \leq \exp \Big (- \frac{\epsilon \hi 2}{128 \lambda \lo 0 \hi 2 n} \Big ).
\end{align*}
Furthermore, $ E ( \sum \lo {u=1} \hi n \|  X \hi i \lo u \| \hi 2) \leq n \lambda \lo 0$.  Thus, for all $\epsilon /2 \geq n \lambda \lo 0 $,
\begin{align}\label{sumofnorms}
P \Big (  \sum \lo {u=1} \hi {n}   \|  X \hi i \lo u \| \hi 2 \geq \epsilon \Big) \leq P \Big (  \sum \lo {u=1} \hi {n}  ( \|  X \hi i \lo u \| \hi 2- E   \|  X \hi i \lo u \| \hi 2) \geq \epsilon \Big ) &\leq \exp \Big (- \frac{\epsilon \hi 2}{128 \lambda \lo 0 \hi 2 n} \Big ).
\end{align}
Now, we obtain for the first probability on the right-hand side of \eqref{t1}
\begin{align}\label{eq1}
& P \Big ( |(\hat \lambda \lo r \hi i) \hi {-1/2} -  (\lambda \lo r \hi i) \hi {-1/2}   | \hi 2 \sum \lo {u=1} \hi {n}   \|  X \hi i \lo u \| \hi 2 \geq  \frac{n \delta \hi 2}{16 M \hi 2 L \lo n \hi 2 } \Big) \nonumber \\
&\quad  \leq P \Big ( |(\hat \lambda \lo r \hi i) \hi {-1/2} -  (\lambda \lo r \hi i) \hi {-1/2}   | \geq  \frac{ \delta }{2 M \hi {1/2}  L \lo n  } \Big)+P \Big (  \sum \lo {u=1} \hi {n}   \|  X \hi i \lo u \| \hi 2 \geq  \frac{n }{4 M  } \Big).
\end{align}

Define the event 
%\begin{align*}
$
\Omega \lo {m \lo n} \hi i= \{\| \hat \Sigma \lo {X \hi i X \hi i}-\Sigma \lo {X \hi i X \hi i}\| \lo {HS} < 2 \inv \delta \lo {m \lo n} \hi i  \}, 
$
%\end{align*}
where $\delta \lo {m \lo n} \hi i = \min \lo {1 \leq r \leq m \lo n} \{ \lambda \lo r \hi i - \lambda \lo {r+1} \hi i  \}$.  Assumption \ref{as1}(i) implies $\delta \lo {m \lo n} \hi i  \geq d \lo 2 \inv  m \lo n \hi {-(1 + \beta)}$ 
% for some positive constant $d \lo 2$ 
leading to
\begin{equation}\label{omega1}
P ((\Omega \lo {m \lo n} \hi i) \hi {\complement}   )  \leq P (\| \hat \Sigma \lo {X \hi i X \hi i}-\Sigma \lo {X \hi i X \hi i}\| \lo {HS} \geq 2 \inv d \lo 2 \inv m \lo n \hi {-(1 + \beta)} )   \lesssim \exp (- C \lo 1 n  m \lo n \hi {-2(1 + \beta)}  ),
\end{equation}
for some $C \lo 1>0 $, where we have used Lemma \ref{qiaolemma} with $ \delta= 2 \inv d \lo 2 \inv m \lo n \hi {-(1 + \beta)}$.  Furthermore, from Lemma 4.43 of \cite{bosq2012linear}  we have on the event $\Omega \lo {m \lo n} \hi i$
\begin{align*} % \label{bosqeigen}
\sup \lo {r \geq 1 }| \hat \lambda \hi i \lo {r} - \lambda \hi i \lo {r}  | \leq \| \hat \Sigma \lo {X \hi i X \hi i}-\Sigma \lo {X \hi i X \hi i}\| \lo {HS} \leq 2 \inv \delta \hi i \lo {m \lo n} \leq  2 \inv \lambda \hi i \lo {m \lo n}.
\end{align*}
This implies $\hat \lambda \hi i \lo r \geq \frac{\lambda \lo r \hi i}{2} $, $\hat \lambda \hi i \lo r  \leq 2  \lambda \hi i \lo r$ and 
\begin{align*}
 |(\hat \lambda \lo r \hi i) \hi {-1/2} -  (\lambda \lo r \hi i) \hi {-1/2}   |  \leq \frac{(\hat \lambda \hi i \lo r) \hi {-1} |\hat \lambda \hi i \lo r- \lambda \hi i \lo r | (\lambda \hi i \lo r) \hi {-1}}{(\hat \lambda \hi i \lo r) \hi {-1/2} + ( \lambda \hi i \lo r) \hi {-1/2}} \leq 2 ( \lambda \hi i \lo r) \hi {-3/2} |\hat \lambda \hi i \lo r- \lambda \hi i \lo r |.
\end{align*}
This together with     \eqref{omega1} imply that
\begin{align*}
P \Big ( |(\hat \lambda \lo r \hi i) \hi {-1/2} -  (\lambda \lo r \hi i) \hi {-1/2}   | \geq  \frac{ \delta }{2 M \hi {1/2}  L \lo n  } \Big)
 \leq &P \Big ( \Big (|(\hat \lambda \lo r \hi i) \hi {-1/2} -  (\lambda \lo r \hi i) \hi {-1/2}   | \geq  \frac{ \delta }{2 M \hi {1/2}  L \lo n  } \Big ) \cap \Omega \lo {m \lo n} \hi i   \Big )+ P ((\Omega \lo {m \lo n} \hi i ) \hi { \complement} ) \\
 \leq & P \Big ( \Big (|\hat \lambda \lo r \hi i -  \lambda \lo r \hi i   | \geq  \frac{ \delta (\lambda \lo r \hi i) \hi {3/2} }{4 M \hi {1/2}  L \lo n  } \Big )+ P ((\Omega \lo {m \lo n} \hi i ) \hi { \complement} ) \\
  \leq & P \Big ( (\| \hat \Sigma \lo {X \hi i X \hi i}-\Sigma \lo {X \hi i X \hi i}\| \lo {HS} \geq  \frac{ \delta (\lambda \lo r \hi i) \hi {3/2} }{4 M \hi {1/2}  L \lo n  } \Big )+ P ((\Omega \lo {m \lo n} \hi i ) \hi { \complement} ), 
\end{align*} 
where  we used Lemma 4.43 of \cite{bosq2012linear} for the third inequality.  Therefore, from this, Lemma \ref{qiaolemma}, \eqref{sumofnorms} with $\epsilon=\frac{n}{4M}$ and the fact that $d \lo 0  r \hi {-\beta} \leq \lambda  \hi i \lo r$,  the right-hand side of \eqref{eq1} can be upper-bounded by
\begin{align}\label{T11part1}
 & P \Big ( (\| \hat \Sigma \lo {X \hi i X \hi i}-\Sigma \lo {X \hi i X \hi i}\| \lo {HS} \geq  \frac{ \delta (\lambda \lo r \hi i) \hi {3/2} }{4 M \hi {1/2}  L \lo n  } \Big )+ P \Big (  \sum \lo {u=1} \hi {n}   \|  X \hi i \lo u \| \hi 2 \geq  \frac{n }{4 M  } \Big)+P ((\Omega \lo {m \lo n} \hi i ) \hi { \complement} ) \\
   & \quad   \lesssim  \exp (-C \lo 1 n L \lo n \hi {-2} m \lo n \hi {-3 \beta} \delta \hi 2)+\exp (-C \lo 2 n )+ \exp (-C \lo 3 n  m \lo n \hi {-2 (\beta+1)}),  \nonumber 
\end{align} 
 for  all $0< \delta L \lo n \hi {-1} m \lo n \hi {-3 \beta/2} \leq C \lo 1$ and some positive constants $C \lo 1, C \lo 2$ and $C \lo 3$.
 
For the second term on the right-hand side of the inequality \eqref{t1} we use Lemma \ref{qiaolemma}, \eqref{sumofnorms} and the fact that $(d \lo r \hi i) \inv  \geq \frac{d \lo 2}{2\sqrt{2}} m \lo n \hi {-1-\beta}$, to obtain,
    \begin{align}\label{T11part2}
  & P \Big (( \lambda \lo r \hi i) \hi {-1} (d \hi i \lo r)  \hi 2 \| \hat \Sigma \lo {X \hi i X \hi i}-\Sigma \lo {X \hi i X \hi i}\| \hi 2 \lo {HS}  \sum \lo {u=1} \hi {n} \|  X \hi i \lo u \| \hi 2 \geq  \frac{n\delta \hi 2}{16 M \hi 2 L \lo n \hi 2 } \Big) \nonumber \\
 & \quad \leq P \Big ( \| \hat \Sigma \lo {X \hi i X \hi i}-\Sigma \lo {X \hi i X \hi i}\|  \lo {HS} \geq \frac{\delta  (\lambda \lo r \hi i) \hi {1/2} (d \hi i \lo r)  \hi {-1}}{2 M \hi {1/2} L \lo n  } \Big) + P \Big (  \sum \lo {u=1} \hi {n} \|  X \hi i \lo u \| \hi 2 \geq  \frac{n}{4 M  } \Big) \\
   & \quad   \lesssim \exp (-C \lo 4 n ) + \exp (-C \lo 5 n L \lo n \hi {-2} m \lo n \hi {- (2+3\beta)} \delta \hi 2) \nonumber, 
   \end{align}
   for $0< \delta L \lo n \hi {-1}m \lo n \hi {- (2+3\beta)/2}  \leq C \lo 5$ and $C \lo 4>0$, $C \lo 5>0$.

Combining \eqref{T11part1} and \eqref{T11part2} we obtain for all $0 < \delta L \lo n \hi {-1} m \lo n \hi {-3\beta/2}   \leq C \lo 3$
the inequality
\begin{align}\label{t11final}
 P \left(T \lo {11} \geq  {n\delta} \right) \lesssim \exp (-C \lo 1 n )+ \exp (-C \lo 2 n  m \lo n \hi {-2 (\beta+1)})+ \exp (-C \lo 3 n L \lo n \hi {-2} m \lo n \hi {- (2+3\beta)} \delta \hi 2).
\end{align}
We now consider the probability $P \left(T \lo {12} \geq  {n\delta} \right) $.  Using the fact that  $| h \lo \ell ({\xi \hi j \lo {uq}})  | \leq 1$ and Taylor's expansion \eqref{taylor} yield 
\begin{eqnarray*}
P \left(T \lo {12} \geq  {n\delta} \right) & \leq &  P \Big( \Big | \sum \lo {u=1} \hi n (\hat \xi \hi i \lo {ur} - \xi \hi i \lo {ur}) \Big | \geq \frac{n \delta}{2 M L \lo n}  \Big ) 
\leq  P \Big(  \sum \lo {u=1} \hi n | \hat \xi \hi i \lo {ur} - \xi \hi i \lo {ur} | \hi 2  \geq \frac{n \delta \hi 2}{4 M \hi 2 L \lo n \hi 2}  \Big ) ,
\end{eqnarray*}
where we have used the Cauchy-Schwarz inequality.
Therefore, by  similar arguments as  used in the derivation of the bound for  $P \left(T \lo {11} \geq  {n\delta} \right)$, there exist positive constants $C \lo 4$, $C \lo 5$ and $C \lo 6$ such that for all $0< \delta L \lo n \hi {-1} m \lo n \hi {-3 \beta /2}  \leq C \lo 6$ 
  \begin{equation}\label{t12final}
  %\begin{split}
  P \left(T \lo {12} \geq  {n\delta} \right)  
%\leq P \Big( \Big | \sum \lo {u=1} \hi n (\hat \xi \hi i \lo {ur} - \xi \hi i %\lo {ur}) \Big | \geq \frac{n \delta}{2 M L \lo n}  \Big )  \\
   \lesssim \exp (-C \lo 4 n )+ \exp (-C \lo 5 n  m \lo n \hi {-2 (\beta+1)})+ \exp (-C \lo 6 n L \lo n \hi {-2} m \lo n \hi {-(2+3 \beta)} \delta \hi 2).
%\end{split}
\end{equation}
  Combining \eqref{t11final} and  \eqref{t12final}  and choosing suitable constants, we obtain for $0< \delta L \lo n \hi {-1} m \lo n \hi {-3 \beta/2} \leq C \lo 1$, 
  \begin{align}\label{t1finala}
P (T \lo 1 \geq n \delta)   \lesssim  \exp (-C \lo 1 n L \lo n \hi {-2} m \lo n \hi {-(2+3 \beta)} \delta \hi 2). 
\end{align}
 %%%%%
 %(b)
 %%%%%
(b) To derive a bound for the term $P \left(T \lo {2} \geq  {n\delta} \right)$ in \eqref{theorem1eq1} we use the decomposition
\begin{align*}
P \left(T \lo {2} \geq  n\delta \right)  \leq & P \Big( \Big| \sum \lo {u=1} \hi {n} \left(\tilde h \lo {n k }(\hat \xi \hi { i} \lo {u r }) -  \tilde{h} \lo { nk }( \xi \hi { i} \lo {u r }) \right) \left( \tilde h \lo {n \ell }(\hat \xi \hi { j} \lo {u q }) -  \tilde{h} \lo {n \ell }( \xi \hi { j} \lo {u q }) \right)\Big|  \geq  \frac{n\delta}{2} \Big)\\
&+ P \Big(   \Big|  \sum \lo {u=1} \hi {n} \tilde{h} \lo { k }( \xi \hi { i} \lo {u r }) \left( \tilde h \lo {n \ell }(\hat \xi \hi { j} \lo {u q }) - \tilde{ h} \lo { \ell }( \xi \hi { j} \lo {u q }) \right) \Big|  \geq  \frac{n\delta}{2} \Big) \\
=& P \Big( T_{21} \geq \frac{n \delta}{2}\Big) + P \Big( T_{22} \geq \frac{n \delta}{2} \Big),
\end{align*}
where the last inequality defines the terms $T \lo {21}$ and $T \lo {22}$ in an obvious manner.  For the second term, we obtain by the same arguments as used to estimate   $ P \left(T \lo {1} \geq  {n\delta} \right)$    for any $0< \delta L \lo n \hi {-1} m \lo n \hi {-3 \beta/2} \leq C \lo 2$, 
\begin{align}\label{t22final}
P \Big( T \lo {22} \geq  \frac{n\delta}{2} \Big)  \lesssim  \exp (- C \lo 2 n L \lo n \hi {-2} m \lo n \hi {-(2+3 \beta)} \delta \hi 2 ).
\end{align}
For the first term, we use the definition of the centred B-splines in \eqref{centersplines} and \eqref{centersamplesplines}  to obtain
for any $\delta>0$
\begin{align*} %\label{t2}
P \Big( T \lo {21} \geq  \frac{n\delta}{2} \Big) \leq  P \Big(T \lo {211} \geq  \frac{n \delta}{4 } \Big) + P \Big(T \lo {212} \geq  \frac{n \delta}{4 } \Big)~,
\end{align*}
where       
\begin{align*}
 T \lo {211}& = \Big| \sum \lo {u=1} \hi {n}  \left( h \lo {k }(\hat \xi \hi { i} \lo {u r }) -  h \lo { k }( \xi \hi { i} \lo {u r }) \right)  \left( h \lo { \ell }(\hat \xi \hi { j} \lo {u q }) -  h \lo { \ell }( \xi \hi { j} \lo {u q }) \right) \Big|, \\
 T \lo {212} & =  n \inv  \Big| \sum \lo {u=1} \hi {n}  \left(   h \lo {k }(\hat \xi \hi { i} \lo {u r }) -  h \lo { k }( \xi \hi { i} \lo {u r }) \right) \Big| \Big| \sum \lo {u=1} \hi {n}  \left(   h \lo { \ell }(\hat \xi \hi { j} \lo {u q }) -h \lo { \ell }( \xi \hi { j} \lo {u q }) \right) \Big|.  
\end{align*}
{
To derive a concentration bound for the first term, we use \eqref{taylor}, \eqref{boundderivativesplines} and the Cauchy-Schwarz inequality to  obtain
\begin{align*}
P \Big(T \lo {211} \geq  \frac{n \delta}{4} \Big) &\leq P \Big( \sum \lo {u=1} \hi {n} |  \hat \xi \hi { i} \lo {u r } -   \xi \hi { i} \lo {u r } | |  \hat \xi \hi { j} \lo {u q } -   \xi \hi { j} \lo {u q } | \geq  \frac{n \delta}{4 M \hi 2 L \lo n \hi 2} \Big) 
 \\
 &  \leq 2 P \Big( \sum \lo {u=1} \hi {n} |  \hat \xi \hi { i} \lo {u r } -   \xi \hi { i} \lo {u r } | \hi 2 \geq  \frac{n \delta}{4 M \hi 2 L \lo n \hi 2} \Big)
 \lesssim       \exp (- C \lo {3} n L \lo n \hi {-2} m \lo n \hi {-(2+3\beta)} \delta),
\end{align*}
for a positive constant $C \lo 3$ such that for  $0 < \delta \hi {1/2}  L \lo n \hi {-1} m \lo n \hi {-3\beta/2}    \leq C \lo {3}$.}
Here the last inequality follows by  the same arguments as used for the bound of $P \left(T \lo {11} \geq  {n\delta} \right)$.     
 Finally,  for the term $P \left(T \lo {212} \geq  \frac{n \delta}{4 } \right)$ we have
\begin{align*}
P \Big(T \lo {212} \geq  \frac{n \delta}{4 } \Big)=P \Big( \Big| \sum \lo {u=1} \hi {n}  \left(   h \lo {k }(\hat \xi \hi { i} \lo {u r }) -  h \lo { k }( \xi \hi { i} \lo {u r }) \right) \Big| \Big| \sum \lo {u=1} \hi {n}  \left(   h \lo { \ell }(\hat \xi \hi { j} \lo {u q }) -h \lo { \ell }( \xi \hi { j} \lo {u q }) \right) \Big|  \geq \frac{n \hi 2 \delta}{4}\Big).
\end{align*}
Thus, we can apply similar techniques as used in the derivation of the bound \eqref{t11final} with  $\frac{n \delta}{2}$ replaced by  $\frac{n \delta \hi {1/2}}{2}$ leading to
 \begin{align}\label{E2final}
  P \Big  ( T \lo {212} \geq  \frac{n\delta}{4} 
  \Big ) \lesssim    \exp (- C \lo 4 n L \lo n \hi {-2} m \lo n \hi {-(2+3\beta)} \delta ),
 \end{align}
for $0 < \delta \hi {1/2}  L \lo n \hi {-1} m \lo n \hi {-3\beta/2} \leq C \lo 4$.
Combining these results with \eqref{theorem1eq1} and \eqref{t1finala} and using the fact that $m \lo n \asymp n \hi \alpha$ and $k \lo n > L \lo n$, we   obtain the assertion of the Theorem.\eop
 
 \noindent
Theorem \ref{theorem1} states that  the elements of the matrix $ \hat \Sigma \lo {\sf N \lo n \hi i \sf N \lo n \hi i}- \Sigma \hi n \lo {\sf N \lo n \hi i \sf N \lo n \hi i}$  exhibit  exponential-type probability tails.  We also observe that the decay rate $\beta$ of the eigenvalues appears in the tail behaviour.  A similar condition of exponential tails is imposed on the elements of the sample covariance matrix of scalar and functional Gaussian data for the analysis of high-dimensional Gaussian graphical models \citep[see, for example,][]{ravikumar2011high,qiao2018functional}.

{ \begin{proposition}\label{splineapproximation}  Suppose that Assumptions \ref{as3}, \ref{as4} and condition \eqref{fbounded} are satisfied.  Then, there exist functions $\tilde{f} \lo {nqr} \hi {ij}= \tsum \lo {k=1} \hi {k \lo n}  \beta \lo {qrk} \hi {ij} \tilde {h} \lo {nk} $ and  positive constants $c_1,C_1$, such that
{\begin{align*}
%\label{eventspline}
P \big ( \Omega^\complement \big) \leq 2 \exp \Big (-C \lo 1 \frac{n k \lo n \hi {-2d}}{n \hi i  m \lo n \hi 2} + \log (n \hi i m \lo n \hi 2) \Big ),
\end{align*}
where 
\begin{equation} \label{omega}
\Omega = \Big \{ \max \lo {j \in \sf N \lo n \hi i  } \max \lo {1 \leq q, r \leq m \lo n}  \frac{1}{\sqrt{n}} \| 
\mathbf {f \hi {ij} \lo {qr}}
- \mathbf {\tilde f \hi {ij} \lo {qr}} \| \lo {2} < c \lo 1 k \lo n \hi {-d}  \Big \}  ~,
\end{equation}}
and 
$\mathbf {f \hi {ij} \lo {qr}}  =  \big ( f \hi {ij} \lo {qr}(\xi \hi {j} \lo {1r}) ,  \ldots , f \hi {ij} \lo {qr}(\xi \hi {j} \lo {nr})  \big)^{\top}$,
$\mathbf {\tilde f \hi {ij} \lo {qr}}  =  \big ( \tilde  f \hi {ij} \lo {qr}(\xi \hi {j} \lo {1r}) ,  \ldots , \tilde  f \hi {ij} \lo {qr}(\xi \hi {j} \lo {nr})  \big)^{\top}$.
\end{proposition} 
}
\proof
By Assumptions \ref{as3} and \ref{as4} for any $f \hi {ij} \lo {qr} \in \sten {F} \lo  {{\kappa,\rho}}$ there exists a $B$-spline $ g^{ij}_{qr}=\tsum \lo {k=1} \hi {k \lo n} \beta \lo {qrk} \hi {ij} {h} \lo {k}  \in \mathcal{S} \lo {\ell L \lo n}$ and a positive constant $c \lo 1$ such that
\begin{align*}
 \| f \hi {ij} \lo {qr}  - g \hi {ij} \lo {qr} \| \lo {\infty} \leq c \lo 1 k \lo n \hi {-d},
 \end{align*}
 (see Lemma 5 in \cite{stone1985additive}). 
 Let  
 $\tilde{f} \lo {qr} \hi {ij} (\xi^j_{ur}) =g \hi {ij} \lo {qr} (\xi^j_{ur})  -\frac{1}{n} \tsum \lo {u=1} \hi {n} g \hi {ij} \lo {qr} (\xi \hi {j} \lo {ur})$.  Then recalling the notation \eqref{omega1} we have $\tilde{f} \lo {qr} \hi {ij} (\xi^j_{ur})= \tsum \lo {k=1} \hi {k \lo n} \beta \lo {qrk} \hi {ij} \tilde {h} \lo {nk} (\xi^j_{ur}) $, and
 we obtain
\begin{align*}
 \frac{1}{n} \| \mathbf {f \hi {ij} \lo {qr}}
- \mathbf {\tilde f \hi {ij} \lo {qr}} \| \lo {2} \hi 2 
& 
  \leq 2 c \lo 1 \hi 2 k \lo n \hi {-2d}+ 2 \Big( \frac{1}{n} \sum \lo {u=1} \hi {n} g \hi {ij} \lo {qr}(\xi \hi {j} \lo {ur}) \Big) \hi 2 \\
& \leq 2 c \lo 1 \hi 2 k \lo n \hi {-2d}+ 4 \Big( \frac{1}{n} \sum \lo {u=1} \hi {n} (g \hi {ij} \lo {qr}(\xi \hi {j} \lo {ur}) -  f \hi {ij} \lo {qr}(\xi \hi {j} \lo {ur})) \Big) \hi 2 +4 \Big( \frac{1}{n} \sum \lo {u=1} \hi {n} f \hi {ij} \lo {qr}(\xi \hi {j} \lo {ur}) \Big) \hi 2 \\
& \leq 6 c \lo 1 \hi 2 k \lo n \hi {-2d}+4 \Big( \frac{1}{n} \sum \lo {u=1} \hi {n} f \hi {ij} \lo {qr}(\xi \hi {j} \lo {ur}) \Big) \hi 2. 
\end{align*}
{From this, condition \eqref{fbounded}, Hoeffding's inequality and the union bound, it follows with an appropriate constant $c_2>0$
\begin{align*} 
P \Big( \max \lo {j \in \sf N \lo n \hi i  } \max \lo {1 \leq q, r \leq m \lo n}  \frac{1}{n} \| f \hi {ij} \lo {qr}(\xi \hi {j} \lo {r} )- \tilde{f} \lo {nqr} \hi {ij} (\xi \hi {j} \lo {r} ) \| \lo {2} \hi 2 \geq  c \lo 1 \hi 2 k \lo n \hi {-2d} \Big) 
&\leq P \Big( \max \lo {j \in \sf N \lo n \hi i  } \max \lo {1 \leq q, r \leq m \lo n} \Big | \frac{1}{n} \sum \lo {u=1} \hi {n} f \hi {ij} \lo {qr}(\xi \hi {ij} \lo {ur}) \Big|  \geq c \lo 2  k \lo n \hi {-d} \Big) \\ \nonumber
&\leq 2 \exp \Big(-  \frac{n k \lo n \hi {-2d}}{{2 M \hi 2 n \hi i} m \lo n \hi 2}  + \log (n \hi i m \lo n \hi 2) \Big)~,
\end{align*}   }
which  completes the proof. \eop}

\noindent
 
\medskip

%%%%%%%%%%%%%%%%PART II%%%%%%%%%%%%%%%%%%%%

\subsection{Rates of convergence for sample design matrices}\label{sec82} 

In this section we show that if   Assumptions \ref{boundedeigen} and \ref{irrepr} hold, then with high probability, the assumptions hold also for the  corresponding {\it {sample matrices}} 
\begin{align}\label{samplematrices} 
 \Sigma \hi n \lo {\sf N \lo n \hi i \sf N \lo n \hi i} 
 = \frac{1}{n} \mathbf{ \tilde{H} \lo n  (\xi \hi { \sf N \hi i \lo n} ) \tilde{H} \lo n \trans (\xi \hi { \sf N \hi i \lo n} ) } \in \mathbb{R}^{n^ik_nm_n \times n^i k_nm_n},  ~   \Sigma \hi n \lo {{\xi \hi j}  \sf N \lo n \hi i  }
  =  \frac {1}{n} \tilde{H} \trans \lo n ( \xi \hi  j) \mathbf{\tilde{H}\trans \lo n (\xi \hi { \sf N \hi i \lo n} )  } \in \mathbb{R}^{k_nm_n \times   k_nm_n n^i},
\end{align}
where $\tilde{H}  \lo n ( \xi \hi  j) $ and $\mathbf{ \tilde{H} \lo n  (\xi \hi { \sf N \hi i \lo n} ) }$ 
 are defined in \eqref{hol1} and \eqref{hol2}, respectively. Note that the matrices in \eqref{samplematrices} are
 based on   the unobserved scores  and are the sample analogs of the matrices $ \Sigma \hi {*} \lo {\sf N \lo n \hi i \sf N \lo n \hi i}$  and   $ \Sigma \hi * \lo {{\xi \hi j}  \sf N \lo n \hi i }$  in \eqref{truemat1} and \eqref{truemat2}, respectively.

   \begin{lemma}\label{sampleboundedeigen} Suppose that Assumption \ref{boundedeigen} holds.  Then, there exists a constant $C \lo 1>0$ such that for any $\delta>0$, 
      \begin{align}\label{boundtruesigma}
      P \left(  \|  \Sigma \hi n \lo {  {\sf N \hi i \lo n }  {\sf N \hi i \lo n } } -  \Sigma \hi * \lo {  {\sf N \hi i \lo n }  {\sf N \hi i \lo n} }  \| \lo F \geq \delta \right) \leq 2 \exp \Big(- C \lo 1 \frac{n \delta \hi 2}{(n \hi i m \lo n k \lo n) \hi 2} + 2 \log (n \hi i m \lo n k \lo n)\Big).
         \end{align}
            \begin{align}\label{sampleboundedeigencondition}
   P \left(  \Lambda \lo {\min} ( \Sigma \hi n \lo {\sf N \lo n \hi i \sf N \lo n \hi i}) \leq C \lo {\min}  - \delta \right) \leq 2 \exp \Big( -C \lo 1  \frac{n\delta \hi 2}{(n \hi i m \lo n k \lo n) \hi 2}+ 2 \log (n \hi i m \lo n k \lo n)\Big).
   \end{align}
   \end{lemma}
   \proof    Weyl's Lemma yields 
\begin{equation*}
%\label{weyl}
\Lambda \lo {\min} (\Sigma \hi * \lo {  {{\sf N \hi i \lo n}}  {{\sf N \hi i \lo n} } })- \Lambda \lo {\min} ( \Sigma \hi n  \lo {  {\sf N \hi i \lo n}   {\sf N \hi i \lo n}}  )  \leq \|  \Sigma \hi n  \lo {  {\sf N \hi i \lo n}   {\sf N \hi i \lo n} } -  \Sigma \hi * \lo {  {\sf N \hi i \lo n}   {\sf N \hi i  \lo n} }  \| \lo 2   \leq \|  \Sigma \hi n \lo {  {\sf N \hi i \lo n}   {\sf N \hi i \lo n }  } -  \Sigma \hi * \lo {  {\sf N \hi i \lo n }  {\sf N \hi i \lo n }}  \| \lo F, 
\end{equation*}
  and   by Assumption \ref{boundedeigen} we have,
   \begin{align}\label{weyls}
P \left( \Lambda \lo {\min} ( \Sigma \hi n  \lo {  {\sf N \hi i \lo n }  {\sf N \hi i \lo n }} ) \leq C \lo {\min} - \delta \right)\leq  P \left(  \|  \Sigma \hi n \lo {  {\sf N \hi i  \lo n} {\sf N \hi i \lo n }} -  \Sigma \hi * \lo {  {\sf N  \hi i \lo n } {\sf N \hi i \lo n}  }  \| \lo F \geq \delta \right).
\end{align}
By definition the $n \hi i m \lo n k \lo n \times n \hi i m \lo n k \lo n$ matrix $\Sigma \hi n  \lo {  {\sf N} \hi i \lo n  {\sf N \hi i \lo n} } -\Sigma \hi * \lo {  {\sf N \hi i \lo n}   {\sf N \hi i \lo n} }$ contains elements of the form
\begin{align*}
W \lo {k \ell, rq} \hi {j \lo 1, j \lo 2} =\frac{1}{n} \sum \lo {u=1} \hi n \tilde h \lo {nk} (\xi \hi {j \lo 1} \lo {uq})  \tilde h \lo {n\ell} (\xi \hi {j \lo 2} \lo {ur})- E(\tilde {h \lo k} (\xi \hi {j \lo 1} \lo q) \tilde {h \lo \ell} (\xi \hi {j \lo 2} \lo r)),
\end{align*}
which can be rewritten (recalling the notation \eqref{centersamplesplines}) as $A \lo 1- A \lo 2$, where
\begin{align*}
A \lo 1 &=\frac{1}{n} \sum \lo {u=1} \hi n  h \lo {k} (\xi \hi {j \lo 1} \lo {uq})   h \lo {\ell} (\xi \hi {j \lo 2} \lo {ur})-E( {h \lo k} (\xi \hi {j \lo 1} \lo q) {h \lo \ell} (\xi \hi {j \lo 2} \lo r)),\\
A \lo 2 &=\frac{1}{n \hi 2} \sum \lo {u \lo 1=1} \hi n \sum \lo {u \lo 2=1} \hi n h \lo {k} (\xi \hi {j \lo 1} \lo {u \lo 1 q})   h \lo {\ell} (\xi \hi {j \lo 2} \lo {u \lo 2 r})-E( {h \lo k} (\xi \hi {j \lo 1} \lo q)) E({h \lo \ell} (\xi \hi {j \lo 2} \lo r)).
\end{align*}
Next, observe that the summands  of $A \lo 1$ have expectation $0$ and are bounded in absolute value by $2$.   Therefore, by Hoeffding's inequality, we have  
% \begin{align*}
$
 P (|A \lo {1}| \geq \epsilon) \leq 2 \exp (-\frac{n \epsilon \hi 2}{128})
 $
% \end{align*}
for any $\epsilon >0$.
Moreover, the term $A \lo 2$ can be written as $ \frac{n-1}{n}A \lo {21} + A \lo {22}$, where
\begin{align*}
 A \lo {21} = \frac{1}{n (n-1)}  \sum \lo {u \lo 1 \ne u \lo 2} \hi n   h \lo k (\xi \hi {j \lo 1} \lo {u \lo 1 q}) h \lo \ell (\xi \hi {j \lo 2} \lo {u \lo 2 r})- E(h \lo k (\xi \hi {j \lo 1} \lo {u \lo 1 q}) h \lo \ell (\xi \hi {j \lo 2} \lo {u \lo 2 r}),  
\end{align*}
is a $U$-statistic and
$ A \lo {22}=\frac{1}{n \hi 2} \sum  \lo {u \lo 1=1} \hi {n}     h \lo k (\xi \hi {j \lo 1} \lo {u \lo 1 q}) h \lo \ell (\xi \hi {j \lo 2} \lo {u \lo 1 r})- E(h \lo k (\xi \hi {j \lo 1} \lo {u \lo 1 q}) h \lo \ell (\xi \hi {j \lo 2} \lo {u \lo 1 r}))$.
Consequently, by Hoeffding's inequality for $U$-statistics \citep{hoeffding1963probability}
% \begin{align*}
$
 P (|A \lo {21}| \geq \epsilon) \leq 2 \exp (-\frac{n \epsilon \hi 2}{128})
 $
% \end{align*}
  for any $\epsilon >0$, 
 and it is easy to see  (due to the additional factor $1/n$) that $A \lo {22}$  satisfies an even stronger concentration inequality.  Therefore, it follows that for any $\epsilon >0$
\begin{align}\label{hoef}
P \left( | W \lo {k \ell, rq} \hi {j \lo 1, j \lo 2} | \geq \epsilon \right) \leq 2 \exp \left ( - C \lo 1 {n \epsilon \hi 2} \right),
\end{align}
for some constant $C \lo 1>0$.  Thus,  the union bound over the $(n \hi i m \lo n k \lo n) \hi 2$ indices and the choice of $\epsilon =\frac{\delta }{n \hi i m \lo n k \lo n}$ in \eqref{hoef} 
 yields \eqref{boundtruesigma}. Finally, the assertion \eqref{sampleboundedeigencondition} follows from relation \eqref{weyls} at the beginning of the proof. \eop

\noindent
The next Lemma guarantees that the matrices defined in \eqref{samplematrices}  satisfy the irrepresentable  condition in Assumption \ref{irrepr} with  high probability.

\begin{lemma}\label{sampleirrepr}  If Assumption \ref{boundedeigen} and 
\ref{irrepr} are satisfied for some $0< \eta \leq 1$, then 
\begin{align*}
&P \Big( \max \lo { j \notin \sf N \hi i \lo n } \| \Sigma \hi n  \lo {\xi \hi j \sf N \hi i \lo n} (\Sigma \hi n \lo { \sf N \hi i \lo n \sf N \hi i \lo n}) \hi {-1} \| \lo {F}  \geq \frac{1-\frac{\eta}{2}}{  \sqrt{n \hi i}} \Big)  \lesssim \exp \Big(- C \lo 1 \frac{n}{((n \hi i) \hi {5/4} m \lo n k \lo n) \hi 2}+ 2 \log (p m \lo n k \lo n) \Big),
\end{align*}
where $C \lo 1$ is a positive constant that depends only on  $C \lo {\min}$ and $\eta$.
\end{lemma}

\proof
First, we decompose 
$$
%\begin{split}
 \max \lo { j \notin \sf N \hi i \lo n }  \| \Sigma \hi n  \lo {\xi \hi j \sf N \hi i \lo n} (\Sigma \hi n \lo { \sf N \hi i \lo n \sf N \hi i \lo n}) \hi {-1}   \| \lo {F} 
  \leq  \max \lo { j \notin \sf N \hi i \lo n }  \| \Sigma \hi n \lo {\xi \hi j \sf N \hi i \lo n} (\Sigma \hi n  \lo { \sf N \hi i \lo n \sf N \hi i \lo n}) \hi {-1} - \Sigma \hi * \lo {\xi \hi j \sf N \hi i \lo n} (\Sigma \hi * \lo { \sf N \hi i \lo n \sf N \hi i \lo n}) \hi {-1}     \| \lo {F}+  \max \lo { j \notin \sf N \hi i \lo n } \| \Sigma \hi * \lo {\xi \hi j \sf N \hi i \lo n} (\Sigma \hi * \lo { \sf N \hi i \lo n \sf N \hi i \lo n}) \hi {-1}  \| \lo {F}.
% \end{split}
$$
 By Assumption \ref{irrepr} we have
 $ \max \lo { j \notin \sf N \hi i \lo n } \| \Sigma \hi * \lo {\xi \hi j \sf N \hi i \lo n} (\Sigma \hi * \lo { \sf N \hi i \lo n \sf N \hi i \lo n}) \hi {-1}  \| \lo F \leq \frac{1- \eta}{\sqrt{n \hi i}}$ and
   it suffices to consider
\begin{align*}
&P \Big( \max \lo { j \notin \sf N \hi i \lo n } \| \Sigma \hi n  \lo {\xi \hi j \sf N \hi i \lo n} (\Sigma \hi n \lo { \sf N \hi i \lo n \sf N \hi i \lo n}) \hi {-1} -\Sigma \hi *  \lo {\xi \hi j \sf N \hi i \lo n} (\Sigma \hi * \lo { \sf N \hi i \lo n \sf N \hi i \lo n}) \hi {-1}  \| \lo {F}  \geq \frac{\eta}{2  \sqrt{n \hi i}} \Big).
\end{align*}
For this purpose we use the decomposition $ \Sigma \hi n  \lo {\xi \hi j \sf N \hi i \lo n} (\Sigma \hi n \lo { \sf N \hi i \lo n \sf N \hi i \lo n}) \hi {-1} -\Sigma \hi *  \lo {\xi \hi j \sf N \hi i \lo n} (\Sigma \hi * \lo { \sf N \hi i \lo n \sf N \hi i \lo n}) \hi {-1} =T \lo 1 \hi j + T \lo 2 \hi j + T \lo 3\hi j $ where 
\begin{align*}
T \lo 1 \hi j &=\Sigma \hi * \lo {\xi \hi j \sf N \hi i \lo n} \left( (\Sigma \hi n \lo { \sf N \hi i \lo n \sf N \hi i \lo n}) \hi {-1}- (\Sigma \hi * \lo { \sf N \hi i \lo n \sf N \hi i \lo n}) \hi {-1}\right)~, ~~
T \lo 2 \hi j =\left( \Sigma \hi n \lo {\xi \hi j \sf N \hi i \lo n} -\Sigma \hi * \lo {\xi \hi j \sf N \hi i \lo n}   \right) (\Sigma \hi * \lo { \sf N \hi i \lo n \sf N \hi i \lo n}) \hi {-1},\\
T \lo 3 \hi j &=\left( \Sigma \hi n \lo {\xi \hi j \sf N \hi i \lo n} -\Sigma \hi * \lo {\xi \hi j \sf N \hi i \lo n}   \right)\left( (\Sigma \hi n \lo { \sf N \hi i \lo n \sf N \hi i \lo n}) \hi {-1}- (\Sigma \hi * \lo { \sf N \hi i \lo n \sf N \hi i \lo n}) \hi {-1}\right), 
\end{align*}
and control
  the probabilities 
$ P  \big ( \max \lo { j \notin \sf N \hi i \lo n } \| T \lo h \hi j \| \lo {F} \geq \frac{\eta}{6 \sqrt{n \hi i}} \big )$
separately.

\noindent
(a) For the first term $T \lo 1 \hi j$,  we use the identity   $A^{-1}-B^{-1}=A^{-1}(B-A)B^{-1}$
and obtain from  Assumption \ref{irrepr}
\begin{align*}
\max \lo { j \notin \sf N \hi i \lo n } \| T \lo 1 \hi j \| \lo {F} &\leq \max \lo { j \notin \sf N \hi i \lo n } \| \Sigma \hi * \lo {\xi \hi j \sf N \hi i \lo n} (\Sigma \hi * \lo { \sf N \hi i \lo n \sf N \hi i \lo n}) \inv \| \lo F \| (\Sigma \hi n  \lo { \sf N \hi i \lo n \sf N \hi i \lo n}-\Sigma \hi * \lo { \sf N \hi i \lo n \sf N \hi i \lo n}) (\Sigma \hi n \lo { \sf N \hi i \lo n \sf N \hi i \lo n}) \inv \| \lo F \\
& \leq \frac{(1-\eta)}{\sqrt{n \hi i}} \| \Sigma \hi n  \lo { \sf N \hi i \lo n \sf N \hi i \lo n}-\Sigma \hi * \lo { \sf N \hi i \lo n \sf N \hi i \lo n} \| \lo F \| (\Sigma \hi n \lo { \sf N \hi i \lo n \sf N \hi i \lo n}) \inv \| \lo 2.
\end{align*}
 %  and the property $\| A B \| \lo F \leq \| A \| \lo 2 \|B \| \lo F$. 
 Thus, defining the event 
%\begin{align}\label{conditionalevent}
$\mathcal{T}= \{ \|(\Sigma \hi n \lo { \sf N \hi i \lo n \sf N \hi i \lo n}) \inv \| \lo {2} \leq \frac{2}{C \lo {\min}} \} $
we obtain 
%\begin{align*}
%P \Big ( \max \lo { j \notin \sf N \hi i \lo n } \| T \lo 1 \hi j \| \lo {F} \geq \frac{\eta}{6  \sqrt{n \hi i}} \Big ) &\leq P \Big(  \| \Sigma \hi n  \lo { \sf N \hi i \lo n \sf N \hi i \lo n}-\Sigma \hi * \lo { \sf N \hi i \lo n \sf N \hi i \lo n}\| \lo F \|(\Sigma \hi n \lo { \sf N \hi i \lo n \sf N \hi i \lo n}) \inv \| \lo {2} \geq \frac{\eta}{6 (1- \eta)} \Big ).
%\end{align*}
{
\begin{equation}\label{T1}
\begin{split}
%&P \left(   \| (\Sigma \hi n  \lo { \sf N \hi i \lo n \sf N \hi i \lo n}-\Sigma \hi * \lo { \sf N \hi i \lo n \sf N \hi i \lo n})\| \lo F \|(\Sigma \hi n \lo { \sf N \hi i \lo n \sf N \hi i \lo n}) \inv \| \lo {2} \geq \frac{\eta}{6 (1- \eta)} \right)\\
P \Big ( \max \lo { j \notin \sf N \hi i \lo n } \| T \lo 1 \hi j \| \lo {F} \geq \frac{\eta}{6  \sqrt{n \hi i}} \Big ) 
&\leq P \Big(  \| \Sigma \hi n  \lo { \sf N \hi i \lo n \sf N \hi i \lo n}-\Sigma \hi * \lo { \sf N \hi i \lo n \sf N \hi i \lo n}\| \lo F \|(\Sigma \hi n \lo { \sf N \hi i \lo n \sf N \hi i \lo n}) \inv \| \lo {2} \geq \frac{\eta}{6 (1- \eta)} \Big )
\\
%& \leq  P \Big(   \| \Sigma \hi n  \lo { \sf N \hi i \lo n \sf N \hi i \lo n}-\Sigma \hi * \lo { \sf N \hi i \lo n \sf N \hi i \lo n}\| \lo F \|(\Sigma \hi n \lo { \sf N \hi i \lo n \sf N \hi i \lo n}) \inv \| \lo {2} \geq \frac{\eta}{6 (1- \eta)} \; \mbox{and} \;\mathcal{T}  \Big) +  P \left( \mathcal{T} \hi {\complement} \right). \\ 
&  \leq  P \Big( \| \Sigma  \hi n \lo { \sf N \hi i \lo n \sf N \hi i \lo n}-\Sigma \hi * \lo { \sf N \hi i \lo n \sf N \hi i \lo n}\| \lo F \geq \frac{\eta C \lo {\min}}{12 (1- \eta)} \Big)  +  P \left( \mathcal{T} \hi {\complement} \right) \\
& \leq 4 \exp \Big(- C \lo 1 \frac{n }{(n \hi i m \lo n k \lo n) \hi 2} + 2 \log (n \hi i m \lo n k \lo n)\Big) , 
\end{split}
\end{equation}
where  we used  Lemma \ref{sampleboundedeigen} with $\delta =\frac{\eta C \lo {\min}}{12 (1- \eta)}$
and $\delta=\frac{C \lo {\min}}{2}$ for the last inequality.}

\noindent
(b) For the second term $T \lo 2 \hi j$, we have
\begin{align*}
\max \lo { j \notin \sf N \hi i \lo n } \| T \lo 2 \hi j \| \lo {F}& \leq \|  (\Sigma \hi * \lo { \sf N \hi i \lo n \sf N \hi i \lo n}) \hi {-1} \| \lo 2 \max \lo { j \notin \sf N \hi i \lo n } \|  \Sigma \hi n \lo {\xi \hi j \sf N \hi i \lo n} -\Sigma \hi * \lo {\xi \hi j \sf N \hi i \lo n}    \| \lo F \leq  C \lo {\min} \inv \max \lo { j \notin \sf N \hi i \lo n } \|  \Sigma \hi n \lo {\xi \hi j \sf N \hi i \lo n} -\Sigma \hi * \lo {\xi \hi j \sf N \hi i \lo n}    \| \lo F,
\end{align*}
where we used Assumption \ref{boundedeigen} in the second inequality.  Thus,  
\begin{align*}
P \Big(\max \lo { j \notin \sf N \hi i \lo n } \| T \lo 2  \hi j\| \lo {F} \geq \frac{\eta}{6 \sqrt{n \hi i}} \Big) 
&\leq P \Big( \max \lo { j \notin \sf N \hi i \lo n } \|  \Sigma \hi n \lo {\xi \hi j \sf N \hi i \lo n} -\Sigma \hi * \lo {\xi \hi j \sf N \hi i \lo n}    \| \lo F \geq \frac{\eta C \lo {\min}}{6 \sqrt{n \hi i}} \Big). 
\end{align*}
Now, using similar arguments as in the proof of \eqref{boundtruesigma} in Lemma \ref{sampleboundedeigen}, we can show
\begin{align}\label{bound2}
P \left(\|  \Sigma \hi n \lo {\xi \hi j \sf N \hi i \lo n} -\Sigma \hi * \lo {\xi \hi j \sf N \hi i \lo n}    \| \lo F \geq \delta \right) \leq 2 \exp \Big(- C \lo 1 \frac{n \delta \hi 2}{n \hi i m \lo n \hi 2 k \lo n \hi 2} +  \log (n \hi i m \lo n \hi 2 k \lo n \hi 2)\Big),
\end{align}
for some positive constant $C \lo 1>0$.
This bound with $\delta=\frac{\eta C \lo {\min}}{6 \sqrt{n \hi i} }$ and the union bound yield
\begin{align}\label{T2}
P \Big(\max \lo { j \notin \sf N \hi i \lo n } \| T \lo 2 \hi j\| \lo {F} \geq \frac{\eta}{6 \sqrt{n \hi i}} \Big) 
% &\leq P \Big( \max \lo { j \notin \sf N \hi i \lo n } \|  \Sigma \hi n \lo {\xi \hi j \sf N \hi i \lo n} -\Sigma \hi * \lo {\xi \hi j \sf N \hi i \lo n}    \| \lo F \geq \frac{\eta C \lo {\min}}{6 \sqrt{n \hi i}} \Big) \nonumber \\
&\leq (p- n \hi i) P \Big( \|  \Sigma \hi n \lo {\xi \hi j \sf N \hi i \lo n} -\Sigma \hi * \lo {\xi \hi j \sf N \hi i \lo n}    \| \lo F \geq \frac{\eta C \lo {\min}}{6 \sqrt{n \hi i}} \Big) \nonumber \\
& \leq  2 \exp \Big(- C \lo 1 \frac{n}{(n \hi i m \lo n k \lo n) \hi 2} +  \log (n \hi i m \lo n \hi 2 k \lo n \hi 2)+ \log (p-n \hi i)\Big).
\end{align}
%where we used the fact that  $ \log (p-n \hi i) \leq  \log ((p-n \hi i)m \lo n \hi 2 k \lo n \hi 2)$.

\noindent
(c) For the third term $T \lo 3 \hi j$, we have
\begin{align}\label{t3fist}
P \Big(\max \lo { j \notin \sf N \hi i \lo n } \| T \lo 3 \hi j \| \lo {F}  \geq \frac{\eta}{6 \sqrt{n \hi i} } \Big)  & 
\leq  P \Big( \max \lo { j \notin \sf N \hi i \lo n } \|  \Sigma \hi n \lo {\xi \hi j \sf N \hi i \lo n} -\Sigma \hi * \lo {\xi \hi j \sf N \hi i \lo n}    \| \lo F \geq \sqrt{ \frac{\eta }{6 \sqrt{n \hi i}}} \Big) \nonumber \\
& \quad + P \Big( \|( \Sigma \hi n \lo { \sf N \hi i \lo n \sf N \hi i \lo n}) \hi {-1}- (\Sigma \hi * \lo { \sf N \hi i \lo n \sf N \hi i \lo n}) \hi {-1} \| \lo F \geq  \sqrt{ \frac{\eta }{6 \sqrt{n \hi i}}}  \Big).
\end{align}
Using \eqref{bound2} with $\delta= \sqrt{ \frac{\eta }{6 \sqrt{n \hi i}}} $ we obtain for the first term on the right-hand side of \eqref{t3fist}
\begin{align}\label{T3a}
 P \Big( \max \lo { j \notin \sf N \hi i \lo n } \|  \Sigma \hi n \lo {\xi \hi j \sf N \hi i \lo n} -\Sigma \hi * \lo {\xi \hi j \sf N \hi i \lo n}    \| \lo F \geq \sqrt{ \frac{\eta }{6 \sqrt{n \hi i}}} \Big)   \lesssim  \exp \Big(- C \lo 1 \frac{n}{(n \hi i) \hi {3/2} m \lo n \hi 2 k \lo n \hi 2} +  \log (n \hi i m \lo n \hi 2 k \lo n \hi 2)+ \log ((p-n \hi i)m \lo n \hi 2 k \lo n \hi 2)\Big).
\end{align}
To derive a bound for the second term in \eqref{t3fist} we apply the same arguments as used for the term $T \lo 1 \hi j$   
\begin{align}\label{T3eq2}
& P \Big ( \|( \Sigma \hi n \lo { \sf N \hi i \lo n \sf N \hi i \lo n}) \hi {-1}- (\Sigma \hi * \lo { \sf N \hi i \lo n \sf N \hi i \lo n}) \hi {-1} \| \lo F \geq  \sqrt{ \frac{\eta }{6 \sqrt{n \hi i}}}  \Big) \nonumber \\
 & \lesssim  \exp \Big (- C \lo 1 \frac{n }{(n \hi i) \hi {5/2} m \lo n \hi 2 k \lo n \hi 2}  2 \log (n \hi i m \lo n k \lo n)\Big )+  \exp \Big (- C \lo 1 \frac{n }{(n \hi i  m \lo n   k \lo n) \hi 2} + 2 \log (n \hi i m \lo n k \lo n)\Big ).
 \end{align}
Thus, from \eqref{t3fist}, \eqref{T3a} and \eqref{T3eq2} we obtain
\begin{align}\label{T3}
P \Big(\max \lo { j \notin \sf N \hi i \lo n } \| T \lo 3 \hi j \| \lo {F}  \geq \frac{\eta}{6 \sqrt{n \hi i}} \Big)  &  \lesssim   \exp \Big(- C \lo 1 \frac{n}{(n \hi i) \hi {3/2}  (m \lo n k \lo n) \hi 2} +  \log (n \hi i m \lo n \hi 2 k \lo n \hi 2)+ \log ((p-n \hi i)m \lo n \hi 2 k \lo n \hi 2)\Big) \nonumber \\
& \quad+ 2 \exp \Big(- C \lo 1 \frac{n }{(n \hi i) \hi {5/2}  m \lo n  \hi 2 k \lo n \hi 2} + 2 \log (n \hi i m \lo n k \lo n)\Big)\\
& \quad+ 2 \exp \Big(- C \lo 1 \frac{n }{(n \hi i  m \lo n   k \lo n) \hi 2} + 2 \log (n \hi i m \lo n k \lo n)\Big).  \nonumber
\end{align}
Putting together \eqref{T1}, \eqref{T2} and \eqref{T3} and using the fact  $\log n \hi i \leq \log (p-n \hi i) \leq \log p$ (since $n \hi i \leq p$) we conclude 
\begin{align*}
&P \Big( \max \lo { j \notin \sf N \hi i \lo n } \| \Sigma \hi n  \lo {\xi \hi j \sf N \hi i \lo n} (\Sigma \hi n \lo { \sf N \hi i \lo n \sf N \hi i \lo n}) \hi {-1}  \| \lo F \geq \frac{1-\frac{\eta}{2}}{\sqrt{n \hi i}} \Big)   \leq 2 \exp \Big(- C \lo 1 \frac{n}{(n \hi i) \hi {5/2} m \lo n \hi 2 k \lo n \hi 2}+ 2 \log (p m \lo n k \lo n) \Big),
\end{align*}
for some positive constant $C \lo 1$ that depends on $C \lo {\min}$ and $\eta$.  This completes the proof. \eop
%%%%%%%%%%%%%%%%PART III%%%%%%%%%%%%%%%%%%%%
\subsection{ Proof of Theorem \ref{theorem2section4}}\label{sec83}

{We begin  establishing the model selection consistency given that Assumptions \ref{boundedeigen} and \ref{irrepr} are satisfied by the sample matrices defined in \eqref{samplematrices}.} 
 In particular we define the event
 \begin{align}\label{conditionalevent}
\mathcal{N}=\{  \Sigma \hi n \lo {\sf N \lo n \hi i \sf N \lo n \hi i}, \Sigma \hi n \lo {{\xi \hi j}  \sf N \lo n \hi i  }  \; \text{satisfy Assumptions} \;  \ref{boundedeigen} \; \mathrm{and} \; \ref{irrepr}   \}
\end{align}
and  state the following  result, which is the essential  step in the proof of Theorem \ref{theorem2section4} and will be 
proved in Section \ref{sec 84} below.

\begin{proposition}\label{sampletheorem} {If the  assumptions of Theorem \ref{theorem2section4} are satisfied}. 
  Then, 
    \begin{equation}
    \label{det5}
\begin{split}
&{ P ( \hat {\sf N} \hi i \lo n \ne {\sf N } \hi i \lo n 
\cap {\cal N} ) 
  \lesssim \exp \Big( - C \lo 1 \frac{n \hi {1- \alpha (2+3\beta)} (\lambda \lo n \tsum \lo {j \in  \sf N \hi i \lo n } \|  B \hi {*ij} \lo {m \lo n k \lo n} \| \lo F) \hi 2  }{n \hi i m \lo n \hi 2  k \lo n \hi {4} }+2 \log (n \hi i m \lo n k \lo n) \Big),}
\end{split}
  \end{equation}
  where $C \lo 1$ is a positive constant.
\end{proposition}

\noindent 
We have
 \begin{align*}
P( \hat {\sf N \hi i \lo n} \ne \sf N \hi i \lo n )  \leq  P( \hat {\sf N \hi i \lo n} \ne \sf N \hi i \lo n \text{ and } \mathcal{N} ) + P (\mathcal{N} \hi {\complement}) ~,
 \end{align*}
 where the first probability on the right hand side can be estimated by \eqref{det5}.
  Moreover, by Lemmas \ref{sampleboundedeigen} and \ref{sampleirrepr},
  \begin{align*}
 P (\mathcal{N} \hi {\complement})&\lesssim \exp \Big(- C \lo 1 \frac{n}{((n \hi i) \hi {5/4} m \lo n k \lo n) \hi 2}+ 2 \log (p m \lo n k \lo n) \Big)
 %\\ & \quad 
 + \exp \Big(- C \lo 1 \frac{n}{(n \hi i) \hi {2} m \lo n k \lo n) \hi 2}+ 2 \log (p m \lo n k \lo n) \Big),
 \end{align*}
 and this proves Theorem \ref{theorem2section4}.
 
\subsection{Proof of Proposition \ref{sampletheorem}}
\label{sec 84}
%%% proof of main theorem %%%%%
%\subsection{Proof of Theorem \ref{theorem2section4}}
We follow a similar strategy as in \cite{bach2008consistency} and  \cite{lee2016variable}, who showed consistency of the group lasso in a reproducing kernel Hilbert space framework.  First, we consider the following alternative  form of the group lasso problem \eqref{finalpenfunction}
  \begin{align}\label{penfunction2} 
\widehat {PL} \lo i (B \hi i, \hat \xi )=\frac{1}{2n} \| \hat \xi \hi i -  \mathbf{\tilde H \lo n} \trans (\hat \xi \hi {-i}) B \hi {i} \| \lo {F} \hi 2 +  \frac{\lambda \lo n}{2}  (\tsum \lo {j \ne i} \hi p \| B \hi {ij} \| \lo F ) \hi 2.
 \end{align}
Because the function $x \to x \hi 2$, $x \geq 0$ is  monotone, problem \eqref{penfunction2} leads to the same regularisation paths as problem \eqref{finalpenfunction} \citep[see][, page 1187 for more details]{bach2008consistency}.  To derive the Karush-Kuhn-Tucker (KKT) conditions, we recall the notations \eqref{hol2}, \eqref{samplemat1} and   define the matrices
   \begin{align}
  \label{hd8}
   \hat {\Sigma} \hi n \lo { \sf N \hi i \lo n {\xi \hi j}} & = \frac {1}{n} \mathbf{\tilde H \lo n} ( \hat \xi \hi { \sf N \hi i \lo n} )  \tilde H \lo n  (\hat \xi \hi { j} )    {~ \in \R \hi {n \hi i k \lo n m \lo n \times k \lo n m \lo n} }, \\ 
    \label{hd9a}
   \hat \Sigma \hi n  \lo { { \xi \hi j  \xi \hi i}} & = \frac {1}{n} \tilde H \lo n  \trans ( \hat \xi \hi { j} ) \hat \xi \hi { i}     {~ \in \R \hi {k \lo n m \lo n \times m \lo n} }
   \end{align} 
    when  $j \ne i$ and 
    \begin{align}
 \label{hd10}
 \hat {\Sigma} \hi n  \lo { \sf N \hi i \lo n {\xi \hi i}}= \frac {1}{n} \mathbf{\tilde{ H \lo n  }  ( \hat \xi \hi { \sf N \hi i \lo n} )} \hat \xi \hi { i}     {~ \in \R \hi {n \hi i  k \lo n m \lo n \times  m \lo n} },
 \end{align}  
  where the matrices $ \tilde H \lo n  ( \hat \xi \hi { j} )$ and $\mathbf{\tilde{ H \lo n  }  ( \hat \xi \hi { \sf N \hi i \lo n} )}$ have been defined in \eqref{estimatedmat2} and \eqref{estimatedmat1} respectively.   We also denote by ${\Sigma} \hi n \lo { \sf N \hi i \lo n {\xi \hi j}}$, $  \Sigma \hi n  \lo { { \xi \hi j  \xi \hi i}} $, $  {\Sigma} \hi n \lo { \sf N \hi i \lo n {\xi \hi i}}$  the versions of \eqref{hd8}, \eqref{hd9a}, \eqref{hd10} that use the true scores $\xi \hi i \lo {ur}$ instead of the estimated  $\hat \xi \hi i \lo {ur}$ (see also equation \eqref{estimatedmat3}).

%%%%%%%%%%%%%%%
%KKT Conditions
%%%%%%%%%%%%%%%%%
\begin{lemma}\label{lemmakkt} (KKT conditions) A matrix $ B \hi i =( B \hi {ij} , j \in \sf V \setminus \{i\}) \in \R \hi {(p-1) k \lo n m \lo n \times m \lo n}$ with support $ {\sf N \hi i \lo n} $ is optimal for problem $\eqref{penfunction2}$ if and only if
\begin{subequations}
\begin{align}
& (\hat \Sigma \hi n  \lo { {\sf N \hi i \lo n} {\sf N \hi i  \lo n}} + \lambda \lo n \hat D \lo { \sf N \hi i \lo n})  B  \hi { {\sf N \hi i \lo n} }-\hat {\Sigma} \hi n  \lo   {\sf N  \hi i \lo n \xi \hi i}=0  , \quad \mbox{for all } j \in {\sf N \hi i  \lo n}, \label{kktfirst}\\
 &\|\hat  \Sigma \hi n  \lo { \xi \hi j {\sf N \hi i \lo n}} B  \hi { {\sf N  \hi i \lo n}} - \hat {\Sigma} \hi n  \lo  {\xi \hi j {\xi \hi i}} \| \lo {F} \leq \lambda \lo n \tsum \lo {j \ne i} \hi {p} \|  B \hi {ij}    \| \lo F,   \quad \mbox{for all } j \notin {{\sf N \lo n \hi i}  }    \label{kktsecond}
\end{align}
\end{subequations}
where $\hat \Sigma \hi n  \lo { {\sf N \hi i \lo n} {\sf N \hi i  \lo n}} $ is defined in \eqref{estimatedmat3}, $B  \hi { {\sf N  \hi i \lo n}}=(B \hi {ij}, j \in {\sf N \hi i \lo n}) \in \R^{n^ik_nm_n \times m_n}$, $B =(\beta \hi {ij} \lo {qrk}: 1 \leq q, r \leq m \lo n, 1 \leq k \leq k \lo n)$ and
\begin{align*} 
%\label{hd12}
\hat D \lo {\sf N \hi i \lo n}
  = \mathrm{diag}\big ( ( \hat D \lo {\sf N \hi i \lo n})\lo {jj} \colon j  \in {\hat {\sf N \hi i \lo n} } ) \big ) 
\end{align*}
 is a  block diagonal  matrix with $n \hi i$ elements
 $ 
 (\hat D \lo { {\sf N \hi i \lo n}}) \lo {jj}=\frac{\tsum \lo {\ell \ne i} \hi {p} \| \hat B \hi {i\ell}   \| \lo F}{\| \hat B \hi {ij} \| \lo F  } I \lo {k \lo n m \lo n} \in \R \hi {k \lo n m \lo n \times k \lo n m \lo n}.
 $ 
\end{lemma}
{The idea of the proof is to first construct an estimator $\hat B \hi {\sf N \hi i \lo n} \lo n$ by minimizing the following restricted problem given the true support $ \sf N \hi i \lo n$.  That is,  
\begin{align}\label{restricted}
\hat B \hi {\sf N \hi i \lo n} \lo n&=\mathrm{argmin} \big \{ \widehat {PL} \lo { \sf N \hi i \lo n } (  B  , \hat \xi ):  B  \in \R \hi {n^i k \lo n m \lo n \times m \lo n} \big \}, 
 \end{align}
where
\begin{align}\label{restricted2}
\widehat {PL} \lo { \sf N \hi i \lo n } (  B  , \hat \xi )= \frac{1}{2n} \| \hat \xi \hi i -  \mathbf{\tilde H\lo n \trans (\hat \xi \hi {\sf N \hi i \lo n})}  B  \| \lo {F} \hi 2 +  \frac{\lambda \lo n}{2} \Big (\tsum \lo {j \in \sf N \hi i \lo n } \hi p \| B \hi {ij} \| \lo F \Big ) \hi 2,
\end{align}
{(note that $\widehat {PL} \lo { \sf N \hi i \lo n } (  B  , \hat \xi )$ corresponds to the function \eqref{penfunction2}, where we  put   $B \hi {ij}=0$ whenever $j \notin \sf N \hi i$).}
and to show  that the minimizer in  \eqref{restricted} is ``close" to the true matrix $B \hi {* \sf N \hi i \lo n} \lo n$ defined in \eqref{hd2}.  To achieve this
we use similar arguments as in    \cite{bach2008consistency} 
and construct another auxiliary estimator $\tilde B \hi {\sf N \hi i \lo n} \lo n$ that minimizes the restricted penalized function, where the group lasso penalty in \eqref{restricted} is replaced by an $\ell \lo 2$-type penalty. 
More precisely, 
 $\tilde B \hi {\sf N \hi i \lo n} \lo n$ is defined by 
 \begin{align}
 \label{restrictedridge}
&\tilde B \hi {\sf N \hi i \lo n} \lo n 
=\mathrm{argmin} \big \{ \widetilde {PL} \lo { \sf N \hi i \lo n } (  B  , \hat \xi ):  B  \in \R \hi {n^i k \lo n m \lo n \times m \lo n} \big  \},
\end{align}
where 
\begin{align*}
\widetilde {PL} \lo { \sf N \hi i \lo n } (  B  , \hat \xi )=  \frac{1}{2n} \| \hat \xi \hi i -   \mathbf{ \tilde H \lo n \trans (\hat \xi \hi {-i}) } B  \| \lo {F} \hi 2 +  \frac{\lambda \lo n}{2} \Big (\tsum \lo {\ell \in \sf N \hi {i} \lo {n}}  \| B \hi {*i\ell}  \lo {m \lo n k \lo n}  \| \lo F \Big ) \Big  (\tsum \lo {j \in \sf N \hi {i} \lo {n}  }  \frac{\| B \hi {ij} \| \hi 2 \lo F}{\| B \hi {*ij}  \lo {m \lo n k \lo n}  \| \lo F} \Big ). 
 \end{align*} 
  We now proceed in the following steps:
 \begin{itemize}
\item[(1)]  {In Proposition \ref{ridgeconsi} we show that the distance $\| \tilde B \hi {\sf N \hi i \lo n} \lo n- B \hi {*\sf N \hi i \lo n} \lo n\| \lo {F}$ is small with high probability.
\item[(2)] In Proposition \ref{bhatconsi} we show  that $\hat B \hi {\sf N \hi i \lo n} \lo n$ is close to $\tilde B \hi {\sf N \hi i \lo n} \lo n$ with high probability.
\item[(3)]  In Proposition \ref{propconsistency} we use this result to derive  a concentration bound for $\| \hat B \hi {\sf N \hi i \lo n} \lo n- B \hi {*\sf N \hi i \lo n} \| \lo {F}$.  
\item[(4)] We then construct the oracle minimiser $(\hat B \hi {\sf N \hi i \lo n}, \mathbf{0})$, where $\hat B \hi {\sf N \hi i \lo n}$ is the minimiser of \eqref{restricted} and $\mathbf{0}$ consists of $(p-1-n \hi i)$ zero $k \lo n m \lo n \times m \lo n$ matrices. 
\item[(5)]
Finally, in Proposition  \ref{bhatoptim} we show that the oracle minimiser is optimal for the restricted problem \eqref{restricted} given the true support $\sf N \hi i \lo n$; that is, it satisfies \eqref{kktsecond}}. 
 \end{itemize}
The minimisation problem \eqref{restricted} is convex; however, for $p >n$, it need not to be strictly convex,  so that there may not be a unique solution. Nevertheless, the next lemma shows that the matrix $ \hat \Sigma \hi n \lo {  {\sf N \hi i \lo n}  {\sf N \hi i \lo n}} $ defined in \eqref{estimatedmat3} is strictly positive definite with high probability, and hence the objective function \eqref{restricted} is strictly convex, and thus $\hat B \hi {\sf N \hi i \lo n}$  is the unique optimal solution.  
{\begin{lemma}    There exists a constant $C \lo 1>0$ such that
\begin{align*} 
P \Big( \Lambda \lo {\min} ( \hat \Sigma \hi n \lo {  {\sf N \hi i \lo n}  {\sf N \hi i \lo n}}) \geq \frac{C \lo {\min}}{4} \Big) \gtrsim 1-  \exp \Big(- C \lo 1 \frac{n \hi {1-\alpha(2+3\beta)}}{ (n \hi i) \hi 2 m \lo n \hi 2 k \lo n \hi 4}+2 \log (n \hi i m \lo n k \lo n) \Big).
\end{align*}
\end{lemma}}
\proof   {By Weyl's Lemma, we have
$
\Lambda \lo {\min} (\Sigma \hi n \lo   {  {\sf N \hi i \lo n}  {\sf N \hi i \lo n}} ) \leq \Lambda \lo {\min} (\hat \Sigma \hi n  \lo {  {\sf N \hi i \lo n}  {\sf N \hi i \lo n}} ) + \|  \hat \Sigma \hi n  \lo {  {\sf N \hi i \lo n} {\sf N \hi i \lo n} } -  \Sigma \hi n \lo {  {\sf N  \hi i \lo n} {\sf N \hi i \lo n} }  \| \lo 2,
$  and     we get
\begin{align*}
P \Big(  \Lambda \lo {\min} ( \hat \Sigma \hi n \lo {  {\sf N \hi i \lo n}  {\sf N \hi i \lo n} }) \leq \frac{C \lo {\min}}{4} \mbox{ and } \Lambda \lo {\min} (\Sigma \hi n \lo   {  {\sf N \hi i \lo n}  {\sf N \hi i \lo n}} ) > \frac{C \lo {\min}}{2} \Big)  \leq & P \Big( \|  \hat \Sigma \hi n \lo {  {\sf N  \hi i \lo n}  {\sf N \hi i \lo n} } -  \Sigma \hi n \lo {  {\sf N  \hi i \lo n} {\sf N \hi i \lo n} }  \| \lo 2 \geq \frac{C \lo {\min}}{4} \Big).
\end{align*}
Furthermore, using $\delta \hi 2=\frac{1}{ (n \hi i m \lo n  k \lo n)  \hi 2 }\frac{C \lo {\min}}{4}$ in Theorem \ref{theorem1} with the union bound over the $(n \hi i m \lo n k \lo n) \hi 2$ index pairs of the matrix $\hat \Sigma \hi n  \lo {  {\sf N \hi i  \lo n} {\sf N \hi i  \lo n}} -  \Sigma  \hi n  \lo {  {\sf N \hi i  \lo n} {\sf N \hi i  \lo n}} $, yields for some positive constant $C_1$
\begin{align*} %\label{positivelemmaeq2}
P \Big( \|  \hat \Sigma \hi n  \lo {  {\sf N \hi i  \lo n} {\sf N \hi i  \lo n}}  - \Sigma \hi n  \lo {  {\sf N \hi i  \lo n} {\sf N \hi i  \lo n}}   \| \lo F \geq \frac{C \lo {\min}}{4} \Big) \lesssim \exp \Big(- C \lo 1 \frac{n \hi {1-\alpha(2+3\beta)}}{  (n \hi i) \hi 2 m \lo n \hi 2 k \lo n \hi 4}+2 \log (n \hi i m \lo n k \lo n) \Big).
\end{align*}
The assertion now follows by the same arguments as given in  the proof of Lemma \ref{sampleboundedeigen}.}} 
\eop

%%%%%%%%%%%%%%%%PROPOSITION 1%%%%%%%%%%%%%%%%%%%%%%%%%%%%%
 % Recalling the notations \eqref{restrictedridge} and \eqref{hd2}, we have the following proposition.

 \begin{proposition} \label{ridgeconsi} Suppose Assumptions \ref{as1}-\ref{boundedeigen} hold and the regularization parameter $\lambda \lo n$  satisfies
  \begin{align}\label{as1ridgeconsi1}
  \frac{n \hi i m \lo n \hi {3/2}}{k \lo n \hi d} \lesssim \sqrt{\frac{2}{C \lo {\min}}} \lambda \lo n  \tsum \lo {j \in  \sf N \hi i  \lo n } \| B \hi {*  i j} \lo {m \lo n k \lo n} \| \lo F.
  \end{align}  
  Then,  there exists a constant $c \lo 2 \in (0,1/2) $ such that, for any $\delta>0$ satisfying
   \begin{align}\label{as1ridgeconsi2}
   \frac{2}{C \lo {\min}} \sqrt{n \hi i} \lambda \lo n \tsum \lo {j \in  \sf N \hi i  \lo n } \| B \hi {*  i j} \lo {m \lo n k \lo n} \| \lo F \leq c \lo 2 \delta,
     \end{align}
    we have for the minimizer of \eqref{restrictedridge}
   \begin{align*}
P \big(\| \tilde B \lo n \hi {\sf N \hi i  \lo n}- B \lo n \hi {*\sf N \hi i \lo n} \| \lo {F} \geq \delta \big) \lesssim&  \exp \Big( - C \lo 1 \frac{n \hi {1- \alpha (2+3\beta)} \delta \hi 2}{(n \hi i) \hi 2 m \lo n \hi 2 k \lo n \hi {4}}+2 \log (n \hi i m \lo n k \lo n ) \Big),
 \end{align*}
 where $B \lo n \hi {*\sf N \hi i \lo n}$ is defined in \eqref{hd2} and  the constant $C \lo 1$ satisfies  $0 < \delta \leq C \lo 1$.
\end{proposition}

\proof
{Before we start with the proof we note that  condition \eqref{as1ridgeconsi1} refers to the spline approximation error from including only $k \lo n$ terms and the second condition \eqref{as1ridgeconsi2} represents the bias due to ridge penalisation.}

For the proof we use similar arguments as given in
the proof of  Proposition 2 of \cite{leefaro}. The main change that we need to consider is the approximation error of the additive regression functions by splines.
First, the minimizer $\tilde B \hi {\sf N \hi i} \lo n$  defined in \eqref{restrictedridge}  is of the form
 \begin{align*}
\tilde B \lo n \hi {\sf N \hi i  \lo n}= \big(\hat \Sigma \hi n  \lo {\sf N \hi i \lo n \sf N \hi i \lo n}+ \lambda \lo n D \hi * \lo {\sf N \hi i \lo n}\big ) \hi {-1}\hat \Sigma \hi n  \lo {\sf N \hi i \lo n \xi \hi i}.
 \end{align*}
where $D \hi * \lo {\sf N \hi i \lo n}$ is a block diagonal matrix with $(D \hi * \lo {\sf N \hi i \lo n}) \lo {jj}=\tsum \lo {\ell \ne i} \hi p \| B \hi {*i \ell}  \lo {m \lo n k \lo n}  \| \lo F /\| B \hi {*ij}  \lo {m \lo n k \lo n}  \| \lo F I \lo {k \lo n m \lo n}$, $j \in \sf N \hi i \lo n$ as diagonal blocks, and the matrices $\hat \Sigma \hi n  \lo {\sf N \hi i \lo n \sf N \hi i \lo n}$ and $\hat \Sigma \hi n  \lo {\sf N \hi i \lo n \xi \hi i}$ are defined in \eqref{estimatedmat3} and \eqref{hd8}, respectively.  

A simple calculation shows that  
$$
%\begin{equation}
%\label{hr3} 
\|  \tilde B \hi {\sf N \hi i \lo n} \lo n-  B \hi {*\sf N \hi i \lo n} \lo n  \| \lo {F} \leq T \lo 1 +  T \lo 2 +  T \lo 3,
%\end{equation}
$$
where the terms $T \lo {1}$, $T \lo {2}$ and $T \lo {3}$ are defined by 
   \begin{align*}
 T \lo {1} & = \big  \| (\hat \Sigma \hi n  \lo {\sf N \hi i \lo n \sf N \hi i \lo n}+ \lambda \lo n D \hi * \lo {\sf N \hi i \lo n}) \hi {-1} (\hat \Sigma \hi n  \lo {\sf N \hi i \lo n \xi \hi i} 
 - \Sigma \hi n \lo {\sf N \hi i \lo n \xi \hi i}) \big  \| \lo {F} ,  \\
T \lo {2} & =  \big  \| \{ (\hat \Sigma \hi n  \lo {\sf N \hi i \lo n \sf N \hi i \lo n}+ \lambda \lo n D \hi * \lo {\sf N \hi i \lo n}) \hi {-1}
 -( \Sigma \hi n  \lo {\sf N \hi i \lo n \sf N \hi i \lo n}+ \lambda \lo n D \hi * \lo {\sf N \hi i \lo n}) \hi {-1} \} \Sigma \hi n  \lo {\sf N \hi i \lo n \xi \hi i} \big  \| \lo {F}  , \\
 T \lo {3}& = \big  \|  ( \Sigma \hi n  \lo {\sf N \hi i \lo n \sf N \hi i \lo n}+ \lambda \lo n D \hi * \lo {\sf N \hi i \lo n}) \hi {-1}  \Sigma \hi n  \lo {\sf N \hi i \lo n \xi \hi i}- B \hi {*\sf N \hi i \lo n} \lo n \big  \| \lo {F}.
  \end{align*}
  Thus,  
  $  
  P \big (\| \tilde B \hi {\sf N \hi i \lo n} \lo n- B \hi {*\sf N \hi i \lo n} \lo n \| \lo {F} \geq 3 \delta \big) \leq \tsum \lo {i=1} \hi 3 P (T \lo i \geq {\delta}),
  $
  and it is sufficient to    derive  bounds for the three probabilities corresponding to the random variables
  $T \lo 1$, $T \lo 2$ and $T \lo 3$. Starting with $T \lo 1$
  we have
 \begin{align*}
 T \lo 1 &\leq   
  \| (\hat \Sigma \hi n  \lo {\sf N \hi i \lo n \sf N \hi i \lo n}+ \lambda \lo n D \hi * \lo {\sf N \hi i \lo n}) \hi {-1} \| \lo {2}  \| \hat \Sigma \hi n  \lo {\sf N \hi i \lo n \xi \hi i} - \Sigma \hi n  \lo {\sf N \hi i \lo n \xi \hi i} \| \lo {F}
  % \\  &
  \leq  \frac{2}{C \lo {\min}}  \| \hat \Sigma \hi n  \lo {\sf N \hi i \lo n \xi \hi i} - \Sigma \hi n  \lo {\sf N \hi i \lo n \xi \hi i} \| \lo {F},
  \end{align*}
   {where we use the fact that
  \begin{equation*} %\label{hol4}
  \| \hat \Sigma \lo {\sf N \hi i \lo n \sf N \hi i \lo n} \| \lo 2 \geq \frac{C \lo \min}{2}
  \end{equation*}
   {on the event $\mathcal{N}$  and that $(\hat \Sigma \hi n  \lo {\sf N \hi i \lo n \sf N \hi i \lo n}+ \lambda \lo n D \hi * \lo {\sf N \hi i \lo n}) \hi {-1} \preceq \big ( \hat \Sigma \hi n  \lo {\sf N \hi i \lo n \sf N \hi i \lo n} \big)  \inv$.
Therefore, using Lemma \ref{sampleboundedeigen}, similar arguments as given in the proof of Theorem \ref{theorem1} and applying the union bound over the $n \hi i m \lo n \hi 2 k \lo n$ pairs, we obtain 
 \begin{align}\label{t1final}
P \left(T \lo 1 \geq {\delta} \right)  \leq   P \Big(  \| \hat \Sigma \hi n \lo {\sf N \hi i \lo n \xi \hi i} - \Sigma \hi n  \lo {\sf N \hi i \lo n \xi \hi i} \| \lo {F} \geq \frac{ C \lo {\min} \delta}{2}  \Big)  
\lesssim    \exp \Big( -C \lo 1 \frac{n \hi {1-  \alpha (2 +3 \beta)} \delta \hi 2 }{n \hi i   m \lo n \hi 2 k \lo n \hi {3}} + \log (n \hi i m \lo n \hi 2 k \lo n)  \Big),
  \end{align}
for $0 < \frac{\delta}{n \hi i m \lo n \hi 2 k \lo n} \leq C \lo 1$  with $C \lo 1>0$ depending on $C \lo {\min}$.
%T2 
To derive the bound  for the probability $P \left(T \lo 2 \geq {\delta} \right)$ we use the identity $A \hi {-1} - B \hi {-1}=A \hi {-1}(B-A)B \hi {-1}$ to obtain (on the event ${\cal N}$)
  \begin{align} \label{t2first}
  T \lo 2&
   \leq  \| (\hat \Sigma \hi n  \lo {\sf N \hi i \lo n \sf N \hi i \lo n}+ \lambda \lo n D \hi * \lo {\sf N \hi i  \lo n}) \hi {-1}\| \lo {2} \| ( \Sigma \hi n  \lo {\sf N \hi i \lo n \sf N \hi i \lo n}- \hat \Sigma \hi n  \lo {\sf N \hi i \lo n \sf N \hi i \lo n}) ( \Sigma \hi n  \lo {\sf N \hi i \lo n \sf N \hi i \lo n}+ \lambda \lo n D \hi * \lo {\sf N \hi i \lo n}) \hi {-1} \Sigma \hi n  \lo {\sf N \hi i \lo n \xi \hi i} \| \lo {F} \nonumber \\
  & \leq  \frac{2}{C \lo {\min}} \|  \Sigma \hi n  \lo {\sf N \hi i \lo n \sf N \hi i \lo n}- \hat \Sigma \hi n  \lo {\sf N \hi i \lo n \sf N \hi i \lo n} \| \lo {F}  \| ( \Sigma \hi n  \lo {\sf N \hi i \lo n \sf N \hi i \lo n}+ \lambda \lo n D \hi * \lo {\sf N \hi i \lo n}) \hi {-1} \Sigma \hi n  \lo {\sf N \hi i \lo n \xi \hi i} \| \lo {2}.
    \end{align}
   Recall the relation \eqref{optpredictor}, the notation  $ f  \lo {qr} \hi {ij}=\sum \lo {k=1} \hi {\infty} \beta \lo {qrk} \hi {*ij} {h} \lo {k}  \in \sten {F} \lo {\kappa,\rho}$, and let
    \begin{align}\label{errorspline}
w \hi i \lo {uq}= \tsum \lo {j \in  \sf N \hi i \lo n} \tsum \lo {r =1}  \hi {m \lo n} (f \hi {ij} \lo {qr} (\xi \hi {j} \lo {ur})- \tilde f  {\hi { ij } \lo {nqr}}    ( \xi \hi j \lo {ur})), \quad q=1, \ldots, m \lo n, u=1, \ldots, n,
\end{align}
where $\tilde f  {\hi {ij } \lo {qr}} $
%=\sum \lo {k=1} \hi {k \lo n} \beta \lo {qrk} \hi {*ij} \tilde {h} \lo {nk} $ 
denotes the function from Proposition \ref{splineapproximation}.  Then, we can rewrite relation \eqref{optpredictor} in the form
\begin{align*}
 \xi \hi i= \mathbf{\tilde{H \lo n  \trans} (\xi \hi {\sf N \hi i \lo n})} B \hi {* \sf N \hi i \lo n} \lo n+ w \hi i + \epsilon \hi i \in \R^{n \times m_n},
\end{align*}
where $w \hi i=(w \hi i \lo {uq}) \lo {1 \leq u \leq n, 1 \leq q \leq m \lo n} $,   $\epsilon \hi i=(\epsilon \hi i \lo {uq}) \lo {1 \leq u \leq n, 1 \leq q \leq m \lo n} \in \mathbb{R}^{n \times m_n} $ 
and $\mathbf{\tilde{H^\top \lo n } (\xi \hi {\sf N \hi i \lo n})}$
is defined in \eqref{hol2}.
 Furthermore, by multiplying from the left the above equation  with $\frac{ \mathbf{\tilde{H} \lo n (\xi \hi {\sf N \hi i})}}{n}$ %and the matrix $\tilde H^\top_n (\xi^j)$, 
 we obtain
 \begin{align}\label{sigmazetai}
  \Sigma \hi n  \lo {\sf N \hi i \lo n \xi \hi i}
  & =\Sigma \hi n  \lo {\sf N \hi i \lo n \sf N \hi i \lo n}B \hi {* \sf N \hi i \lo n} \lo n+ \frac{\mathbf{\tilde{H} \lo n  (\xi \hi {\sf N \hi i \lo n})}}{n}  w \hi i + \frac{\mathbf{\tilde{H} \lo n  (\xi \hi {\sf N \hi i \lo n})}}{n}\epsilon \hi i,\\
 \Sigma \hi n  \lo {\xi^j  \xi \hi i}
  &  =   
   \Sigma \hi n  \lo { \xi^j \sf N \hi i \lo n }B \hi {* \sf N \hi i \lo n } \lo n- \frac{\tilde H_n \trans (\xi \hi {j})}{n}  w \hi i - \frac{\tilde H_n \trans (\xi \hi {j})}{n}  \epsilon \hi i. \label{m9}
  \end{align} 
  where the matrix $ \Sigma \hi n  \lo {\sf N \hi i \lo n \xi \hi i}$ is defined in Section \ref{sec 84} and $   \Sigma \hi n  \lo { \xi^j \sf N \hi i \lo n }=  ({\Sigma} \hi n \lo { \sf N \hi i \lo n {\xi \hi j}}  ) \trans$.
 Using this representation and the triangle inequality we get
     \begin{align*} 
\| ( \Sigma \hi n \lo {\sf N \hi i \lo n \sf N \hi i \lo n}+ \lambda \lo n D \hi * \lo {\sf N \hi i \lo n}) \hi {-1} \Sigma \hi n \lo {\sf N \hi i \lo n \xi \hi i}  \| \lo {2}   
  \leq& 
\| ( \Sigma \hi n \lo {\sf N \hi i \lo n \sf N \hi i \lo n}+ \lambda \lo n D \hi * \lo {\sf N \hi i \lo n}) \hi {-1} \Sigma \hi n \lo {\sf N \hi i \lo n \sf N \hi i \lo n}B \hi {* \sf N \hi i \lo n}  \lo n \| \lo 2 + \| ( \Sigma \hi n \lo {\sf N \hi i \lo n \sf N \hi i \lo n}+ \lambda \lo n D \hi * \lo {\sf N \hi i \lo n}) \hi {-1} \frac{\mathbf{\tilde{H} \lo n  (\xi \hi {\sf N \hi i  \lo n})}}{n}  w \hi i  \| \lo {2} \\
&+{ \| ( \Sigma \hi n \lo {\sf N \hi i  \lo n \sf N \hi i  \lo n}+ \lambda \lo n D \hi * \lo {\sf N \hi i  \lo n}) \hi {-1} \frac{\mathbf{\tilde{H} \lo n (\xi \hi {\sf N \hi i  \lo n})}}{n}  \epsilon \hi i  \| \lo {2}}. 
    \end{align*}
As a result   from this and \eqref{t2first}, it follows that  for all $\delta >0$ (on the event ${\cal N}$)
     \begin{align} \label{hr2}
P (T \lo 2 \geq {\delta}) \leq P \Big(T \lo {21} \geq \frac{\delta}{3} \Big)+ P \Big(T \lo {22} \geq \frac{\delta}{3} \Big)+ P \Big(T \lo {23} \geq \frac{\delta}{3} \Big),
\end{align}
 where
    \begin{align*}
    T \lo {21} = &  \frac{2} {C \lo {\min}}   \| \hat \Sigma \hi n  \lo {\sf N \hi i  \lo n \sf N \hi i  \lo n}-  \Sigma \hi n  \lo {\sf N \hi i  \lo n \sf N \hi i  \lo n} \| \lo {F}  \| ( \Sigma \hi n  \lo {\sf N \hi i  \lo n \sf N \hi i  \lo n}+ \lambda \lo n D \hi * \lo {\sf N \hi i  \lo n}) \hi {-1} \Sigma \hi n  \lo {\sf N \hi i  \lo n \sf N \hi i  \lo n}B \hi {* \sf N \hi i  \lo n} \lo n \| \lo 2, \\
     T \lo {22} = &  \frac{2} {C \lo {\min}}  \| \hat \Sigma \hi n  \lo {\sf N \hi i  \lo n \sf N \hi i  \lo n}-  \Sigma  \hi n \lo {\sf N \hi i  \lo n \sf N \hi i  \lo n} \| \lo {F}  \|  ( \Sigma \hi n  \lo {\sf N \hi i  \lo n \sf N \hi i  \lo n}+ \lambda \lo n D \hi * \lo {\sf N \hi i  \lo n}) \hi {-1} \frac{ \mathbf{\tilde H \lo n  (\xi \hi {\sf N \hi i  \lo n})}}{n}  w \hi i  \| \lo {2}, \\
      T \lo {23} = & {  \frac{2} {C \lo {\min}}   \| \hat \Sigma \hi n  \lo {\sf N \hi i  \lo n \sf N \hi i  \lo n}-  \Sigma \hi n  \lo {\sf N \hi i  \lo n \sf N \hi i  \lo n} \| \lo {F}  \|  ( \Sigma \hi n  \lo {\sf N \hi i  \lo n \sf N \hi i  \lo n}+ \lambda \lo n D \hi * \lo {\sf N \hi i  \lo n}) \hi {-1} \frac{ \mathbf{\tilde H \lo n  (\xi \hi {\sf N \hi i  \lo n})}}{n}  \epsilon \hi i  \| \lo {2}}. 
    \end{align*}
Next we derive upper bounds for the probabilities in \eqref{hr2}.  For $T \lo {21} $ observe that
 \begin{align*} 
T \lo {21} \leq&   \frac{2} {C \lo {\min}}   \| \hat \Sigma \hi n \lo {\sf N \hi i   \lo n \sf N \hi i  \lo n}-  \Sigma \hi n  \lo {\sf N \hi i  \lo n \sf N \hi i  \lo n} \| \lo {F}  \|  \| ( \Sigma \hi n  \lo {\sf N \hi i  \lo n \sf N \hi i  \lo n}+ \lambda \lo n D \hi * \lo {\sf N \hi i  \lo n}) \hi {-1} \Sigma \hi n  \lo {\sf N \hi i  \lo n \sf N \hi i  \lo n} \| \lo 2 \| B \hi {* \sf N \hi i  \lo n} \lo n\| \lo F 
%\\
 \leq 
% & 
 \frac{2} {C \lo {\min}}  \| \hat \Sigma \hi n  \lo {\sf N \hi i  \lo n \sf N \hi i  \lo n}-  \Sigma \hi n  \lo {\sf N \hi i  \lo n \sf N \hi i  \lo n} \| \lo {F}, 
 \end{align*}
 where the second inequality uses the fact that $\| B \hi {* \sf N \hi i  \lo n} \lo n\| \lo F < \infty$ by assumption \eqref{splinesasu1} and that the norm $\| ( \Sigma \hi n  \lo {\sf N \hi i  \lo n \sf N \hi i  \lo n}+ \lambda \lo n D \hi * \lo {\sf N \hi i  \lo n}) \hi {-1} \Sigma \hi n  \lo {\sf N \hi i  \lo n \sf N \hi i  \lo n} \| \lo 2$ is bounded by one.
Therefore, it follows from Theorem \ref{theorem1} with $\delta$ replaced by $  \frac{C \lo {\min} \delta}{6 n \hi i m \lo n  k \lo n}$  and  the union bound over the $(n \hi i m \lo n  k \lo n) \hi 2$ pairs that
\begin{equation}\label{t21final}
%\begin{split}
P \Big(T \lo {21} \geq \frac{\delta}{3} \Big) \leq   P \Big( \| \hat \Sigma \hi n  \lo {\sf N \hi i \lo n \sf N \hi i \lo n}-  \Sigma \hi n  \lo {\sf N \hi i \lo n \sf N \hi i \lo n } \| \lo {F}  \geq \frac{ \delta  C \lo {\min} }{6 } \Big)   
\lesssim   \exp \Big(- C \lo 2 \frac{n \hi {1- \alpha (2+3\beta)}   \delta \hi 2}{(n \hi i) \hi 2 m \lo n \hi 2 k \lo n \hi 4}+ 2 \log (n \hi i m \lo n k \lo n) \Big),
%\end{split}
\end{equation}
for $0 <   \frac{ \delta}{ n \hi i m \lo n  k \lo n}  \leq C \lo 2$.
% T22

For the term $T \lo {22}$ 
  note that
  \begin{align}\label{t22eq1}
 \| ( \Sigma \hi n  \lo {\sf N \hi i  \lo n \sf N \hi i  \lo n}+ \lambda \lo n D \hi * \lo {\sf N \hi i  \lo n}) \hi {-1} \frac{\mathbf{\tilde{H} \lo n  (\xi \hi {\sf N \hi i  \lo n})}}{\sqrt{n}} \| \lo {2} \leq \| ( \Sigma \hi n  \lo {\sf N \hi i  \lo n \sf N \hi i  \lo n}) \hi {-1} \frac{\mathbf{\tilde{H} \lo n  (\xi \hi {\sf N \hi i  \lo n})}}{\sqrt{n}} \| \lo {2}=\Lambda \lo {\mathrm{min}}( \Sigma \hi n  \lo {\sf N \hi i  \lo n \sf N \hi i  \lo n}) \hi {-1/2} ,
  \end{align}
  where we used   Lemma \ref{sampleboundedeigen} for the last inequality with $\delta=\frac{C \lo {\min}}{2}$.
Thus,  (on the event ${\cal N}$) the term $T \lo {22}$ can be bounded by
   \begin{align*}
 T \lo {22} \leq &  \Big (\frac{2} {C \lo {\min}} \Big ) \hi {3/2} \| \hat \Sigma \hi n  \lo {\sf N \hi i  \lo n \sf N \hi i  \lo n}-  \Sigma \hi n  \lo {\sf N \hi i  \lo n \sf N \hi i  \lo n} \| \lo {F}    \|\frac {w \hi i} {\sqrt{n}} \| \lo {F}.
    \end{align*}
{Recall the notation of $w \lo {uq}^i$ in \eqref{errorspline}
and the definition of the event $\Omega$ in Proposition
\ref{splineapproximation}.
Then, if the event $\Omega$ holds, we have,
\begin{equation}\label{bounderrorspine}
\begin{split}
 \|  \frac{1}{\sqrt{n}} w \hi i \| \lo {F} \hi 2& = \frac{1}{{n}} \tsum \lo {u=1} \hi {n} \tsum \lo {q=1} \hi {m \lo n} (w \hi i \lo {uq}) \hi 2  
 %\\ &
 = \frac{1}{{n}}  \tsum \lo {u=1} \hi {n} \tsum \lo {q=1} \hi {m \lo n} \Big ( \tsum \lo {j \in  \sf N \hi i}  \tsum \lo {r =1}  \hi {m \lo n} (f \hi {ij} \lo {qr} (\xi \hi {j} \lo {ur})-  \tilde f \hi {ij} \lo {nqr}  ( \xi \hi j \lo {ur})) \Big ) \hi 2  , \\
 & \leq \frac{n \hi i m \lo n}{{n}}  \tsum \lo {j \in  \sf N \hi i}   \tsum \lo {q=1} \hi {m \lo n}  \tsum \lo {r =1}  \hi {m \lo n} \tsum \lo {u=1} \hi {n}(f \hi {ij} \lo {qr} (\xi \hi {j} \lo {ur})-  \tilde f \hi {ij} \lo {nqr} ( \xi \hi j \lo {ur})) \hi 2  \\
 & 
 {\leq \frac{n \hi i m \lo n}{{n}}  \tsum \lo {j \in  \sf N \hi i}   \tsum \lo {q=1} \hi {m \lo n}  \tsum \lo {r =1}  \hi {m \lo n} 
 \max \lo {j \in  \sf N \hi i} \max \lo {1 \leq q, r \leq m \lo n}
 \| 
\mathbf {f \hi {ij} \lo {qr}}
- \mathbf {\tilde f \hi {ij} \lo {qr}} \| \lo {2}}
 ~\leq ~ c \lo 1  (n \hi i)  \hi 2m \lo n \hi {3} k \lo n \hi {-2d},
 \end{split}
  \end{equation}
 and by assumption \eqref{as1ridgeconsi1} it follows that on the 
 event $\Omega$
\begin{align}\label{m1}
 \|  \frac{1}{\sqrt{n}} w \hi i \| \lo {F}  \lesssim  \sqrt{\frac{2}{C \lo {\min}}} \lambda \lo n \tsum \lo {j \in  \sf N \hi i  \lo n } \| B \hi {*  i j} \lo {m \lo n k \lo n} \| \lo F.
\end{align}
 As a result, 
\begin{align} \nonumber 
 P \Big(T \lo {22} \geq \frac{\delta}{3}  \Big)
& \leq  P \Big( \Big(\frac{2}{C \lo {\min}} \Big) \hi {2} c \lo 1  \lambda \lo n \tsum \lo {j \in  \sf N \hi i  \lo n } \| B \hi {*  i j} \lo {m \lo n k \lo n} \| \lo F \| \hat \Sigma \hi n \lo {\sf N \hi i  \lo n \sf N \hi i  \lo n}-  \Sigma \hi n \lo {\sf N \hi i  \lo n \sf N \hi i  \lo n} \| \lo {F} \geq \frac{\delta }{3} \Big)+P(\Omega \hi {\complement}) \\
\label{t22final}
&\lesssim \exp \Big( - C \lo 3 \frac{n \hi {1- \alpha (2+3\beta)} (\lambda \lo n \tsum \lo {j \in  \sf N \hi i  \lo n } \| B \hi {*  i j} \lo {m \lo n k \lo n} \| \lo F )  \hi {-2} \delta \hi 2 }{(n \hi i) \hi 2 m \lo n \hi 2  k \lo n \hi {4}}+2 \log (n \hi i m \lo n k \lo n) \Big) \\
& \quad + \exp \Big(-C \lo 3 \frac{n}{n \hi i m \lo n \hi 2 k \lo n \hi{2d}}+\log(n \hi i m \lo n \hi {2}) \Big )
\nonumber 
 \end{align}
for $0 <  \delta (\lambda \lo n \tsum \lo {j \in  \sf N \hi i  \lo n } \| B \hi {*  i j} \lo n \| \lo F ) \inv  \leq C \lo 3$, where we have used 
Theorem \ref{theorem1} and   Proposition
\ref{splineapproximation}.  
 
% T23

We next derive an upper bound for the probability corresponding to the term $T \lo {23}$ in \eqref{hr2}   noting that (on the event ${\cal N}$)
{\begin{align}\label{det1}
T \lo {23} \leq \left(\frac{2}{C \lo {\min}} \right) \hi 2 \| \hat \Sigma \hi n \lo {\sf N \hi i \lo n \sf N \hi i \lo n}-  \Sigma \hi n \lo {\sf N \hi i \lo n \sf N \hi i \lo n} \| \lo {F}  \|  \frac{ \mathbf{\tilde{H} \lo n  (\xi \hi {\sf N \hi i \lo n})}}{n}  \epsilon \hi i  \| \lo {F}.
\end{align}}
 The $(i, j)$ element of the matrix $\frac{\mathbf{\tilde{H} \lo n  (\xi \hi {\sf N \hi i \lo n})}}{n}  \epsilon \hi i$ can be written as an i.i.d sum of the form $\frac{1}{n}\tsum \lo {u=1} \hi n \tilde{h} \lo {nk} ({\xi \hi {j} \lo {ur}) \epsilon \hi i \lo {uq}}$. 
Thus, by Assumption \ref{aserrors} it follows that 
{{\begin{align*}
P \Big( \Big| \frac{1}{n}\tsum \lo {u=1} \hi n \tilde{h} \lo {nk} ({\xi \hi {j} \lo {ur} ) \epsilon \hi i \lo {uq}} \Big| \geq \epsilon \Big) \leq 2 \exp (-C \lo 5 {n \epsilon \hi 2}),
\end{align*}}}
for any $\epsilon >0$.
Therefore, by applying the union bound over the $n \hi i m \lo n \hi 2 k \lo n$ gives
{{\begin{align}\label{t23erroreq1}
P \Big(   \|  \frac{\mathbf{\tilde{H} \lo n  (\xi \hi {\sf N \hi i \lo n})}}{n}  \epsilon \hi i  \| \lo {F} \geq \epsilon   \Big) \leq 2 \exp \Big(-C \lo 5 \frac{n  \epsilon \hi 2}{n \hi i m \lo n \hi 2 k \lo n}+ \log (n \hi i m \lo n \hi 2 k \lo n) \Big).
\end{align}}}
Using   this inequality with $\epsilon= C \lo {\min}  / {6} $, 
\eqref{t21final} and \eqref{det1} gives 
\begin{align}\label{t23final}
P \Big(T \lo {23} \geq \frac{\delta}{3} \Big) \lesssim \exp \Big( - C \lo 4 \frac{n \hi {1- \alpha (2+3\beta)} \delta \hi {2} }{(n \hi i) \hi 2 m \lo n \hi 2 k \lo n \hi 4}+2 \log (n \hi i m \lo n k \lo n) \Big)  +  \exp \Big(-C \lo 5 \frac{n   }{n \hi i m \lo n \hi 2 k \lo n}+ \log (n \hi i m \lo n \hi 2 k \lo n) \Big),
\end{align}
for $0 <  \frac{ \delta}{ n \hi i m \lo n  k \lo n}  \leq C \lo 4$.
Therefore from \eqref{t21final}, \eqref{t22final} and \eqref{t23final} it follows that
\begin{equation}\label{t2final}
%\begin{split}
P (T \lo 2 \geq {\delta})   \lesssim \exp \Big(- C \lo 2 \frac{n \hi {1- \alpha (2+3\beta)}   \delta \hi 2}{(n \hi i) \hi 2 m \lo n \hi 2 k \lo n \hi 4}+ 2 \log (n \hi i m \lo n k \lo n) \Big )  
    + 2 \exp \Big(-C \lo 3 \frac{n}{n \hi i m \lo n \hi 2 k \lo n \hi{2d}}+\log(n \hi i m \lo n \hi {2}) \Big),
%\end{split}
\end{equation}
for $0 <   \frac{ \delta}{ n \hi i m \lo n  k \lo n}  \leq C \lo 2$, {where the first term dominates the second one because
of Assumption \ref{assnu}.}

Finally, we  derive an upper bound for the probability involving $T \lo {3}$. Using representation \eqref{sigmazetai} we obtain
\begin{align} \label{t3eq1}
T \lo 3 \leq &  T_{31} + T_{32} + T_{33}~, 
\end{align}
where
\begin{align*} 
T_{31} = & \|  ( \Sigma \hi n \lo {\sf N \hi i \lo n \sf N \hi i \lo n} + \lambda \lo n D \hi * \lo {\sf N \hi i \lo n}) \hi {-1} \Sigma \hi n \lo {\sf N \hi i \lo n \sf N \hi i \lo n}B \hi {* \sf N \hi i \lo n} \lo n- B \hi {*\sf N \hi i \lo n } \lo n  \| \lo {F}, \\
T_{32}= & \| ( \Sigma \hi n \lo {\sf N \hi i \lo n \sf N \hi i \lo n}+ \lambda \lo n D \hi * \lo {\sf N \hi i \lo n }) \hi {-1} \frac{\mathbf{\tilde{H} \lo n  (\xi \hi {\sf N \hi i \lo n })}}{n}  w \hi i  \| \lo {F}, \\
T_{33}= & {\| ( \Sigma \hi n \lo {\sf N \hi i \lo n  \sf N \hi i \lo n }+ \lambda \lo n D \hi * \lo {\sf N \hi i \lo n }) \hi {-1} \frac{\mathbf{\tilde{H} \lo n  (\xi \hi {\sf N \hi i \lo n })}}{n}  \epsilon \hi i  \| \lo {F}.}
\end{align*}
For the first term on the right-hand side of the above inequality we have (on the event ${\cal N}$)
\begin{align*}
% \|  ( \Sigma \hi n  \lo {\sf N \hi i \lo n \sf N \hi i \lo n}+ \lambda \lo n D \hi * \lo {\sf N \hi i \lo n }) \hi {-1} \Sigma \hi n  \lo {\sf N \hi i \lo n \sf N \hi i \lo n}B \hi {* \sf N \hi i \lo n } \lo n- B \hi {*\sf N \hi i \lo n } \lo n   \| \lo {F}  &
 T_{31} =\lambda \lo n  \|  ( \Sigma \hi n  \lo {\sf N \hi i \lo n \sf N \hi i \lo n}+ \lambda \lo n D \hi * \lo {\sf N \hi i \lo n}) \hi {-1} D \hi * \lo {\sf N \hi i \lo n } B \hi {* \sf N \hi i \lo n} \lo n\| \lo {F} \leq
 %\\&= 
 \frac{2}{C \lo {\min}} \lambda \lo n (n \hi i) \hi {1/2} \sum \lo {j  \in  \sf N \hi i \lo n}\|B \hi {*ij}  \lo {m \lo n k \lo n} \| \lo {F} \leq \frac{\delta}{2(1+c \lo 1)},
\end{align*}
where we used condition  \eqref{as1ridgeconsi2} with $c_1 = \frac{1}{2c_2}-1$ and the fact that $\| \mathrm{diag} ( \frac{B \hi {*ij}  \lo {m \lo n k \lo n}}{\|B \hi {*ij}  \lo {m \lo n k \lo n} \| \lo {F}}: j \in {\sf N } \hi i \lo n  ) \| \lo F =(n \hi i) \hi {1/2}$.
Moreover, by applying the same arguments for deriving the bound of $T \lo {22}$ and by using \eqref{t22eq1}, conditions \eqref{as1ridgeconsi1} and  \eqref{as1ridgeconsi2} it follows that on  the event $\Omega$
\begin{align*}
T_{32}
 \leq c \lo 1 \left( \frac{2}{C \lo {\min}} \right)  \hi {1/2} \frac{n \hi i m \lo n \hi {3/2}}{k \lo n \hi d} 
 \leq {c \lo 1} \frac{2}{C \lo {\min}}  \lambda \lo n \sum \lo {j  \in  \sf N \hi i \lo n}\|B \hi {*ij}  \lo {m \lo n k \lo n} \| \lo {F} 
  \leq \frac{c \lo 1}{2 (1+ c \lo 1)}  \delta.
\end{align*}
 {Therefore, inequalities \eqref{t3eq1} and \eqref{t23erroreq1} imply that for all $\delta>0$
{\begin{align}\label{t3final}
P \left(  T \lo 3 \geq  {\delta} \right) &\leq P \Big(  \frac{2}{C \lo {\min}} \| \frac{ \mathbf{\tilde H \lo n  (\xi \hi {\sf N \hi i \lo n})} \epsilon \hi i}{n} \| \lo {F} \geq \frac{\delta}{2}      \Big)+ P(\Omega \hi {\complement}) \nonumber \\
&\leq 2 \exp \Big( - C \lo 3  \frac{n  \delta \hi 2 }{n \hi i m \lo n \hi 2 k \lo n}+ \log (n \hi i m \lo n \hi 2 k \lo n)  \Big)+2 \exp \Big( - C \lo 3  \frac{n   }{n \hi i m \lo n \hi 2 k \lo n \hi {2d}}+ \log (n \hi i m \lo n \hi 2)  \Big)
\end{align}}
Thus, by \eqref{t1final}, \eqref{t2final} and \eqref{t3final},  we have shown, for any $\delta>$ such that $0 <  \delta \leq C \lo 1$ and  $0 <  \delta    \leq C \lo 2$
\begin{align*}
P (\| \tilde B \hi {\sf N \hi i \lo n} \lo n- B \hi {*\sf N \hi i \lo n} \lo n \| \lo {F} \geq \delta) \lesssim &   \exp \Big ( -C \lo 1 \frac{n \hi {1- \alpha (2+3 \beta)}  \delta \hi 2}{n \hi i    m \lo n \hi 2 k \lo n \hi {3}} + \log (n \hi i m \lo n \hi 2 k \lo n)\Big) \\
&+ \exp \Big( - C \lo 2 \frac{n \hi {1- \alpha (2+3\beta)}   \delta \hi 2}{(n \hi i) \hi 2 m \lo n  \hi 2k \lo n \hi 4}+2 \log (n \hi i m \lo n k \lo n) \Big) \\
 &+  \exp \Big( - C \lo 3 \frac{n  \delta \hi 2  }{n \hi i m \lo n \hi 2 k \lo n}+ \log (n \hi i m \lo n \hi 2 k \lo n) \Big)  .
\end{align*}
Since the second term dominates the first and the third, the assertion in Proposition \ref{ridgeconsi} follows.  \eop

%%%%%%%%%proposition 2%%%%%%%%%%%%%%%%%%
The next proposition brings  $\hat B \hi {\sf N \hi i \lo n} \lo n=( \hat B \hi {ij} \lo n, j \in \sf N \hi i \lo n)$ close to $\tilde B \hi {\sf N \hi i \lo n} \lo n=(\tilde B \hi {ij} \lo n,  j \in \sf N \hi i \lo n)$, from which we can establish the concentration inequality for $\hat B \hi {\sf N \hi i \lo n} \lo n$.
%%%%%%%%%%%%%%%proposition 3%%%%%%%%%%%%%%%%%%%%%%%%%%%%%
\begin{proposition}\label{bhatconsi} Let $\hat B \hi {\sf N \hi i \lo n} \lo n$ be the minimiser of  \eqref{restricted} and $\tilde B \lo n \hi {\sf N \hi i \lo n}$ be the minimiser of\eqref{restrictedridge}.  If $\Lambda \lo {\min} (\hat \Sigma \hi n  \lo {\sf N \hi i \lo n \sf N \hi i \lo n}) \geq \frac{C \lo {\min}}{4}$ then,
\begin{align*}
\| \hat B \hi {\sf N \hi i \lo n} \lo n-  \tilde B \hi {\sf N \hi i \lo n} \lo n \|  \lo {F}  \leq \frac{10}{C \lo {\min}} \lambda \lo n n \hi i  (b \hi {*i} \lo n) \hi {-1} \| \tilde B \hi {\sf N \hi i \lo n} \lo n -  B \hi {*\sf N \hi i \lo n} \lo n \|  \lo {F} \tsum \lo {j \in  \sf N \hi i \lo n } \|  B \hi {*ij} \lo {m \lo n k \lo n} \| \lo F, 
\end{align*}
where $b \hi {*i} \lo n= \min \lo {j \in \sf N \hi i \lo n} \| B \hi {* ij } \lo {m \lo n k \lo n} \| \lo F$.
\end{proposition}

\proof  The idea of the proof is similar as in the proof of Proposition 3 in \cite{leefaro}. Consider the sphere $\mathcal{S} \lo n(\delta \lo n)=\{B \in \R \hi {n \hi i k \lo n m \lo n \times m \lo n} :  \| B- \tilde B \lo n \hi {\sf N \hi i \lo n} \| \lo {F}=\delta \lo n \}$, where $(\delta \lo n) \lo {n \in \mathbb{N}} $ is a positive sequence of real numbers.  For $\epsilon \in [0,1]$ let
$$ 
f (\epsilon)=\widehat {PL}  \lo {\sf N \hi i \lo n}  ( \tilde B \hi {\sf N \hi i \lo n}  \lo n+ \epsilon A  , \hat \xi \hi i  ), 
$$
where the function $\widehat {PL}  \lo {\sf N \hi i \lo n}$ is defined in \eqref{restricted2} and $A=B - \tilde B \lo n \hi {\sf N \hi i \lo n}$. 
 A straightforward calculation gives for the first and the second derivatives of the function $f (\epsilon)$
\begin{align*} %\label{hol5}
 \dot{f} (\epsilon)&=- \langle \hat \Sigma \hi n  \lo {\sf N \hi i \lo n \xi \hi i}, A \rangle \lo F +  \langle A, \hat \Sigma \hi n  \lo {\sf N \hi i \lo n \sf N \hi i \lo n} (  \tilde B \hi {\sf N \hi i \lo n}  \lo { n}+\epsilon A) \rangle \lo F\\ \nonumber 
 &  +  \lambda \lo n  \sum \lo {j \in \sf N \hi i \lo n } \| \tilde B \hi {ij}  \lo { n} + \epsilon A \hi j \| \lo F \sum \lo {k \in \sf N \hi i \lo n } \left( \langle \tilde B \hi {ik}  \lo {n} + \epsilon A \hi k,  A \hi k \rangle \lo F \| \tilde B \hi {ik}  \lo { n} + \epsilon A \hi k \| \lo F \hi {-1}\right). \\ 
 %\label{hol6}
 {\ddot{f}(\epsilon)}&=  \langle A, \hat \Sigma \hi n  \lo {\sf N \hi i \lo n \sf N \hi i \lo n} A  \rangle \lo F +  \lambda \lo n  \Big ( \sum \lo {j \in \sf N \hi i \lo n } \| \tilde B \hi {ij}  \lo { n} + \epsilon A \hi j \| \lo F \hi {-1} \langle \tilde B \hi {ij}  \lo {n} + \epsilon A \hi j,  A \hi j \rangle \lo F \Big ) \hi 2 \\ \nonumber
 &+  \lambda \lo n  \sum \lo {j \in \sf N \hi i \lo n} \| \tilde B \hi {ij}  \lo { n} + \epsilon A \hi j \| \lo F
 \sum \lo {k \in \sf N \hi i \lo n }  \| A \hi k \| \hi 2 \lo F \| \tilde B \hi {ik}  \lo { n} + \epsilon A \hi k \| \hi {-1} \lo F \\ \nonumber
 &- \lambda \lo n  \sum \lo {j \in \sf N \hi i \lo n } \| \tilde B \hi {ij}  \lo { n} + \epsilon A \hi j \| \lo F \sum \lo {k \in \sf N \hi i \lo n } \left( \langle \tilde B \hi {ik}  \lo { n} + \epsilon A \hi k,  A \hi k \rangle \lo F \hi 2 \| \tilde B \hi {ik}  \lo { n} + \epsilon A \hi k \| \lo F \hi {-3}\right).
\end{align*}
where $\langle \cdot , \cdot \rangle \lo F$ denotes the Frobenious inner product and $ \hat \Sigma \hi n  \lo {\sf N \hi i \lo n \xi \hi i}$, $\hat \Sigma \hi n  \lo {\sf N \hi i \lo n \sf N \hi i \lo n}$ are defined in \eqref{hd10} and \eqref{estimatedmat3} respectively.
By construction, $f(0)=\widehat {PL} \lo {\sf N \hi i \lo n} ( \tilde B \hi {\sf N \hi i \lo n}  \lo n, \hat \xi \hi i  )$,   $f(1)=\widehat {PL}  \lo {\sf N \hi i \lo n}  (  B,  \hat \xi \hi i   )$ and by Taylor's theorem, we have for some $\epsilon \in (0, 1)$
\begin{align}\label{tayolorsexp}
\widehat {PL}  \lo {\sf N \hi i \lo n}  (  B, \hat \xi \hi i )-\widehat {PL}  \lo {\sf N \hi i \lo n}  ( \tilde B \hi {\sf N \hi i \lo n} \lo n,  \hat \xi \hi i  )=f(1)-f(0)= \dot{f} (0) +\frac{\ddot{f}(\epsilon)}{2}.
\end{align}
The Cauchy-Schwarz inequality yields for any $B \in \mathcal{S} \lo n(\delta \lo n)$ and   $\epsilon \in [0,1]$
\begin{align}\label{boundsecondderiv}
{\ddot{f}(\epsilon)} & \geq \langle A, \hat \Sigma \hi n  \lo {\sf N \hi i \lo n \sf N \hi i \lo n} A  \rangle \lo F +   \lambda \lo n  \Big ( \sum \limits{\lo {j \in \sf N \hi i \lo n}} \| \tilde B \hi {ij}  \lo { n} + \epsilon A \hi j \| \lo F \hi {-1} \langle \tilde B \hi {ij}  \lo { n} + \epsilon A \hi j,  A \hi j \rangle \lo F \Big ) \hi 2  \geq \frac{C \lo {\min}}{2} \delta \lo n \hi 2.
\end{align}
On the other hand, by Lemma A7 in \cite{lee2016variable} it follows that
 \begin{align*}
 | \dot{f} (0) | &\leq  \lambda \lo n \tsum \lo {j \in \sf N \hi i \lo n} \tsum \lo {k \in \sf N \hi i \lo n} \Big [ \|  \tilde B \hi { ij} \lo { n} - B \hi {* ij}  \lo {m \lo n k \lo n}   \| \lo {F} \|  \tilde B \hi {ik} \lo { n} - B \hi { k}  \| \lo {F}  \\
&  +  \|  B \hi {* ij} \lo {m \lo n k \lo n} \| \lo {F}  \|  B \hi {*ik} \lo {m \lo n k \lo n}  \| \lo {F} \hi {-1} \|  \tilde B \hi { ij} \lo { n} - B \hi {* ij}  \lo {m \lo n k \lo n}   \| \lo {F} \|  \tilde B \hi {ik} \lo { n} - B \hi { k}  \| \lo {F} \Big].
\end{align*}
A further application of the Cauchy-Schwarz  inequality gives
 \begin{align*} 
| \dot{f} (0) |  &\leq  \lambda \lo n \sqrt{n \hi i \tsum \lo {j \in \sf N \hi i \lo n} \|  \tilde B \hi { ij}  \lo n - B \hi {* ij}  \lo {m \lo n k \lo n}   \| \hi 2 \lo {F} n \hi i \tsum \lo {k \in \sf N \hi i \lo n} \|  \tilde B \hi {ik} \lo {n} - B \hi { k}  \| \hi 2 \lo {F} } \nonumber \\
&  +  \lambda \lo n   \tsum \lo {j \in \sf N \hi i \lo n}   \|  B \hi {* ij} \lo {m \lo n k \lo n} \| \lo {F}  \|  \tilde B \hi { ij} \lo { n} - B \hi {* ij}  \lo {m \lo n k \lo n}   \| \lo {F} (b \hi {*i} \lo n) \hi {-1} \tsum \lo {k \in \sf N \hi i \lo n} \|  \tilde B \hi {ik} \lo { \lo n} - B \hi { k}  \| \lo {F} \nonumber \\
&=  \lambda \lo n  n \hi i  \|  \tilde B \hi {\sf N \hi i \lo n} \lo n - B \hi {* \sf N \hi i \lo n   }  \lo n   \|  \lo {F}  \|  \tilde B \hi {\sf N \hi i \lo n} \lo n - B       \|  \lo {F} \nonumber \\
&+  \lambda \lo n  (b \hi {*i} \lo n) \hi {-1}  \tsum \lo {j \in \sf N \hi i \lo n}   \|  B \hi {* ij} \lo {m \lo n k \lo n} \| \lo {F} \|  \tilde B \hi { ij} \lo {\lo n} - B \hi {* ij}  \lo {m \lo n k \lo n}   \| \lo {F}  \tsum \lo {k \in \sf N \hi i \lo n} \|  \tilde B \hi {ik} \lo {\lo n} - B \hi { k}  \| \lo {F}  \nonumber \\
& \leq  \lambda \lo n  n \hi i  \|  \tilde B \hi {\sf N \hi i \lo n } \lo n - B \hi {* \sf N \hi i  \lo n}  \lo n   \|  \lo {F}  \|  \tilde B \hi {\sf N \hi i \lo n} \lo n - B      \|  \lo {F} \nonumber \\
&+  \lambda \lo n  \sqrt{n \hi i}  (b \hi {*i} \lo n) \hi {-1}  \|  \tilde B \hi {\sf N \hi i \lo n} \lo n - B \hi {* \sf N \hi i \lo n  }  \lo n   \|  \lo {F}  \|  \tilde B \hi {\sf N \hi i \lo n} \lo n - B     \|  \lo {F} \| B \hi {* \sf N \hi i \lo n }  \lo n \| \lo {F}.
\end{align*}
Using the fact $b \hi {*i} \lo n \leq   \| B \hi {* \sf N \hi i \lo n }  \lo n \| \lo {F} \leq \tsum \lo {j \in  \sf N \hi i \lo n} \|  B \hi {*ij} \lo {m \lo n k \lo n} \| \lo F$ we obtain
%(observing that $B \in {\cal S}_n(\delta_n)$) 
 \begin{align}\label{boundfirstderiv}
 | \dot{f} (0) |  \leq  2  \lambda \lo n  n \hi i  (b \hi {*i} \lo n)  \hi {-1} \tsum \lo {j \in  \sf N \hi i \lo n } \|  B \hi {*ij} \lo {m \lo n k \lo n} \| \lo F \|  \tilde B \hi {\sf N \hi i \lo n } \lo n - B \hi {* \sf N \hi i  \lo n}  \lo n   \|  \lo {F}  \|  \tilde B \hi {\sf N \hi i \lo n} \lo n - B      \|  \lo {F}.
 \end{align}
 Hence, combining \eqref{tayolorsexp}, \eqref{boundsecondderiv} and \eqref{boundfirstderiv}, we obtain
 \begin{align*}
\widehat {PL}  \lo {\sf N \hi i \lo n}  (  B , \hat \xi \hi i )-\widehat {PL}  \lo {\sf N \hi i \lo n}  ( \tilde B \hi {\sf N \hi i \lo n} \lo n, \hat \xi \hi i ) &\geq -   2 \lambda \lo n  n \hi i  (b \hi {*i} \lo n) \hi {-1} \tsum \lo {j \in  \sf N \hi i \lo n } \|  B \hi {*ij} \lo {m \lo n k \lo n} \| \lo F \|  \tilde B \hi {\sf N \hi i \lo n} \lo n - B \hi {* \sf N \hi i \lo n }  \lo n   \|  \lo {F}   +\frac{C \lo {\min}}{4} \delta \lo n \hi 2
\end{align*}
If we choose $ \delta \lo n \hi 2 =  \frac{10}{C \lo {\min}} \lambda \lo n n \hi i  (b \hi {*i} \lo n) \hi {-1} \tsum \lo {j \in  \sf N \hi i \lo n} \|  B \hi {*ij} \lo {m \lo n k \lo n} \| \lo F \| \tilde B \hi {\sf N \hi i \lo n} \lo n -  B \hi {*\sf N \hi i \lo n} \lo n \|  \lo {F}$,  it follows that
  \begin{align*}
  \widehat {PL}  \lo {\sf N \hi i \lo n}  (  B , \hat \xi \hi i )-\widehat {PL}  \lo {\sf N \hi i \lo n}  ( \tilde B \hi {\sf N \hi i \lo n} \lo n,  \hat \xi \hi i  ) > 0.
  \end{align*}
 Since the function $\widehat {PL}   \lo {\sf N \hi i \lo n}  (  B ,  \hat \xi \hi i   ) $ is convex, the minimizer $\hat B \hi {\sf N \hi i \lo n} \lo n $ of $\widehat {PL}   \lo {\sf N \hi i \lo n}  (  B ,  \hat \xi \hi i   )$ is going be inside the sphere defined by $\mathcal{S} \lo n(\delta \lo n)$,  that is,
 $$
 \| \hat B \hi {\sf N \hi i \lo n} \lo n-  \tilde B \hi {\sf N \hi i \lo n} \lo n \|  \lo {F}  \leq \frac{10}{C \lo {\min}} \lambda \lo n n \hi i  (b \hi {*i} \lo n) \hi {-1} \tsum \lo {j \in  \sf N \hi i \lo n} \|  B \hi {*ij} \lo {m \lo n k \lo n} \| \lo F\| \tilde B \hi {\sf N \hi i \lo n} \lo n -  B \hi {*\sf N \hi i \lo n} \lo n \|  \lo {F}.
 $$ ~ 
\eop

\medskip
%%%%concentration inequality %%%%%%%%%%%%%%%%%%%%%%%%%%%%%%%%%%%%%%%
\noindent
Using the Propositions \ref{ridgeconsi} and \ref{bhatconsi}, we now can establish the concentration bounds for $\| \hat B \hi {\sf N \hi i \lo n} \lo n-  B \hi {*\sf N \hi i \lo n}  \lo n \| \lo {F}$. 
{\begin{proposition}\label{propconsistency}
Suppose Assumptions of Proposition \ref{ridgeconsi} are satisfied      and that $\delta$ satisfies
\begin{align}\label{cond:propconsistency}
\frac{2}{C \lo {\min}} \lambda \lo n (n \hi i) \hi {3/2} (\tsum \lo {j \in  \sf N \hi i \lo n} \|  B \hi {*ij} \lo {m \lo n k \lo n} \| \lo F) \hi 2 \leq c \lo 2 b \hi {*i} \lo n \delta 
\end{align}
for some constant $c_2 >0$. Then,
\begin{align*}
P \Big( \| \hat B \hi {\sf N \hi i \lo n} \lo n-  B \hi {*\sf N \hi i \lo n}  \lo n \| \lo {F} \geq \delta \Big) \lesssim & \exp \Big( - C \lo 1 \frac{n \hi {1-\alpha (2+3\beta)} (b \hi {*i} \lo n) \hi 2 \delta \hi 2}{(n \hi i) \hi 4 m \lo n \hi 2 k \lo n \hi {4} (\tsum \lo {j \in  \sf N \hi i \lo n} \|  B \hi {*ij} \lo {m \lo n k \lo n} \| \lo F) \hi 2}+2 \log (n \hi i m \lo n k \lo n) \Big),
\end{align*}
 where $C \lo 1>0$   such that $0 < \delta \leq C \lo 1$.
\end{proposition}} 
\proof By  Proposition \ref{bhatconsi} and the triangle inequality,
\begin{align*}
\| \hat B \hi {\sf N \hi i \lo n} \lo n -   B \hi {*\sf N \hi i \lo n} \lo n \|  \lo {F} & \leq \| \hat B \hi {\sf N \hi i \lo n} \lo n-  \tilde B \hi {\sf N \hi i \lo n} \lo n \|  \lo {F} +  \|  B \hi {*\sf N \hi i \lo n}  \lo n-  \tilde B \hi {\sf N \hi i \lo n} \lo n \|  \lo {F} \\
&\leq (b \hi {*i} \lo n) \hi {-1}  \|  B \hi {*\sf N \hi i \lo n} \lo n-  \tilde B \hi {\sf N \hi i \lo n} \lo n \|  \lo {F}  \Big (   \frac{10}{C \lo {\min}} \lambda \lo n n \hi i \tsum \lo {j \in  \sf N \hi i \lo n} \|  B \hi {*ij} \lo {m \lo n k \lo n} \| \lo F+  b \hi {*i} \lo n  \Big )\\
&  \lesssim n \hi i (b \hi {*i} \lo n) \hi {-1}  \|  B \hi {*\sf N \hi i \lo n} \lo n-  \tilde B \hi {\sf N \hi i \lo n} \lo n \|  \lo {F}     \tsum \lo {j \in  \sf N \hi i \lo n} \|  B \hi {*ij} \lo {m \lo n k \lo n} \| \lo F,
\end{align*}

where  we have used the fact that $b \hi {*i} \lo n \leq   \tsum \lo {j \in  \sf N \hi i \lo n} \|  B \hi {*ij} \lo {m \lo n k \lo n} \| \lo F$ and   $\lambda \lo n \lesssim 1$. 
The assertion now follows from Proposition \ref{ridgeconsi} with $\delta$ replaced by ${b \hi {*i} \lo n \delta}\big ( n \hi i  \tsum \lo {j \in  \sf N \hi i \lo n} \|  B \hi {*ij} \lo {m \lo n k \lo n} \| \lo F \big)^{-1}$.
 \eop

  %%%%%%%%%%%%%%%%%%%%%%%%%%%%%%%%%%%%%%%%%%%%%%%%
 %%%%%%%%%%%%%%%%last proposition -second KKT%%%%%%%%%%%%%%%%%%%%%%%%%%%%%
%%%%%%%%%%%%%%%%%%%%%%%%%%%%%%%%%%%%%%%%%%%%%%%%

 \noindent
 Let $\hat B \hi {\sf N \hi i \lo n} \lo n$ be the minimizer of the restricted problem \eqref{restricted}.   By construction, the estimator $(\hat B \hi {\sf N \hi i \lo n} \lo n, \mathbf{0})$ obtained from $\hat B \hi {\sf N \hi i \lo n} \lo n$ by adding blocks with 0 elements whenever $j \notin \sf N \hi i \lo n$, satisfies the first KKT-condition \eqref{kktfirst}.  To prove that $(\hat B \hi {\sf N \hi i \lo n} \lo n, \mathbf{0})$ is, with high probability, optimal for problem \eqref{restricted}, it is therefore sufficient to show that the second KKT-condition \eqref{kktsecond} is satisfied. This is the statement of the following proposition.

 \begin{proposition}\label{bhatoptim} The matrix  $(\hat B \hi {\sf N \hi i \lo n } \lo n, \mathbf{0})$ satisfies  \eqref{kktsecond} with high probability, in the sense that  
 \begin{eqnarray*}
&& P ( \max \lo {j \notin \sf N \hi i \lo n  } \|\hat \Sigma \hi n \lo { \xi \hi j  {\sf N \hi i \lo n }} \hat B \hi {{\sf N \hi i \lo n }} \lo n - \hat {\Sigma}  \hi n \lo { \xi \hi j {\xi \hi i}} \| \lo {F} \geq \lambda \lo n \tsum \lo {j \ne i} \hi {p} \| \hat B \hi {ij} \lo n \| \lo F ) 
\\
&&\quad \quad \quad \quad \quad \quad 
 \lesssim  \exp \Big( - C \lo 1 \frac{n \hi {1- \alpha (2+3\beta)} (\lambda \lo n \tsum \lo {j \in  \sf N \hi i \lo n } \|  B \hi {*ij} \lo {m \lo n k \lo n} \| \lo F) \hi 2  }{n \hi i m \lo n \hi 2  k \lo n \hi {4} }+2 \log (n \hi i m \lo n k \lo n) \Big)   ,
  \end{eqnarray*}
   where $C \lo 1$ is a positive constant.
 \end{proposition}

  %%proof
\proof The idea of the proof is similar as in the proof of Proposition 4 in \cite{leefaro}.
 By the first optimality condition \eqref{kktfirst}, we have for all $j \in \sf N \hi i \lo n$,
  \begin{align}\label{firstoptim}
\hat B \hi {\sf N \hi i \lo n} \lo n=(\hat \Sigma  \hi n \lo { \sf N \hi i \lo n \sf N \hi i \lo n} + \lambda \lo n \hat D \lo {\sf N \hi i \lo n}) \hi {-1}  \hat {\Sigma} \hi n  \lo { \sf N \hi i {\xi \hi i}},
 \end{align}
where $\hat D \lo {\sf N \hi i \lo n}$ is defined in Lemma \ref{lemmakkt}.  Using  \eqref{firstoptim} in the expression  at the left-hand side of condition \eqref{kktsecond}  gives
  \begin{align*}
   \hat \Sigma \hi n  \lo { \xi \hi j \sf N \hi i \lo n} \hat B \hi {\sf N \hi i \lo n} \lo n - \hat {\Sigma} \hi n  \lo { \xi \hi j {\xi \hi i}} &=  \hat \Sigma \hi n  \lo { \xi \hi j \sf N \hi i \lo n}(\hat \Sigma \hi n  \lo { \sf N \hi i \lo n \sf N \hi i \lo n} + \lambda \lo n \hat D \lo {\sf N \hi i \lo n}) \hi {-1}  \hat {\Sigma} \hi n  \lo { \sf N \hi i \lo n {\xi \hi i}}- \hat {\Sigma} \hi n  \lo { \xi \hi j {\xi \hi i}}= R \hi j \lo 1+ \ldots +R \hi j \lo 7,
  \end{align*}
where
    \begin{align*}
   R \lo 1 \hi j&=(  \hat \Sigma \hi n  \lo { \xi \hi j \sf N \hi i \lo n} -   \Sigma \hi n  \lo { \xi \hi j \sf N \hi i \lo n} )(\hat \Sigma \hi n  \lo { \sf N \hi i \lo n \sf N \hi i \lo n} + \lambda \lo n \hat D \lo {\sf N \hi i \lo n}) \hi {-1} {\Sigma} \hi n  \lo { \sf N \hi i \lo n {\xi \hi i}}\\
   R \lo 2 \hi j&=\Sigma \hi n  \lo { \xi \hi j \sf N \hi i \lo n} (\hat \Sigma \hi n  \lo { \sf N \hi i \lo n \sf N \hi i \lo n} + \lambda \lo n \hat D \lo {\sf N \hi i \lo n}) \hi {-1}(\hat {\Sigma} \hi n  \lo { \sf N \hi i \lo n {\xi \hi i}}- {\Sigma} \hi n  \lo { \sf N \hi i \lo n {\xi \hi i}}) \\
    R \lo 3 \hi j&= \Sigma \hi n  \lo { \xi \hi j \sf N \hi i \lo n} \{(\hat \Sigma \hi n  \lo { \sf N \hi i \lo n \sf N \hi i \lo n} + \lambda \lo n \hat D \lo {\sf N \hi i \lo n}) \hi {-1}-( \Sigma \hi n  \lo { \sf N \hi i \lo n \sf N \hi i \lo n} + \lambda \lo n \hat D \lo {\sf N \hi i \lo n}) \hi {-1} \} {\Sigma} \hi n  \lo { \sf N \hi i \lo n {\xi \hi i}}\\
    R \lo 4 \hi j&= (  \hat \Sigma \hi n  \lo { \xi \hi j \sf N \hi i \lo n} -   \Sigma \hi n  \lo { \xi \hi j \sf N \hi i \lo n} )(\hat \Sigma \hi n  \lo { \sf N \hi i \lo n \sf N \hi i \lo n} + \lambda \lo n \hat D \lo {\sf N \hi i \lo n}) \hi {-1} (\hat {\Sigma} \hi n  \lo { \sf N \hi i \lo n {\xi \hi i}}- {\Sigma} \hi n  \lo { \sf N \hi i \lo n {\xi \hi i}}) \\
    R \lo 5 \hi j&= {\Sigma} \hi n  \lo { \xi \hi j {\xi \hi i}} -\hat {\Sigma} \hi n  \lo { \xi \hi j {\xi \hi i}} \\
    R \lo 6 \hi j&=  \Sigma \hi n  \lo { \xi \hi j \sf N \hi i \lo n} \{ ( \Sigma \lo { \sf N \hi i \lo n \sf N \hi i \lo n} + \lambda \lo n \hat D \lo {\sf N \hi i \lo n}) \hi {-1}-( \Sigma \hi n  \lo { \sf N \hi i \lo n \sf N \hi i \lo n} + \lambda \lo n D \hi * \lo {\sf N \hi i \lo n}) \hi {-1} \}  {\Sigma} \hi n  \lo { \sf N \hi i \lo n {\xi \hi i}}\\
    R \lo 7 \hi j&=  \Sigma \hi n  \lo { \xi \hi j \sf N \hi i \lo n} ( \Sigma \hi n  \lo { \sf N \hi i \lo n \sf N \hi i \lo n} + \lambda \lo n D \hi * \lo {\sf N \hi i \lo n}) \hi {-1}   {\Sigma} \hi n  \lo { \sf N \hi i \lo n {\xi \hi i}}- {\Sigma} \hi n  \lo { \xi \hi j {\xi \hi i}}.
     \end{align*}
     {
In the following we derive bounds for the probabilities 
   \begin{align}
    P \Big( \max \lo {j \notin  {\sf N \hi i \lo n}}   {\| R \lo r \hi j \| \lo {F}} \geq \frac{ \lambda \lo n}{7} \tsum \lo {j \in  \sf N \hi i \lo n} \| \hat B \hi {ij} \lo {   n} \| \lo F   \Big) 
    \label{deta}
   \end{align}  
   ($r=1,\ldots, 7$). 
  For this purpose we proceed in two steps.   }  
\\
   %%% step 1
{\bf Step 1}:
First, we  define the event, that there exists a constant $0<c \lo 0 <1$, such that
\begin{align*}
\mathcal{A} \lo 0=\Big \{ \| \hat B \hi {\sf N \hi i \lo n} \lo n -  B \hi {*\sf N \hi i \lo n} \lo n \| \lo {F} \leq c \lo 0  (n \hi i) \hi {-1/2}  \tsum \lo {j \in  \sf N \hi i \lo n} \|  B \hi {*ij} \lo {m \lo n k \lo n} \| \lo F   \Big \}.
\end{align*}
  Then, by Proposition \ref{propconsistency} with  
  $\delta=c_0(n^i)^{-1/2} \sum_{j \in \sf N^i}\|B^*_{m_n k_n}\|_F$ we have
\begin{equation}\label{step1}
\begin{split}
P ( \mathcal{A} \lo 0 \hi {\complement}) &\lesssim \exp \Big( - C \lo 1 \frac{n \hi {1- \alpha (2+3\beta)} (b \hi {*i} \lo n) \hi 2  }{(n \hi i) \hi {5} m \lo n \hi 2  k \lo n \hi {4} }+2 \log (n \hi i m \lo n k \lo n) \Big) \\
&\lesssim \exp \Big( - C \lo 1 \frac{n \hi {1- \alpha (2+3\beta)} (\lambda \lo n \tsum \lo {j \in  \sf N \hi i \lo n } \|  B \hi {*ij} \lo {m \lo n k \lo n} \| \lo F) \hi 2  }{n \hi i m \lo n \hi 2  k \lo n \hi {4} }+2 \log (n \hi i m \lo n k \lo n) \Big),
\end{split}
\end{equation}
for some constant $C \lo 1>0$.

Now,  on the event $\mathcal{A} \lo 0$, it follows by  the Cauchy-Schwarz inequality that
\begin{align*}
 \Big | \tsum \lo {j \in  \sf N \hi i \lo n} \|  \hat B \hi {ij} \lo n \| \lo F-\tsum \lo {j \in  \sf N \hi i \lo n} \|  B \hi {*ij} \lo {m \lo n k \lo n} \| \lo F  \Big | \leq c \lo 0  \tsum \lo {j \in  \sf N \hi i \lo n} \| B \hi {*ij} \lo {m \lo n k \lo n}\| \lo F .
 \end{align*}
 Therefore, we have for $r=1,\ldots,7$
\begin{align*}
P \Big(\max \lo {j \notin  {\sf N \hi i \lo n}}   {\| R \lo r \hi j \| \lo {F}} \geq \frac{ \lambda \lo n}{7} \tsum \lo {j \in  \sf N \hi i \lo n} \| \hat B \hi {ij} \lo n \| \lo F      \Big)   \leq P   \Big(\max \lo {j \notin  {\sf N \hi i \lo n}}   {\| R \lo r \hi j \| \lo {F}}    \geq \frac{ \lambda \lo n (1-c \lo 0)}{7} \tsum \lo {j \in  \sf N \hi i \lo n} \|  B \hi {*ij} \lo {m \lo n k \lo n} \| \lo F \Big)  
  + P(\mathcal{A}^c_0),
\end{align*}
where the second probability can be bounded by \eqref{step1}.
Thus, in order to control the probabilities of the terms $ \max \lo {j \notin  {\sf N \hi i \lo n}}   {\| R \lo r \hi j \| \lo {F}}    \geq \frac{ \lambda \lo n}{7} \tsum \lo {j \in  \sf N \hi i \lo n} \| \hat B \hi {ij} \lo n \| \lo F $ it  suffices to derive the probabilities $P (\max \lo {j \notin {\sf N \hi i \lo n}}   {\| R \lo r \hi j \| \lo {F}}    \geq  c \lo 1 \lambda \lo n \tsum \lo {j \in  \sf N \hi i \lo n} \|  B \hi {*ij} \lo {m \lo n k \lo n} \| \lo F   )$ for all $r=1, \ldots, 7$, where $c \lo 1=\frac{1-c \lo 0}{7}$. 

\smallskip
   %%% step 2
\noindent
{\bf Step 2}:  {\it {Term $R \lo 1 \hi j$}}.  Substituting the representation $\eqref{sigmazetai}$ used in the proof of Proposition \ref{ridgeconsi} we can rewrite $R \lo 1 \hi j$ as
 \begin{align*}
      R \lo 1 \hi j = \leq   ( \hat \Sigma \hi n  \lo { \xi \hi j \sf N \hi i \lo n } -   \Sigma \hi n  \lo { \xi \hi j \sf N \hi i \lo n })
 (\hat \Sigma \hi n  \lo { \sf N \hi i \lo n  \sf N \hi i \lo n } + \lambda \lo n \hat D \lo {\sf N \hi i}) \hi {-1} \Big( \Sigma \hi n  \lo {\sf N \hi i \lo n  \sf N \hi i \lo n }B \hi {* \sf N \hi i \lo n } \lo n+ \frac{\mathbf{\tilde{H} \lo n  (\xi \hi {\sf N \hi i \lo n })}}{n}  w \hi i +{\frac{\mathbf{\tilde{H} \lo n  (\xi \hi {\sf N \hi i \lo n })}}{n}  \epsilon \hi i \Big)},
    \end{align*}
    and by the union bound it follows that
    \begin{align} \label{erg1}
    P \Big(  \max \lo {j \notin  {\sf N \hi i \lo n}}   {\| R \lo 1 \hi j\| \lo {F}}  \geq  c \lo 1 \lambda \lo n \tsum \lo {j \in  \sf N \hi i \lo n } \|  B \hi {*ij} \lo {m \lo n k \lo n} \| \lo F  \Big) \leq \tsum \lo {r=1} \hi {3} P \Big  ( \max \lo {j \notin  {\sf N \hi i \lo n}} {\| R \lo {1r} \hi j \| \lo {F}} \geq  \frac{c \lo 1 \lambda \lo n}{3} \tsum \lo {j \in  \sf N \hi i \lo n } \|  B \hi {*ij} \lo {m \lo n k \lo n} \| \lo F \Big ),
    \end{align}
     where $R \lo {11} \hi j$, $R \lo {12} \hi j$, $R \lo {13} \hi j$ are defined in an obvious manner satisfying
 \begin{align*}
  \|   R \lo {11} \hi j\| \lo F \leq &  \|   \hat \Sigma \hi n  \lo { \xi \hi j \sf N \hi i \lo n } -   \Sigma \hi n  \lo { \xi \hi j \sf N \hi i \lo n } \| \lo F  \| (\hat \Sigma \hi n  \lo { \sf N \hi i \lo n  \sf N \hi i \lo n } + \lambda \lo n \hat D \lo {\sf N \hi i \lo n }) \hi {-1} \Sigma \hi n  \lo {\sf N \hi i \lo n  \sf N \hi i \lo n }B \hi {* \sf N \hi i \lo n } \lo n \| \lo 2.  \\
    \|   R \lo {12} \hi j \| \lo F \leq &   \|   \hat \Sigma \hi n  \lo { \xi \hi j \sf N \hi i \lo n } -   \Sigma \hi n  \lo { \xi \hi j \sf N \hi i \lo n } \| \lo F \|(\hat \Sigma \hi n  \lo { \sf N \hi i \lo n  \sf N \hi i \lo n } + \lambda \lo n \hat D \lo {\sf N \hi i \lo n }) \hi {-1} \frac{\mathbf{\tilde{H} \lo n  (\xi \hi {\sf N \hi i \lo n })}}{n}  w \hi i\| \lo 2.  \\
  {  \|   R \lo {13} \hi j \| \lo F} \leq & { \|   \hat \Sigma \hi n  \lo { \xi \hi j \sf N \hi i \lo n } -   \Sigma \hi n  \lo { \xi \hi j \sf N \hi i \lo n } \| \lo F\|(\hat \Sigma  \hi n \lo { \sf N \hi i  \lo n  \sf N \hi i \lo n } + \lambda \lo n \hat D \lo {\sf N \hi i \lo n }) \hi {-1} \frac{\mathbf{\tilde{H} \lo n  (\xi \hi {\sf N \hi i \lo n })}}{n}  \epsilon \hi i\| \lo 2.}
    \end{align*}  
Next we derive bounds for the probabilities 
on the right hand side of \eqref{erg1}
starting 
 with the term $\|   R \lo {11} \hi j \| \lo F$.
 Observing that 
  \begin{align*}
  (\hat \Sigma \hi n  \lo { \sf N \hi i \lo n  \sf N \hi i \lo n } - \lambda \lo n \hat D \lo {\sf N \hi i \lo n }) \hi {-1} \Sigma \hi n  \lo {\sf N \hi i  \lo n  \sf N \hi i \lo n }
  =
   (\hat \Sigma \hi n  \lo { \sf N \hi i \lo n  \sf N \hi i \lo n } + \lambda \lo n \hat D \lo {\sf N \hi i \lo n }) \hi {-1}( \Sigma \hi n  \lo {\sf N \hi i \lo n  \sf N \hi i \lo n }-\hat \Sigma \hi n  \lo {\sf N \hi i \lo n  \sf N \hi i \lo n })+ (\hat \Sigma \hi n  \lo { \sf N \hi i \lo n  \sf N \hi i \lo n } + \lambda \lo n \hat D \lo {\sf N \hi i \lo n }) \hi {-1} \hat  \Sigma \hi n  \lo {\sf N \hi i \lo n  \sf N \hi i \lo n }.
     \end{align*}
 we have 
  \begin{align}\label{ineq1}
 \| (\hat \Sigma \hi n  \lo { \sf N \hi i \lo n  \sf N \hi i \lo n } + \lambda \lo n \hat D \lo {\sf N \hi i \lo n }) \hi {-1} \Sigma \hi n  \lo {\sf N \hi i \lo n  \sf N \hi i \lo n } \| \lo {2} &\leq   \| (\hat \Sigma \hi n \lo { \sf N \hi i \lo n  \sf N \hi i \lo n } + \lambda \lo n \hat D \lo {\sf N \hi i \lo n }) \hi {-1}(\hat \Sigma \hi n  \lo {\sf N \hi i \lo n  \sf N \hi i \lo n }- \Sigma  \hi n  \lo {\sf N \hi i \lo n  \sf N \hi i \lo n })\| \lo 2 +  \| (\hat \Sigma \hi n  \lo { \sf N \hi i \lo n  \sf N \hi i \lo n } + \lambda \lo n \hat D \lo {\sf N \hi i \lo n }) \hi {-1} \hat \Sigma \hi n  \lo {\sf N \hi i \lo n  \sf N \hi i \lo n } \| \lo 2 \nonumber  \\
& \leq \frac{4}{C \lo {\min}} \| \hat \Sigma \hi n  \lo {\sf N \hi i \lo n  \sf N \hi i \lo n }- \Sigma \hi n  \lo {\sf N \hi i \lo n \sf N \hi i \lo n } \| \lo F +1,
     \end{align}
    where we used the fact that 
      $     \hat \Sigma \hi n  \lo { \sf N \hi i \lo n \sf N \hi i \lo n} 
      \preceq        
      \hat \Sigma \hi n  \lo { \sf N \hi i \lo n \sf N \hi i \lo n} + \lambda \lo n I \lo {\sf N \hi i \lo n}
      \preceq 
\hat \Sigma \hi n  \lo { \sf N \hi i \lo n \sf N \hi i \lo n} + \lambda \lo n \hat D \lo {\sf N \hi i \lo n} $.
          % R11
   From this, it follows that
       \begin{align*}
 \max \lo {j \notin  {\sf N \hi i \lo n}}   \|   R \lo {11} \hi j \| \lo F   \lesssim \max \lo {j \notin  {\sf N \hi i \lo n}} \{  \|   \hat \Sigma \hi n  \lo { \xi \hi j \sf N \hi i \lo n } -   \Sigma \hi n  \lo { \xi \hi j \sf N \hi i \lo n } \| \lo F  \}  \Big( \frac{4}{C \lo {\min}}\| \hat \Sigma \hi n  \lo {\sf N \hi i \lo n  \sf N \hi i \lo n }- \Sigma \hi n  \lo {\sf N \hi i \lo n  \sf N \hi i \lo n } \| \lo F +1  \Big), 
\end{align*}
where we use the fact that $\|B \lo n \hi {* \sf N \hi i \lo n}  \|\lo 2$ is bounded.
Consider the event
$$
\mathcal {A} \lo 1=\{   \max \lo {j \notin  {\sf N \hi i \lo n}} \{ \|   \hat \Sigma \hi n  \lo { \xi \hi j \sf N \hi i \lo n } -   \Sigma \hi n  \lo { \xi \hi j \sf N \hi i \lo n } \| \lo F \} \leq  \frac{c \lo 1 \lambda \lo n  \tsum \lo {l \in  \sf N \hi i \lo n } \|  B \hi {*il} \lo {m \lo n k \lo n} \| \lo F}{6} \}.
$$
Then,
      \begin{align*}
&  P  \Big(  \max \lo {j \notin  {\sf N \hi i \lo n}}  \|   R \lo {11} \hi j \| \lo F \geq \frac{ c \lo 1 \lambda \lo n}{3} \tsum \lo {j \in  \sf N \hi i \lo n } \|  B \hi {*ij} \lo {m \lo n k \lo n} \| \lo F    \Big)  
 \leq    P  \Big(  \| \hat \Sigma \hi n  \lo {\sf N \hi i \lo n  \sf N \hi i \lo n }- \Sigma \hi n  \lo {\sf N \hi i \lo n  \sf N \hi i \lo n } \| \lo F   \geq  \frac{C \lo {\min}}{4}  \Big)+ P  \Big( \mathcal {A} \lo 1 \hi {\complement}  \Big), \\
&   \quad \quad   \lesssim  \exp  \Big(  - C \lo 2 \frac{n \hi {1-\alpha(2 +3  \beta  )}}{  (n \hi i) \hi 2 m \lo n \hi 2 k \lo n \hi 4} + 2 \log ( {n \hi i} m \lo n  k \lo n )  \Big)\\ 
&  \quad \quad    
+ \exp  \Big( - C \lo 2 \frac{n \hi {1-\alpha(2 +3  \beta  )} (\lambda \lo n  \tsum \lo {l \in  \sf N \hi i \lo n } \|  B \hi {*il} \lo {m \lo n k \lo n} \| \lo F) \hi 2}{  {n \hi i}  m \lo n \hi 2 k \lo n \hi 4} +  \log ( {n \hi i} m \lo n \hi 2 k \lo n \hi 2 )+\log ( {(p-n \hi i)} m \lo n \hi 2 k \lo n \hi 2 )   \Big) ,
   \end{align*} 
   where we have used Theorem \ref{theorem1} and the union bound.  
  
       For  the term $R \lo {12} \hi j$, we use the same arguments as for the term $T \lo 2$ in the proof of Proposition \ref{ridgeconsi}. Specifically, recall the definition of  the event $ \Omega$ in  \eqref{omega} and the calculation in \eqref{errorspline} to obtain on the event $\Omega$   
   \begin{align*}
 \max \lo {j \notin  {\sf N \hi i \lo n}}  \| R \lo {12} \hi j \| \lo F & \leq  \max \lo {j \notin  {\sf N \hi i \lo n}} \|   \hat \Sigma \hi n  \lo { \xi \hi j \sf N \hi i  \lo n } -   \Sigma \hi n  \lo { \xi \hi j \sf N \hi i \lo n } \| \lo F \Big(\frac{2}{C \lo {\min}}  \|\hat \Sigma \hi n  \lo {\sf N \hi i \lo n  \sf N \hi i \lo n }- \Sigma \hi n \lo {\sf N \hi i \lo n  \sf N \hi i \lo n } \| \lo F + 1 \Big) \Big(\frac{C \lo {\min}}{2} \Big) \hi {-1/2} c \lo 2 \frac{n \hi i m \lo n \hi {3/2}}{ k \lo n \hi d}\\
& \leq  \max \lo {j \notin  {\sf N \hi i \lo n}} \|   \hat \Sigma \hi n  \lo { \xi \hi j \sf N \hi i  \lo n } -   \Sigma \hi n  \lo { \xi \hi j \sf N \hi i \lo n } \| \lo F \Big(\frac{4}{C \lo {\min}}  \|\hat \Sigma \hi n  \lo {\sf N \hi i \lo n  \sf N \hi i \lo n }- \Sigma \hi n \lo {\sf N \hi i \lo n  \sf N \hi i \lo n } \| \lo F + 1 \Big) \Big(\frac{4}{C \lo {\min}} \Big)  c \lo 2 \lambda \lo n   \tsum \lo {j \in  \sf N \hi i \lo n } \|  B \hi {*ij} \lo {m \lo n k \lo n} \| \lo F 
  \end{align*}  
  for some constant $ c \lo 2>0$, where we used condition $\eqref{as2theorem2}$ for the last inequality.

Then, conditioning on the event $\mathcal{A}\lo 2=\{\max \lo {j \notin  {\sf N \hi i \lo n}}\|   \hat \Sigma \hi n  \lo { \xi \hi j \sf N \hi i \lo n } -   \Sigma \hi n  \lo { \xi \hi j \sf N \hi i \lo n } \| \lo F \geq \frac{c \lo 1 C \lo {\min}}{c \lo 212}\} $, we have by Proposition \ref{splineapproximation}
\begin{equation*}
   \begin{split}
&P\Big(  \max \lo {j \notin  {\sf N \hi i \lo n}} \| R \lo {12} \|  \geq  \frac{c \lo 1}{3}{ \lambda \lo n} \tsum \lo {j \in  \sf N \hi i \lo n } \|  B \hi {*ij} \lo {m \lo n k \lo n} \| \lo F  \Big)
\\  & \quad  \quad \quad \leq P\Big( \|\hat \Sigma \hi n  \lo {\sf N \hi i \lo n  \sf N \hi i \lo n }- \Sigma \hi n \lo {\sf N \hi i \lo n  \sf N \hi i \lo n } \| \lo F \geq \frac{C \lo {\min}}{4} \Big)+ P \Big( \max \lo {j \notin  {\sf N \hi i \lo n}}\|   \hat \Sigma \hi n  \lo { \xi \hi j \sf N \hi i \lo n } -   \Sigma \hi n  \lo { \xi \hi j \sf N \hi i \lo n } \| \lo F \geq \frac{c \lo 1 C \lo {\min}}{c \lo 212}  \Big) +P(\Omega \hi {\complement}) \\
&    \quad \quad \quad \lesssim   \exp \Big( - C \lo {3} \frac{n \hi {1-\alpha(2 +3  \beta  )} }{  (n \hi i) \hi 2 m \lo n \hi 2 k \lo n \hi {4}} + 2 \log ( {n \hi i} m \lo n  k \lo n ) \Big) \\
&      \quad \quad \quad \quad +  \exp \Big( - C \lo {3} \frac{n \hi {1-\alpha(2 +3  \beta  )}}{ {n \hi i}  m \lo n \hi 2 k \lo n \hi 4} +  \log ( {n \hi i} m \lo n  \hi 2 k \lo n \hi 2  )+\log ( {(p-n \hi i)} m \lo n  \hi 2 k \lo n \hi 2  ) \Big) \\
& \quad \quad \quad \quad + \exp (-C \lo 3 \frac{n}{n \hi i m \lo n \hi 2 k \lo n \hi{2d}}+\log(n \hi i m \lo n \hi {2}) ).
\end{split}
\end{equation*}
     % R 13
\noindent
For the term $R \lo {13} \hi j$, we have
  { \begin{align*}
 \max \lo {j \notin  {\sf N \hi i \lo n}} \| R \lo {13} \hi j \| \lo F \leq  \frac{2}{ C \lo {\min}}  \max \lo {j \notin  {\sf N \hi i \lo n}}  \|   \hat \Sigma \hi n  \lo { \xi \hi j \sf N \hi i \lo n } -   \Sigma \hi n  \lo { \xi \hi j \sf N \hi i \lo n } \| \lo F \| \frac{\mathbf{\tilde{H}} \lo n  (\xi \hi {\sf N \hi i \lo n })}{n}  \epsilon \hi i\| \lo F, 
\end{align*}  } 
which yields
\begin{equation*}
\begin{split}
& P \Big(  \max \lo {j \notin  {\sf N \hi i \lo n}} \| R \lo {13} \hi j \| \lo F \geq   \frac{ c \lo 1}{3} { \lambda \lo n} \tsum \lo {j \in  \sf N \hi i \lo n } \|  B \hi {*ij} \lo {m \lo n k \lo n} \| \lo F   \Big)\\
 &  \quad  \quad \leq P \Big( \max \lo {j \notin  {\sf N \hi i \lo n}}  \|   \hat \Sigma \hi n  \lo { \xi \hi j \sf N \hi i \lo n } -   \Sigma \hi n  \lo { \xi \hi j \sf N \hi i \lo n } \| \lo F \geq  \frac{ c \lo 1 \lambda \lo n \tsum \lo {j \in  \sf N \hi i \lo n } \|  B \hi {*ij} \lo {m \lo n k \lo n} \| \lo F}{3}  \Big)+ {P \Big( \| \frac{\mathbf{\tilde{H}} \lo n  (\xi \hi {\sf N \hi i \lo n })}{n}  \epsilon \hi i\| \lo F \geq \frac{C \lo {\min}}{2}  \Big) }\\
  &   \quad  \quad   \lesssim    \exp \Big( - C \lo {4} \frac{n \hi {1-\alpha(2 +3  \beta  )}  (\lambda \lo n \tsum \lo {j \in  \sf N \hi i \lo n } \|  B \hi {*ij} \lo {m \lo n k \lo n} \| \lo F) \hi 2 }{ {n \hi i}  m \lo n \hi 2 k \lo n \hi 4} +  \log ( {n \hi i} m \lo n  \hi 2 k \lo n \hi 2  )+\log ( {(p-n \hi i)} m \lo n  \hi 2 k \lo n \hi 2  ) \Big) \\
  & \quad  \quad \quad  + \exp \Big( - C \lo {4}  \frac{ n  }{n \hi i m \lo n \hi 2 k \lo n}+ \log (n \hi i m \lo n \hi 2 k \lo n)  \Big),  
 \end{split}
\end{equation*}   
where we used Theorem \ref{theorem1} and \eqref{t23erroreq1}.  
 Combining together the results for the terms $\| R \lo {13} \hi j \| $, $\| R \lo {12} \hi j \| $ and $\| R \lo {13} \hi j \| $ we conclude that
\begin{equation} \nonumber 
 P \Big(\max \lo {j \notin {\sf N} \hi i \lo n }   {\| R \lo 1 \hi j \| \lo {F}} \geq    c \lo 1 { \lambda \lo n} \tsum \lo {j \in  \sf N \hi i \lo n } \|  B \hi {*ij} \lo {m \lo n k \lo n} \| \lo F    \Big)
  \lesssim \exp \Big ( -C \lo {2} \frac{n \hi {1 - \alpha (2+3 \beta)}(\lambda \lo n  \tsum \lo {l \in  \sf N \hi i \lo n } \|  B \hi {*il} \lo {m \lo n k \lo n} \| \lo F) \hi 2 }{n \hi i m \hi 2 \lo n k \hi 4 \lo n}  + \log ( pm \hi 2 \lo n k \hi 2  \lo n) \Big ). 
\end{equation}

 %%% R2
 \noindent
{\it {Term $R \lo 2 \hi j$}}. First we write
     \begin{align*}
\max \lo {j \notin {\sf N} \hi i \lo n }  \| R \lo 2  \hi j \| \lo F &\leq  \max \lo {j \notin {\sf N} \hi i \lo n } \big \{ \|  \Sigma \hi n  \lo { \xi \hi j \sf N \hi i \lo n } (\Sigma  \hi n \lo { \sf N \hi i \lo n  \sf N \hi i \lo n }) \hi {-1} \| \lo F  \big \} \| \Sigma \hi n \lo { \sf N \hi i \lo n  \sf N \hi i \lo n } (\hat \Sigma \hi n  \lo { \sf N \hi i \lo n  \sf N \hi i \lo n } + \lambda \lo n \hat D \lo {\sf N \hi i \lo n }) \hi {-1}\| \lo 2  \|\hat {\Sigma} \hi n  \lo { \sf N \hi i \lo n  {\xi \hi i}}- {\Sigma} \hi n  \lo { \sf N \hi i \lo n  {\xi \hi i}} \| \lo F\\
& \leq \frac{1-\eta}{\sqrt{n \hi i}}  \| \Sigma \hi n  \lo { \sf N \hi i \lo n  \sf N \hi i \lo n } (\hat \Sigma \hi n  \lo { \sf N \hi i \lo n  \sf N \hi i \lo n } + \lambda \lo n \hat D \lo {\sf N \hi i \lo n }) \hi {-1}\| \lo 2  \|\hat {\Sigma}  \hi n \lo { \sf N \hi i \lo n  {\xi \hi i}}- {\Sigma} \hi n  \lo { \sf N \hi i \lo n  {\xi \hi i}} \| \lo F\\
& \leq \frac{1-\eta}{\sqrt{n \hi i}} (\frac{4}{C \lo {\min}}\|\hat \Sigma \hi n  \lo {\sf N \hi i \lo n  \sf N \hi i \lo n }- \Sigma \hi n  \lo {\sf N \hi i \lo n  \sf N \hi i \lo n } \| \lo F + 1 )   \|\hat {\Sigma}  \hi n \lo { \sf N \hi i \lo n  {\xi \hi i}}- {\Sigma} \hi n  \lo { \sf N \hi i \lo n  {\xi \hi i}} \| \lo F,
     \end{align*}
where we used \eqref{ineq1} and the second inequality holds with high probability by   Lemma \ref{sampleirrepr}.
As a result, 
it follows that
   \begin{equation}
   \nonumber 
   \begin{split}
& P \Big(\max \lo {j \notin {\sf N \hi i \lo n }} \{ {\| R \lo 2 \| \lo {F}} \geq    c \lo 1 \lambda \lo n \tsum \lo {j \in  \sf N \hi i \lo n } \|  B \hi {*ij} \lo {m \lo n k \lo n} \| \lo F  \} \Big) \\
& \quad \quad \quad  \leq P \Big (  \frac{1-\eta}{\sqrt{n \hi i}}  (\frac{4}{C \lo {\min}}\|\hat \Sigma \hi n  \lo {\sf N \hi i \lo n  \sf N \hi i \lo n }- \Sigma \hi n  \lo {\sf N \hi i \lo n  \sf N \hi i \lo n } \| \lo F + 1 )  \|\hat {\Sigma} \hi n  \lo { \sf N \hi i \lo n  {\xi \hi i}}- {\Sigma} \hi n  \lo { \sf N \hi i \lo n  {\xi \hi i}}  \| \lo F \geq  c \lo 1 \lambda \lo n \tsum \lo {j \in  \sf N \hi i \lo n } \|  B \hi {*ij} \lo {m \lo n k \lo n} \| \lo F  \Big ) \\
& \quad \quad \quad \leq   P \Big(  \| \hat \Sigma \hi n  \lo {\sf N \hi i \lo n \sf N \hi i \lo n}- \Sigma \hi n  \lo {\sf N \hi i \lo n \sf N \hi i \lo n} \| \lo F   \geq  \frac{C \lo {\min}}{4} \Big) 
    +P (\mathcal {A} \lo 2 \hi {\complement})  \\
   & \quad \quad \quad \lesssim   \exp  \left( - C \lo {3} \frac{n \hi {1-\alpha(2 +3  \beta  )}  }{  (n \hi i) \hi 2 m \lo n \hi 2  k \lo n \hi {4}}  + 2 \log ( {n \hi i} m \lo n  k \lo n )  \right)   \\
 & \quad \quad \quad \quad \quad + \exp  \Big(- C \lo {3} \frac{n \hi {1-\alpha(2 +3  \beta  )} (\lambda \lo n \tsum \lo {j \in  \sf N \hi i \lo n } \| B \hi {*ij} \lo {m \lo n k \lo n} \| \lo F) \hi 2}{ m \lo n \hi 2 k \lo n \hi {3}  } 
 +  \log ( n \hi i  m \lo n \hi 2 k \lo n )\Big),
 \end{split}
 \end{equation}
 where
 \begin{align*}
\mathcal{A} \lo 2= \Big \{ \max \lo {j \notin  {\sf N \hi i \lo n }} \|\hat {\Sigma} \hi n  \lo { \sf N \hi i \lo n  {\xi \hi i}}- {\Sigma} \hi n  \lo { \sf N \hi i \lo n  {\xi \hi i}}  \| \lo F \leq \frac{c \lo 1}{2} \frac{\sqrt{n \hi i}}{1-\eta  } \lambda \lo n \tsum \lo {j \in  \sf N \hi i \lo n } \|  B \hi {*ij} \lo {m \lo n k \lo n} \| \lo F \Big \}, 
 \end{align*}
 and we used  Theorem \ref{theorem1} and the union bound for the last inequality.
 
 %%% R3
             {\it {Term $R \lo 3 \hi j$}}. 
       Using the identity $A \hi {-1} - B \hi {-1}=A \hi {-1} (B-A ) B \hi {-1} $ and following similar arguments used to obtain bounds  for the terms $R \lo 1$ and $R \lo 2$, we  get (note that we are working on the event ${\cal N}$)
       \begin{align*}
       P \Big(\max \lo {j \notin {\sf N \hi i \lo n }}  {\| R \lo 3 \hi j \| \lo {F}} \geq  c \lo 1{ \lambda \lo n} \tsum \lo {j \in  \sf N \hi i \lo n } \|  B \hi {*ij} \lo {m \lo n k \lo n} \| \lo F   \Big) \leq \tsum \lo {r=1} \hi {3} P \Big ( \max \lo {j \notin {\sf N \hi i  }} {\| R \hi j \lo {3r} \| \lo {F}} \geq  \frac{c \lo 1}{3}{ \lambda \lo n} \tsum \lo {j \in  \sf N \hi i \lo n } \|  B \hi {*ij} \lo {m \lo n k \lo n} \| \lo F   \Big ),
       \end{align*}
       where,
\begin{align*}
\| R \lo {31} \hi j \| \lo F & \leq  \frac{1- \eta}{\sqrt{n \hi i}} \Big(1+\frac{4}{C \lo {\min}}  \| \Sigma \hi n  \lo { \sf N \hi i \lo n  \sf N \hi i \lo n } - \hat  \Sigma \hi n  \lo { \sf N \hi i \lo n  \sf N \hi i \lo n } \| \lo {F}\Big)    \| \Sigma \hi n  \lo { \sf N \hi i \lo n  \sf N \hi i \lo n } - \hat  \Sigma \hi n  \lo { \sf N \hi i \lo n  \sf N \hi i \lo n } \| \lo {F} \|(  \Sigma \hi n  \lo { \sf N \hi i \lo n  \sf N \hi i \lo n } + \lambda \lo n \hat D \lo {\sf N \hi i \lo n }) \hi {-1}  \Sigma \hi n  \lo {\sf N \hi i \lo n  \sf N \hi i \lo n }B \hi {* \sf N \hi i \lo n } \lo n \| \lo F,\\
\| R \lo {32} \hi j \| \lo F & \leq \frac{1- \eta}{\sqrt{n \hi i}} \Big(1+\frac{4}{C \lo {\min}}   \| \Sigma \hi n  \lo { \sf N \hi i \lo n  \sf N \hi i \lo n } - \hat  \Sigma \hi n  \lo { \sf N \hi i \lo n  \sf N \hi i \lo n } \| \lo {F}\Big)  \| \Sigma \hi n  \lo { \sf N \hi i \lo n  \sf N \hi i \lo n } -    \hat \Sigma \hi n  \lo { \sf N \hi i \lo n  \sf N \hi i \lo n } \| \lo {F}   \| (  \Sigma \hi n  \lo { \sf N \hi i \lo n  \sf N \hi i \lo n } + \lambda \lo n \hat D \lo {\sf N \hi i \lo n }) \hi {-1}   \frac{\mathbf{\tilde{H}} \lo n  (\xi \hi {\sf N \hi i \lo n })}{n}  w \hi i \| \lo F, \\
\| R \lo {33} \hi j \| \lo F& \leq \frac{1- \eta}{\sqrt{n \hi i}}  \Big(1+\frac{4}{C \lo {\min}}  \| \Sigma \hi n  \lo { \sf N \hi i \lo n  \sf N \hi i \lo n } - \hat  \Sigma \hi n  \lo { \sf N \hi i \lo n  \sf N \hi i \lo n } \| \lo {F}\Big)  \| \Sigma \hi n  \lo { \sf N \hi i \lo n  \sf N \hi i \lo n } - \hat  \Sigma \hi n  \lo { \sf N \hi i \lo n  \sf N \hi i \lo n } \| \lo {F} { \| (  \Sigma \hi n \lo { \sf N \hi i \lo n  \sf N \hi i \lo n } + \lambda \lo n \hat D \lo {\sf N \hi i \lo n }) \hi {-1}   \frac{\mathbf{\tilde{H}} \lo n  (\xi \hi {\sf N \hi i \lo n })}{n}  \epsilon \hi i \| \lo F},
\end{align*}
(note that the terms on the right-hand side are independent of $j$).  We next derive bounds for the probabilities $P \left( {\| R \lo {3r} \hi j \| \lo {F}} \geq  \frac{c \lo 1}{3}{ \lambda \lo n} \tsum \lo {j \in  \sf N \hi i \lo n } \|  B \hi {*ij} \lo {m \lo n k \lo n} \| \lo F   \right)$ for $r=1,2,3$. For the term $\| R \lo {31} \| \lo F$, note that
 \begin{align*}
\| R \lo {31} \hi j \| \lo F & \lesssim \frac{(1- \eta)}{\sqrt{n \hi i}} \| \Sigma \hi n  \lo { \sf N \hi i \lo n  \sf N \hi i \lo n } - \hat  \Sigma \hi n  \lo { \sf N \hi i \lo n  \sf N \hi i \lo n } \| \lo {F}  \Big( \frac{4}{C \lo {\min}} \| \hat \Sigma \hi n \lo {\sf N \hi i \lo n  \sf N \hi i \lo n }- \Sigma \hi n  \lo {\sf N \hi i \lo n  \sf N \hi i \lo n } \| \lo F +1 \Big)~,
\end{align*}
which implies (using the same arguments as before)
\begin{equation}\label{R31}
\begin{split}
&P \Big (\max \lo {j \notin \sf N^i} {\| R \lo {31} \hi j \| \lo {F}} \geq  \frac{c \lo 1}{3} { \lambda \lo n} \tsum \lo {j \in  \sf N \hi i \lo n } \|  B \hi {*ij} \lo {m \lo n k \lo n} \| \lo F   \Big )\\ 
& \quad \leq P \Big( \| \Sigma \hi n  \lo { \sf N \hi i \lo n  \sf N \hi i \lo n } - \hat  \Sigma \hi n  \lo { \sf N \hi i \lo n  \sf N \hi i \lo n } \| \lo {F} \geq  \frac{ \sqrt{n \hi i} \lambda \lo n \tsum \lo {j \in  \sf N \hi i \lo n } \|  B \hi {*ij} \lo {m \lo n k \lo n} \| \lo F  }{1-\eta} \frac{c \lo 1 }{6} \Big) +P \Big( \| \Sigma \hi n  \lo { \sf N \hi i \lo n  \sf N \hi i \lo n } - \hat  \Sigma \hi n  \lo { \sf N \hi i \lo n  \sf N \hi i \lo n } \| \lo {F} \geq  \frac{C \lo {\min}}{4} \Big)\\
& \quad \lesssim    \exp  \Big( - C \lo {4} \frac{n \hi {1-\alpha(2 +3  \beta  )} (\lambda \lo n \tsum \lo {j \in  \sf N \hi i \lo n } \|  B \hi {*ij} \lo {m \lo n k \lo n} \| \lo F ) \hi 2 }{  n \hi i  m \lo n \hi 2  k \lo n \hi {4}}  + 2 \log ( {n \hi i} m \lo n  k \lo n )\Big) \nonumber 
\\
& \quad +   \exp  \Big( - C \lo {4} \frac{n \hi {1-\alpha(2 +3  \beta  )} }{  (n \hi i) \hi 2 m \lo n \hi 2  k \lo n \hi {4}}  + 2 \log ( {n \hi i} m \lo n  k \lo n )\Big)  ,
\end{split}
 \end{equation}
for some positive constant $C \lo 4$ that depends on $\eta$ and $C \lo {\min}$.
%R32
Similarly, for the term $\| R \lo {32} \hi j \| \lo F$ have
\begin{align*}
\| R \lo {32} \hi j \| \lo F \leq  \frac{1- \eta}{\sqrt{n \hi i}}   \| \Sigma \hi n  \lo { \sf N \hi i \lo n  \sf N \hi i \lo n } - \hat  \Sigma \hi n  \lo { \sf N \hi i \lo n  \sf N \hi i \lo n } \| \lo {F}  (\frac{4}{C \lo {\min}} \| \hat \Sigma \hi n \lo {\sf N \hi i \lo n  \sf N \hi i \lo n }- \Sigma \hi n  \lo {\sf N \hi i \lo n  \sf N \hi i \lo n } \| \lo F +1 )  \| \frac{w \hi i}{\sqrt{n}} \| \lo F.
\end{align*}
 Hence, recalling the definition of the set $\Omega$ in \eqref{omega}  and using conditions \eqref{as2theorem2} and \eqref{m1}, Proposition \ref{splineapproximation} and Theorem \ref{theorem1},   we obtain for some positive constant
 \begin{equation*}
 \begin{split}
&P \Big( \max \lo {j \notin \sf N \hi i}  \| R \lo {32} \hi j \| \lo {F} \geq  \frac{c \lo 1}{3} { \lambda \lo n} \tsum \lo {j \in  \sf N \hi i \lo n } \|  B \hi {*ij} \lo {m \lo n k \lo n} \| \lo F   \Big) \\
& \quad \leq P \Big( \| \Sigma \hi n  \lo { \sf N \hi i \lo n  \sf N \hi i \lo n } - \hat  \Sigma \hi n  \lo { \sf N \hi i \lo n  \sf N \hi i \lo n } \| \lo {F} \geq  \frac{ \sqrt{n \hi i} c \lo 1    }{12 c\lo 2(1-\eta)} \Big) +P \Big( \| \Sigma \hi n  \lo { \sf N \hi i \lo n  \sf N \hi i \lo n } - \hat  \Sigma  \hi n \lo { \sf N \hi i \lo n  \sf N \hi i \lo n } \| \lo {F} \geq  \frac{C \lo {\min}}{4} \Big) + P(\Omega \hi {\complement})\\
& \quad\lesssim   \exp  \Big( - C \lo {5} \frac{n \hi {1-\alpha(2 +3  \beta  )} }{  n \hi i  m \lo n \hi 2  k \lo n \hi {4}}  + 2 \log ( {n \hi i} m \lo n  k \lo n )\Big)+   \exp  \Big( - C \lo {5} \frac{n \hi {1-\alpha(2 +3  \beta  )}  }{  (n \hi i) \hi 2 m \lo n \hi 2  k \lo n \hi {4}}  + 2 \log ( {n \hi i} m \lo n  k \lo n )\Big)   \\
& \quad \quad + \exp (-C \lo 3 \frac{n}{n \hi i m \lo n \hi 2 k \lo n \hi{2d}}+\log(n \hi i m \lo n \hi {2}) ),
\end{split}
 \end{equation*}
for some positive constant $C \lo {5}$ that depends on $\eta$ and $C \lo {\min}$.
%R33
For the term $\| R \lo {33} \hi j \| \lo F$, we write
 \begin{equation*}
 \begin{split}
{\| R \lo {33} \hi j \| \lo F  \leq \frac{(1-\eta)}{\sqrt{n \hi i}} (\frac{2}{C \lo {\min}}) \hi {1/2}  \Big(1+\frac{4}{C \lo {\min}}  \| \Sigma \hi n  \lo { \sf N \hi i \lo n  \sf N \hi i \lo n } - \hat  \Sigma \hi n  \lo { \sf N \hi i \lo n  \sf N \hi i \lo n } \| \lo {F}\Big)   \| \Sigma \hi n  \lo { \sf N \hi i \lo n  \sf N \hi i \lo n } - \hat  \Sigma \hi n  \lo { \sf N \hi i \lo n  \sf N \hi i \lo n } \| \lo {F}   \|    \frac{\mathbf{\tilde{H} \lo n}  (\xi \hi {\sf N \hi i \lo n })}{n}  \epsilon \hi i \| \lo F,}
\end{split}
 \end{equation*}
and consider the event $ \mathcal {A} \lo 3 =\{  \| \Sigma \hi n  \lo { \sf N \hi i \lo n  \sf N \hi i \lo n } - \hat  \Sigma \hi n  \lo { \sf N \hi i \lo n  \sf N \hi i \lo n } \| \lo {F}  \leq \frac{c \lo 1 \sqrt{n \hi i}}{ 1-\eta} \lambda \lo n   \tsum \lo {j \in  \sf N \hi i \lo n } \|  B \hi {*ij} \lo {m \lo n k \lo n} \| \lo F  \}$. Then
 \begin{equation*}
 \begin{split}
&P \Big( \max \lo {j \notin \sf N \hi i} \| R \lo {33} \hi j \| \lo {F} \geq   \frac{c \lo 1}{3} { \lambda \lo n} \tsum \lo {j \in  \sf N \hi i \lo n } \|  B \hi {*ij} \lo {m \lo n k \lo n} \| \lo F \Big) \\
& \quad \quad  \leq  P \Big  ( \max \lo {j \notin \sf N \hi i}  \| R \lo {33} \hi j \| \lo {F} \geq   \frac{c \lo 1}{3} { \lambda \lo n} \tsum \lo {j \in  \sf N \hi i \lo n } \|  B \hi {*ij} \lo {m \lo n k \lo n} \| \lo F \text{ and } \mathcal {A} \lo 3  \Big ) 
+P(\mathcal {A} \lo 3 \hi {\complement})\\
& \quad \quad  \leq P \Big(  \|    \frac{\mathbf{\tilde{H} \lo n}  (\xi \hi {\sf N \hi i \lo n })}{n}  \epsilon \hi i \| \lo F
\geq  \frac{\sqrt{C \lo {\min}}}{6\sqrt{2}} \Big) +P \Big( \| \Sigma \hi n  \lo { \sf N \hi i \lo n  \sf N \hi i \lo n } - \hat  \Sigma  \hi n \lo { \sf N \hi i \lo n  \sf N \hi i \lo n } \| \lo {F} \geq  \frac{C \lo {\min}}{4} \Big) +P(\mathcal {A} \lo 3 \hi {\complement})\\
& \quad  \quad \lesssim    \exp  \Big( - C \lo {6} \frac{n \hi {1-\alpha(2 +3  \beta)} (\lambda \lo n \tsum \lo {j \in  \sf N \hi i \lo n } \|  B \hi {*ij} \lo {m \lo n k \lo n} \| \lo F) \hi 2 }{  n \hi i  m \lo n \hi 2  k \lo n \hi {4}}  + 2 \log ( {n \hi i} m \lo n  k \lo n )\Big)  \\
& \quad \quad  
+   \exp  \Big( - C \lo {6} \frac{n \hi {1-\alpha(2 +3  \beta  )}  }{  (n \hi i) \hi 2 m \lo n \hi 2  k \lo n \hi {4}}  + 2 \log ( {n \hi i} m \lo n  k \lo n )\Big)
+ \exp \Big( - C \lo {6}  \frac{n    }{n \hi i m \lo n \hi 2 k \lo n}+ \log (n \hi i m \lo n \hi 2 k \lo n)  \Big)  ,
\end{split}
 \end{equation*}
where we used \eqref{t23erroreq1} and Theorem \ref{theorem1}.
Combining    the results of $\| R \lo {31} \hi j \| \lo {F}$, $\| R \lo {32} \hi j \| \lo {F}$ and $\| R \lo {33} \hi j \| \lo {F}$ we conclude,
\begin{equation*}
P \Big( \max \lo {j \notin {\sf N \hi i \lo n }}  {\| R \lo {3} \hi j \| \lo {F}} \geq    c \lo 1  { \lambda \lo n} \tsum \lo {j \in  \sf N \hi i \lo n } \|  B \hi {*ij} \lo {m \lo n k \lo n} \| \lo F  \Big)\\
% & \quad \quad \quad 
\lesssim    \exp  \Big( - C \lo {4} \frac{n \hi {1-\alpha(2 +3  \beta)} (\lambda \lo n \tsum \lo {j \in  \sf N \hi i \lo n } \|  B \hi {*ij} \lo {m \lo n k \lo n} \| \lo F) \hi 2 }{  n \hi i  m \lo n \hi 2  k \lo n \hi {4}}  + 2 \log ( {n \hi i} m \lo n  k \lo n )\Big) .
 \end{equation*}
% with the first term dominating the second term.\\
 %%% R4
{\it {Term $R \lo 4 \hi j$}}. We have  
 \begin{align*}
 \max \lo {j \notin {\sf N \hi i \lo n }} \| R \lo 4 \hi j \| \lo {F}& \leq  \max \lo {j \notin {\sf N \hi i \lo n }}  \|  \hat \Sigma \hi n  \lo { \xi \hi j \sf N \hi i \lo n } -   \Sigma \hi n  \lo { \xi \hi j \sf N \hi i \lo n } \| \lo {F}  \|  (\hat \Sigma \hi n  \lo { \sf N \hi i \lo n  \sf N \hi i \lo n } + \lambda \lo n \hat D \lo {\sf N \hi i \lo n }) \hi {-1}  \|  \lo {2}   \| \hat {\Sigma} \hi n  \lo { \sf N \hi i \lo n  {\xi \hi i}}- {\Sigma} \lo { \sf N \hi i \lo n  {\xi \hi i}}  \|  \lo F \\
& \leq \frac{2}{C \lo {\min}} \max \lo {j \notin {\sf N \hi i \lo n }} \|  \hat \Sigma \lo { \xi \hi j \sf N \hi i \lo n } -   \Sigma \hi n  \lo { \xi \hi j \sf N \hi i \lo n } \| \lo {F}   \| \hat {\Sigma} \lo { \sf N \hi i \lo n  {\xi \hi i}}- {\Sigma} \hi n  \lo { \sf N \hi i \lo n {\xi \hi i}}  \|  \lo F.
          \end{align*}
Recall the definition of event $\mathcal {A} \lo 1=\{  \max \lo {j \notin {\sf N \hi i \lo n }}   \|   \hat \Sigma \hi n  \lo { \xi \hi j \sf N \hi i \lo n } -   \Sigma \hi n  \lo { \xi \hi j \sf N \hi i \lo n } \| \lo F \leq  \lambda \lo n \tsum \lo {j \in  \sf N \hi i \lo n } \|  B \hi {*ij} \lo {m \lo n k \lo n} \| \lo F  \}$, then we have
 \begin{align*}
 & P \Big( \max \lo {j \notin {\sf N \lo n  \hi i}} \| R \hi j \lo 4 \| \lo {F} \geq   c \lo 1{ \lambda \lo n} \tsum \lo {j \in  \sf N \hi i \lo n } \|  B \hi {*ij} \lo {m \lo n k \lo n} \| \lo F  \Big) 
 \leq    P \Big( \| \hat {\Sigma} \hi n \lo { \sf N \hi i \lo n {\xi \hi i}}- {\Sigma} \hi n  \lo { \sf N \hi i \lo n  {\xi \hi i}}  \|  \lo F \geq   c \lo 1 \frac{C \lo {\min}}{2}    \Big) +  P (\mathcal {A} \lo 1 \hi {\complement}) \\
 &  \quad  \quad  \lesssim   \exp  \Big (- C \lo {5} \frac{n \hi {1-\alpha(2+3  \beta  )}}{n \hi i   {m \lo n} \hi 2 k \lo n \hi {3}  } 
 +  \log ( n \hi i  m \lo n \hi 2 k \lo n )  \Big )  \\
 &   \quad \quad   +  \exp \Big  ( - C \lo {5} \frac{n \hi {1-\alpha(2 +3  \beta  )} (\lambda \lo n \tsum \lo {j \in  \sf N \hi i \lo n } \|  B \hi {*ij} \lo {m \lo n k \lo n} \| \lo F  ) \hi 2}  {n \hi i  m \lo n \hi 2 k \lo n \hi {4}} 
 + \log \left( p  m \lo n \hi 2  k \lo n \hi 2 \right)\Big )  , 
   \end{align*}
for some positive constant $C \lo 5$, where the estimates follow from Theorem  \ref{theorem1}.
 
\noindent
  %%% R5
{\it {Term $R \lo 5 \hi j$}}.
 For  the term $R \lo 5 \hi j$ we obtain  by the same argument for some positive constant $C \lo 6$
  \begin{equation}
  \label{R5}
  \begin{split}
&  P  \Big( \max \lo {j \notin {\sf N \hi i \lo n }} \| R \lo 5 \hi j \| \lo {F} \geq    c \lo 1 { \lambda \lo n} \tsum \lo {j \in  \sf N \hi i \lo n } \|  B \hi {*ij} \lo {m \lo n k \lo n} \| \lo F   \Big) 
 \leq  (p- n \hi i )P  \Big( \| {\Sigma} \hi n  \lo { \xi \hi j {\xi \hi i}} -\hat {\Sigma} \hi n  \lo { \xi \hi j {\xi \hi i}}\| \lo F \geq   c \lo 1 \lambda \lo n \tsum \lo {j \in  \sf N \hi i \lo n } \| B \hi {*ij} \lo {m \lo n k \lo n} \| \lo F  
     \Big) \\
   & \quad \quad \quad \lesssim \exp  \Big(- C \lo {6} \frac{n \hi {1-\alpha(2+ 3 \beta  )}( \lambda \lo n \tsum \lo {j \in  \sf N \hi i \lo n } \| B \hi {*ij} \lo {m \lo n k \lo n} \| \lo F) \hi 2}{ m \lo n \hi 2 k \lo n \hi {3}} 
   %+   \log (  m \lo n \hi 2 k \lo n)
   + \log (p   m \lo n \hi 2 k \lo n )   \Big) ~. 
     \end{split}
     \nonumber 
  \end{equation}

  \noindent
    %%% R6
    {\it {Term $R \lo 6 \hi j$}}.  Using again the  identity $A \hi {-1} - B \hi {-1}=A \hi {-1} (B-A ) B \hi {-1}$ and \eqref{sigmazetai}, we have 
   {  \begin{align*}
 R \lo 6 \hi j =\lambda \lo n \Sigma \hi n  \lo { \xi \hi j \sf N \hi i \lo n } ( \Sigma \hi n  \lo { \sf N \hi i \lo n  \sf N \hi i \lo n } + \lambda \lo n  D \hi * \lo {\sf N \hi i \lo n }) \hi {-1} (D \hi * \lo {\sf N \hi i \lo n } -\hat D  \lo {\sf N \hi i \lo n })  ( \Sigma \hi n  \lo { \sf N \hi i \lo n  \sf N \hi i \lo n } + \lambda \lo n  \hat D  \lo {\sf N \hi i \lo n }) \hi {-1}  \Big( \Sigma \hi n \lo {\sf N \hi i \lo n  \sf N \hi i \lo n }B \hi {* \sf N \hi i \lo n } \lo n+ \frac{\mathbf{\tilde{H} \lo n}(\xi \hi {\sf N \hi i \lo n })}{n}  w \hi i +  \frac{\mathbf{\tilde{H} \lo n }(\xi \hi {\sf N \hi i \lo n })}{n}  \epsilon \hi i  \Big).
          \end{align*}}
Obviously, on the event ${\cal N}$ we have 
 \begin{align*}
       P \Big(\max \lo {j \notin {\sf N  \hi i \lo n }}  {\| R \lo 6 \hi j \| \lo {F}} \geq   c \lo 1 { \lambda \lo n} \tsum \lo {j \in  \sf N \hi i \lo n } \|  B \hi {*ij} \lo {m \lo n k \lo n} \| \lo F   \Big) \leq \tsum \lo {r=1} \hi {3} P \Big ( {\| R \lo {3r} \| \lo {F}} \geq   \frac{c \lo 1}{3} { \lambda \lo n} \tsum \lo {j \in  \sf N \hi i \lo n } \|  B \hi {*ij} \lo {m \lo n k \lo n} \| \lo F   \Big ),
       \end{align*}
  where
     \begin{align*}
\| R \lo {61} \hi j \| \lo {F} & \leq \frac{ 1- \eta}{\sqrt{n \hi i}} \lambda \lo n  \| D \hi * \lo {\sf N \hi i \lo n } -\hat D  \lo {\sf N \hi i \lo n }  \| \lo {2}, \\
\| R \lo {62} \hi j \| \lo {F} &\leq   \frac{\ 1- \eta}{\sqrt{n \hi i}}\lambda \lo n  \| D \hi * \lo {\sf N \hi i \lo n } -\hat D  \lo {\sf N \hi i \lo n }  \| \lo {2} \| (\Sigma \hi n  \lo {\sf N \hi i \lo n  \sf N \hi i \lo n }+ \lambda \lo n \hat D  \lo {\sf N \hi i \lo n } ) \inv  \frac{\mathbf{\tilde{H} \lo n(\xi \hi {\sf N \hi i \lo n })}}{n}  w \hi i \| \lo F , \\
\| R \lo {63} \hi j \|  \lo {F} & \leq  \frac{ 1- \eta}{\sqrt{n \hi i}} \lambda \lo n   \| D \hi * \lo {\sf N \hi i \lo n } -\hat D  \lo {\sf N \hi i \lo n }  \| \lo {2} {\| (\Sigma \hi n  \lo {\sf N \hi i \lo n  \sf N \hi i \lo n }+ \lambda \lo n \hat D  \lo {\sf N \hi i \lo n } ) \inv   \frac{\mathbf{\tilde{H} \lo n (\xi \hi {\sf N \hi i \lo n })}}{n}  \epsilon \hi i \| \lo F.}
\end{align*}

Now, using the same computations as in the proof of Proposition 4 in \cite{leefaro}, we obtain for the operator norm of the matrix $( D \hi {*} \lo {\sf N \hi i \lo n }) \lo {jj} -(\hat D \lo {\sf N \hi i \lo n }) \lo {jj}$  
  \begin{align*}
 \| ( D \hi {*} \lo {\sf N \hi i \lo n }) \lo {jj} -(\hat D \lo {\sf N \hi i \lo n }) \lo {jj} \|  \lo 2 &=\Big |  \frac{\tsum \lo {\ell   \in  \sf N \hi i \lo n } \| \hat B \hi {i \ell} \lo n \| \lo F}{\| \hat B \hi {ij} \lo n \| \lo F} - \frac{\tsum \lo {\ell \in  \sf N \hi i \lo n } \|  B \hi {*i \ell} \lo {m \lo n k \lo n} \| \lo F}{\|  B \hi {*ij} \lo {m \lo n k \lo n} \| \lo F} \Big  | \\
 & \leq  \frac{\|  B \hi {*ij} \lo {m \lo n k \lo n} \| \lo F  \tsum \lo {\ell \in  \sf N \hi i \lo n } \left| \| \hat B \hi {i \ell} \lo n \| \lo F  - \|  B \hi {*i \ell} \lo {m \lo n k \lo n} \| \lo F \right|+  \left|\| B \hi {*i j} \lo {m \lo n k \lo n} \| \lo F-\| \hat B \hi {ij} \lo n \| \lo F\right| \tsum \lo {\ell \in  \sf N \hi i \lo n } \|  B \hi {*i\ell} \lo {m \lo n k \lo n} \| \lo F }{\| \hat B \hi {ij} \lo n \| \lo F \|  B \hi {*i j} \lo {m \lo n k \lo n} \| \lo F} \\
&\leq  \frac{\|  B \hi {*ij} \lo {m \lo n k \lo n} \| \lo F  \tsum \lo {\ell \in  \sf N \hi i \lo n }  \| \hat B \hi {i \ell} \lo n   -  B \hi {*i \ell} \lo n \| \lo F +  \| B \hi {*ij} \lo {m \lo n k \lo n} - \hat B \hi {ij} \lo n \| \lo F  \tsum \lo {\ell \in  \sf N \hi i \lo n } \|  B \hi {*i \ell} \lo {m \lo n k \lo n} \| \lo F }{\| \hat B \hi {ij} \lo n \| \lo F \|  B \hi {*ij} \lo {m \lo n k \lo n} \| \lo F} \\
  &\leq  \frac{  \tsum \lo {\ell \in  \sf N \hi i \lo n }  \| \hat B \hi {i \ell} \lo n   -  B \hi {*i \ell} \lo {m \lo n k \lo n} \| \lo F  \tsum \lo {\ell \in  \sf N \hi i \lo n } \|  B \hi {*i \ell} \lo {m \lo n k \lo n} \| \lo F }{ \|  B \hi {*ij} \lo {m \lo n k \lo n} \| \lo F (\|  B \hi {*ij} \lo {m \lo n k \lo n} \| \lo F-   \|  \hat B \hi {ij} \lo n   -   B \hi {*ij} \lo {m \lo n k \lo n} \| \lo F )}.
   \end{align*}
On the event $\mathcal{A} \lo 4=\{ \| \hat B \hi {\sf N \hi i \lo n } \lo n -  B \hi {*\sf N \hi i \lo n } \lo { n} \| \lo {F} \leq \frac{ b \hi {*i} \lo n}{2} \}$,     we obtain
   \begin{align*}
  \| ( D \hi {*} \lo {\sf N \hi i \lo n }) \lo {jj} -(\hat D \lo {\sf N \hi i \lo n }) \lo {jj} \|  \lo 2  & \leq   2 (b \hi {*i} \lo n) \hi {-2}  \tsum \lo {\ell \in  \sf N \hi i \lo n }  \| \hat B \hi {i \ell} \lo n   -  B \hi {*i \ell} \lo {m \lo n k \lo n} \| \lo F  \tsum \lo {\ell \in  \sf N \hi i \lo n } \|  B \hi {*i \ell} \lo {m \lo n k \lo n} \| \lo F \\
  & \leq   2 (b \hi {*i} \lo n) \hi {-2} (n \hi i) \hi {1/2} \| \hat B \hi { \sf N \hi i \lo n } \lo n   -  B \hi {* \sf N \hi i \lo n } \lo n \| \lo F   \tsum \lo {\ell \in  \sf N \hi i \lo n } \|  B \hi {*i \ell} \lo {m \lo n k \lo n} \| \lo F,
   \end{align*}
   where we used the Cauchy-Schwarz inequality.
Consequently,  
 \begin{equation}
% \begin{split}
P \Big( \max \lo {j \notin \sf N \hi i} \| R \lo {61} \hi j \| \lo {F} \geq \frac{c 
\lo 1}{3} { \lambda \lo n} \tsum \lo {j \in  \sf N \hi i \lo n } \|  B \hi {*ij} \lo {m \lo n k \lo n } \| \lo F  \Big)
\leq P (\| \hat B \hi { \sf N \hi i \lo n } \lo n   -  B \hi {* \sf N \hi i \lo n } \lo n \| \lo F  \geq \frac{c \lo 1}{6} \frac{1}{1-\eta} (b \hi {*i} \lo n) \hi {2})+ P (\mathcal {A} \lo 4 \hi {\complement}).
%\end{split}
\nonumber
\end{equation}

Now,  we apply  Proposition  \ref{propconsistency} with $\delta= \frac{c \lo 1}{6} \frac{1}{1-\eta} (b \hi {*i} \lo n) \hi {2}$ and $\delta = \frac{b^{*i}_n}{2}$  to obtain
 \begin{equation*}\label{R61}
 \begin{split}
&P \Big( \max \lo {j \notin \sf N \hi i} \| R \lo {61} \hi j \| \lo {F} \geq \frac{c \lo 1}{3} { \lambda \lo n} \tsum \lo {j \in  \sf N \hi i \lo n } \|  B \hi {*ij} \lo {m \lo n k \lo n} \| \lo F  \Big)  \\
&\quad \quad \lesssim   \exp \Big( - C \lo {7} \frac{n \hi {1- \alpha (2+3\beta)} (b \hi {*i} \lo n)  \hi {6} }{(n \hi i) \hi {4} m \lo n \hi 2  k \lo n \hi {4} (\tsum \lo {j \in  \sf N \hi i \lo n } \|  B \hi {*ij} \lo {m \lo n k \lo n} \| \lo F) \hi 2}+2 \log (n \hi i m \lo n k \lo n)  \Big)  \\
&\quad \quad  \quad + \exp \Big( - C \lo {7} \frac{n \hi {1- \alpha (2+3\beta)} (b \hi {*i} \lo n)  \hi {4} }{(n \hi i) \hi {4} m \lo n \hi 2  k \lo n \hi {4} (\tsum \lo {j \in  \sf N \hi i \lo n } \|  B \hi {*ij} \lo {m \lo n k \lo n} \| \lo F) \hi 2}+2 \log (n \hi i m \lo n k \lo n)  \Big) \\
&\quad \quad  \lesssim \exp \Big( - C \lo {7} \frac{n \hi {1- \alpha (2+3\beta)} (\lambda \lo n \tsum \lo {j \in  \sf N \hi i \lo n } \|  B \hi {*ij} \lo {m \lo n k \lo n} \| \lo F) \hi 2 }{n \hi i m \lo n \hi 2  k \lo n \hi {4}}+2 \log (n \hi i m \lo n k \lo n)  \Big), \\
\end{split}
\end{equation*}
where we used
$\frac{2}{C \lo {\min}} \lambda \lo n (n \hi i) \hi {3/2} (\tsum \lo {j \in  \sf N \hi i \lo n } \|  B \hi {*ij} \lo {m \lo n k \lo n} \| \lo F) \hi 2 \leq c \lo 2 (b \hi {*i} \lo n)  \hi {s} $ ($s=2,3$)
for the last inequality  (see condition \eqref{cond:propconsistency}).

By the same arguments and using \eqref{t22eq1}, \eqref{bounderrorspine} and  Assumption \eqref{as2theorem2}, we can show the existence of a constant $C \lo {8}> 0$ such that
 \begin{equation*}\label{R62}
 \begin{split}
&P \Big( \max \lo {j \notin \sf N \hi i}\| R \lo {62}  \hi j \| \lo {F} \geq \frac{c \lo 1}{3} { \lambda \lo n} \tsum \lo {j \in  \sf N \hi i \lo n } \|  B \hi {*ij} \lo {m \lo n k \lo n} \| \lo F  \Big)  \\
&\quad \quad \lesssim   \exp \Big( - C \lo {8} \frac{n \hi {1- \alpha (2+3\beta)} \lambda \lo n \hi {-2} (b \hi {*i} \lo n)  \hi {6} }{(n \hi i) \hi {4} m \lo n \hi 2  k \lo n \hi {4} (\tsum \lo {j \in  \sf N \hi i \lo n } \|  B \hi {*ij} \lo {m \lo n k \lo n} \| \lo F) \hi 4}+2 \log (n \hi i m \lo n k \lo n)  \Big)  \\
&\quad \quad  \quad +  \exp \Big( - C \lo {8} \frac{n \hi {1- \alpha (2+3\beta)} (b \hi {*i} \lo n)  \hi {4} }{(n \hi i) \hi {4} m \lo n \hi 2  k \lo n \hi {4} (\tsum \lo {j \in  \sf N \hi i \lo n } \|  B \hi {*ij} \lo {m \lo n k \lo n} \| \lo F) \hi 2}+2 \log (n \hi i m \lo n k \lo n)  \Big) + P(\Omega \hi {\complement})\\
& \quad \quad 
{\lesssim \exp \Big( - C \lo {7} \frac{n \hi {1- \alpha (2+3\beta)} (\lambda \lo n \tsum \lo {j \in  \sf N \hi i \lo n } \|  B \hi {*ij} \lo {m \lo n k \lo n} \| \lo F) \hi 2 }{n \hi i m \lo n \hi 2  k \lo n \hi {4}}+2 \log (n \hi i m \lo n k \lo n)  \Big)
+ P(\Omega \hi {\complement}) ~,}
\end{split}
\end{equation*}
{where the  probability $P(\Omega \hi {\complement})$ can  be estimated by Proposition \ref{splineapproximation} and is dominated by the first term because of Assumption \ref{assnu}.}
 Similarly, using \eqref{t23erroreq1}, we obtain
 \begin{equation*}  %\label{R63}
  \begin{split}
P \Big( \max \lo {j \notin \sf N \hi i}\| R \lo {63} \hi j \| \lo {F} \geq \frac{ c \lo 1}{3} { \lambda \lo n} \tsum \lo {j \in  \sf N \hi i \lo n } \|  B \hi {*ij} \lo {m \lo n k \lo n} \| \lo F  \Big)   
& \lesssim   \exp \Big( - C \lo {9} \frac{n \hi {1- \alpha (2+3\beta)} (b \hi {*i} \lo n)  \hi {6} }{(n \hi i) \hi {4} m \lo n \hi 2  k \lo n \hi {4} (\tsum \lo {j \in  \sf N \hi i \lo n } \|  B \hi {*ij} \lo {m \lo n k \lo n} \| \lo F) \hi 2}+2 \log (n \hi i m \lo n k \lo n)  \Big)  \\
& \quad + \exp \Big( - C \lo {9} \frac{n \hi {1- \alpha (2+3\beta)} (b \hi {*i} \lo n)  \hi {4} }{(n \hi i) \hi {4} m \lo n \hi 2  k \lo n \hi {4} (\tsum \lo {j \in  \sf N \hi i \lo n } \|  B \hi {*ij} \lo {m \lo n k \lo n} \| \lo F) \hi 2}+2 \log (n \hi i m \lo n k \lo n)  \Big) \\
& \quad +  \exp \Big( - C \lo {9} \frac{n  }{n \hi i m \lo n \hi 2  k \lo n }+ \log (n \hi i m \lo n \hi 2 k \lo n)  \Big)\\
& \lesssim \exp \Big( - C \lo {7} \frac{n \hi {1- \alpha (2+3\beta)} (\lambda \lo n \tsum \lo {j \in  \sf N \hi i \lo n } \|  B \hi {*ij} \lo {m \lo n k \lo n} \| \lo F) \hi 2 }{n \hi i m \lo n \hi 2  k \lo n \hi {4}}+2 \log (n \hi i m \lo n k \lo n)  \Big) \\
& \quad +  \exp \Big( - C \lo {9} \frac{n  }{n \hi i m \lo n \hi 2  k \lo n }+ \log (n \hi i m \lo n \hi 2 k \lo n)  \Big).
\end{split}
\end{equation*}
Combining these results, we can conclude that
\begin{equation*}  %\label{R63}
  \begin{split}
P \Big( \max \lo {j \notin \sf N \hi i}\| R \lo {6} \hi j \| \lo {F} \geq \frac{ c \lo 1}{3} { \lambda \lo n} \tsum \lo {j \in  \sf N \hi i \lo n } \|  B \hi {*ij} \lo {m \lo n k \lo n} \| \lo F  \Big)   
& \lesssim \exp \Big( - C \lo {7} \frac{n \hi {1- \alpha (2+3\beta)} (\lambda \lo n \tsum \lo {j \in  \sf N \hi i \lo n } \|  B \hi {*ij} \lo {m \lo n k \lo n} \| \lo F) \hi 2 }{n \hi i m \lo n \hi 2  k \lo n \hi {4}}+2 \log (n \hi i m \lo n k \lo n)  \Big), 
\end{split}
\end{equation*}
}
and  \eqref{deta} holds for $r=6$ (note  that this argument requires condition \eqref{as2theorem2}).

  %%R 7%%
  \noindent
{\it {Term $R \lo 7 \hi j$}}.  Observing \eqref{sigmazetai} and \eqref{m9},
the term $R \lo 7 \hi j$ can be further decomposed as
    \begin{align*}
      R \lo 7 \hi j 
  =& \Sigma \hi n  \lo { \xi \hi j \sf N \hi i \lo n } ( \Sigma \hi n  \lo { \sf N \hi i \lo n  \sf N \hi i \lo n } + \lambda \lo n D \hi * \lo {\sf N \hi i \lo n }) \hi {-1}   \Big( \Sigma \hi n  \lo {\sf N \hi i \lo n  \sf N \hi i \lo n }B \hi {* \sf N \hi i \lo n } \lo n+ \frac{\mathbf{\tilde{H} \lo n} (\xi \hi {\sf N \hi i \lo n })}{n}  w \hi i + \frac{\mathbf{\tilde{H} \lo n } (\xi \hi {\sf N \hi i \lo n })}{n}  \epsilon \hi i   \Big) \\
  &-   \Sigma \hi n  \lo { \xi^j \sf N \hi i \lo n }B \hi {* \sf N \hi i \lo n } \lo n- \frac{\tilde H_n \trans (\xi \hi {j})}{n}  w \hi i - {\frac{\tilde H_n \trans (\xi \hi {j})}{n}  \epsilon \hi i}\\
  =&R \lo {71} \hi j +R \lo {72} \hi j +R \lo {73} \hi j,
     \end{align*}
     where
      \begin{align*}
R \lo {71} \hi j  &= \Sigma  \hi n  \lo { \xi \hi j \sf N \hi i \lo n } ( \Sigma \hi n  \lo { \sf N \hi i \lo n  \sf N \hi i \lo n } + \lambda \lo n D \hi * \lo {\sf N \hi i \lo n }) \hi {-1}  \Sigma \hi n  \lo {\sf N \hi i \lo n \sf N \hi i \lo n }B \hi {* \sf N \hi i \lo n } \lo n-   \Sigma \hi n  \lo { \xi \hi j \sf N \hi i \lo n} B \hi {* \sf N \hi i \lo n } \lo n\\
R \lo {72} \hi j   &= \Sigma \hi n \lo { \xi \hi j \sf N \hi i \lo n } ( \Sigma \hi n  \lo { \sf N \hi i \lo n  \sf N \hi i \lo n } + \lambda \lo n D \hi * \lo {\sf N \hi i \lo n }) \hi {-1}\frac{\mathbf{\tilde{H} \lo n}  (\xi \hi {\sf N \hi i \lo n })}{n}  w \hi i - \frac{{\tilde{H} \lo n} \trans (\xi \hi {j})}{n}  w \hi i \\
R \lo {73} \hi j  &=  \Sigma \hi n  \lo { \xi \hi j \sf N \hi i \lo n } ( \Sigma \hi n  \lo { \sf N \hi i \lo n  \sf N \hi i \lo n } + \lambda \lo n D \hi * \lo {\sf N \hi i \lo n }) \hi {-1}\frac{\mathbf{\tilde{H \lo n} } (\xi \hi {\sf N \hi i \lo n })}{n}  \epsilon \hi i - \frac{{\tilde{H}} \lo n \trans (\xi \hi {j})}{n}  \epsilon \hi i.
      \end{align*}
     
 Turning to $\| R \lo {71} \hi j\|\lo{F}$, we have
      \begin{eqnarray} \label{m5} \nonumber  
\max \lo {j \notin {\sf N \hi i \lo n }}      \| R \lo {71} \hi j \|\lo{F} & \leq & 
%\max \lo {j \notin {\sf N \hi i \lo n }}    \| \Sigma \hi n  \lo { \xi \hi j \sf N \hi i \lo n } ( \Sigma \hi n  \lo { \sf N \hi i \lo n  \sf N \hi i \lo n } + \lambda \lo n D \hi * \lo {\sf N \hi i \lo n }) \hi {-1}  \Sigma \hi n  \lo {\sf N \hi i \lo n  \sf N \hi i \lo n }B \hi {* \sf N \hi i \lo n } \lo n-   \Sigma \hi n  \lo { \xi \hi j  \sf N \hi i \lo n }B \hi {* \sf N \hi i \lo n } \lo n \| \lo {F}\\  \nonumber
%      & = &
      \max \lo {j \notin {\sf N \hi i \lo n }}  \|\Sigma \hi n \lo { \xi \hi j  \sf N \hi i \lo n } ( \Sigma \hi n \lo { \sf N \hi i \lo n  \sf N \hi i \lo n } + \lambda \lo n D \hi * \lo {\sf N \hi i \lo n }) \hi {-1} \left(  \Sigma \hi n  \lo {\sf N \hi i \lo n  \sf N \hi i \lo n }- ( \Sigma \hi n  \lo { \sf N \hi i \lo n  \sf N \hi i \lo n } + \lambda \lo n D \hi * \lo {\sf N \hi i \lo n }) \right) B \hi {* \sf N \hi i \lo n } \lo n \| \lo F \\  \nonumber
     & = & \lambda \lo n \max \lo {j \notin {\sf N \hi i \lo n }} \|\Sigma \hi n  \lo { \xi \hi j \sf N \hi i \lo n }  (\Sigma  \hi n  \lo { \sf N \hi i \lo n  \sf N \hi i \lo n }) \inv D \hi * \lo {\sf N \hi i \lo n } B \hi {* \sf N \hi i \lo n } \lo n \| \lo F \\ 
      & \leq & (1- \eta) \lambda \lo n \tsum \lo {j \in  \sf N \hi i \lo n } \|  B \hi {*ij} \lo {m \lo n k \lo n} \| \lo F ,
\end{eqnarray}    
where we have used $ \Sigma \hi n \lo { \sf N \hi i \lo n  \sf N \hi i \lo n } \preccurlyeq  \Sigma \hi n  \lo { \sf N \hi i \lo n  \sf N \hi i \lo n } +  \lambda \lo n D \hi * \lo {\sf N \hi i \lo n }$ and Assumption \ref{irrepr}.

 For the term $R \lo {72} \hi j$, we first write
    \begin{align*}
  \max \lo {j \notin {\sf N \hi i \lo n }}   \| R \lo {72} \hi j  \|\lo{F} \leq& \max \lo {j \notin {\sf N \hi i \lo n }} \|   \frac{\tilde H \lo n \trans (\xi \hi j)}{ \sqrt{n}} \Big( I - \frac{\mathbf{\tilde{H} \lo n}  (\xi \hi {\sf N \hi i \lo n })}{\sqrt{n}}( \Sigma \hi n  \lo { \sf N \hi i \lo n  \sf N \hi i \lo n } + \lambda \lo n D \hi * \lo {\sf N \hi i \lo n }) \hi {-1}\frac{\mathbf{\tilde{H} \lo n}   (\xi \hi {\sf N \hi i \lo n })}{\sqrt{n}}  \Big) \frac{w \hi i}{\sqrt{n}}  \| \lo F \\ 
      \leq & \max \lo {j \notin {\sf N \hi i \lo n }} \|  \frac{\tilde H \lo n \trans (\xi \hi j)}{ \sqrt{n}} \|\lo{2} \|  \frac{w \hi i}{\sqrt{n}}  \| \lo F \|I - \frac{\mathbf{\tilde{H} \lo n  }(\xi \hi {\sf N \hi i \lo n })}{\sqrt{n}}( \Sigma \hi n  \lo { \sf N \hi i \lo n  \sf N \hi i \lo n } ) \hi {-1}\frac{\mathbf{\tilde{H} \lo n}  (\xi \hi {\sf N \hi i \lo n })}{\sqrt{n}} \| \lo 2.
 \end{align*}  
The matrix $I - \frac{\mathbf{\tilde{H} \lo n}  (\xi \hi {\sf N \hi i \lo n })}{\sqrt{n}}( \Sigma \hi n  \lo { \sf N \hi i \lo n  \sf N \hi i \lo n } ) \hi {-1}\frac{\mathbf{\tilde{H} \lo n}   (\xi \hi {\sf N \hi i} \lo n )}{\sqrt{n}}$ defines a projection  and therefore its operator norm is $1$.
Moreover, by the Lemma  6.2 in \cite{zhou1998local} we have $\| \frac{\tilde H \lo n \trans (\xi \hi j )}{ \sqrt{n}}\frac{\tilde H \lo n (\xi \hi j)}{ \sqrt{n}}) \| \lo 2 \leq \frac{1}{ k \lo n}$ for all $j \in \sf V$.  Thus, observing \eqref{bounderrorspine}, it follows that on the event $\Omega$ defined in \eqref{omega} 
$$ 
\|  \frac{\tilde H \lo n \trans (\xi \hi j)}{ \sqrt{n}} \|\lo{2} \|  \frac{w \hi i}{\sqrt{n}}  \| \lo F \lesssim \frac{n \hi i m \lo n \hi {3/2}}{k \lo n \hi {d+1/2}}.
$$

By condition \eqref{as2theorem2},   we can choose a constant $c \lo 3>0$ such that
\begin{align*}
\max \lo {j \notin {\sf N \hi i \lo n} } \| R^j \lo {72} \| \lo {F} \leq \frac{n \hi i m \lo n \hi {3/2}}{k \lo n \hi {d+1/2}} \leq c \lo 3 \sqrt{\frac{2}{C \lo {\min}}} { \lambda \lo n} \tsum \lo {j \in  \sf N \hi i \lo n } \|  B \hi {*ij} \lo {m \lo n k \lo n} \| \lo F.
\end{align*}
  {Combining this inequality with  the constant $c \lo 3=\sqrt{\frac{C \lo {\min}}{2}} (1-\eta-\frac{c \lo 1}{2})$ and \eqref{m5} now yields
}
{\begin{equation*}
\begin{split}
&P \Big( \max \lo {j \notin {\sf N \lo n  \hi i}} \| R \lo {7} \hi j \| \lo {F} \geq c \lo 1 \lambda \lo n \tsum \lo {j \in  \sf N \hi i \lo n } \|  B \hi {*ij} \lo {m \lo n k \lo n} \| \lo F   \Big) \leq P \Big( \max \lo {j \notin {\sf N \lo n  \hi i}} \| R \lo {73} \hi j \| \lo {F} \geq \frac{c \lo 1}{2} \lambda \lo n \tsum \lo {j \in  \sf N \hi i \lo n } \|  B \hi {*ij} \lo {m \lo n k \lo n} \| \lo F   \Big)\\
& \quad \leq  P  \Big(   \| \frac{\mathbf{\tilde{H} \lo n}  \trans (\xi \hi {\sf N \hi i \lo n })}{n}  \epsilon \hi i \| \lo F \geq \frac{c \lo 1}{4} \frac{ \sqrt{n \hi i}}{(1-\eta)}   { \lambda \lo n} \tsum \lo {j \in  \sf N \hi i \lo n } \|  B \hi {*ij} \lo {m \lo n k \lo n} \| \lo F   \Big) +(p-n \hi i) P  \Big( \| \frac{\tilde H \lo n \trans (\xi \hi { j})}{n}  \epsilon \hi i \| \lo F \geq \frac{c \lo 1}{4} { \lambda \lo n} \tsum \lo {j \in  \sf N \hi i \lo n } \|  B \hi {*ij} \lo {m \lo n k \lo n} \| \lo F     \Big)\\
&  \quad \lesssim     \exp  \Big( -C \lo {8} \frac{n \left( \lambda \lo n \tsum \lo {j \in  \sf N \hi i \lo n } \|  B \hi {*ij} \lo {m \lo n k \lo n} \| \lo F \right) \hi 2}{ n \hi i m \lo n \hi 2 k \lo n}+ \log (n \hi i m \lo n \hi 2 k \lo n)  \Big)\\
& \quad \quad+  \exp  \Big( -C \lo {8} \frac{n \left( \lambda \lo n \tsum \lo {j \in  \sf N \hi i \lo n } \|  B \hi {*ij} \lo {m \lo n k \lo n} \| \lo F \right) \hi 2}{ m \lo n \hi 2 k \lo n}+ \log ((p-n \hi i)  m \lo n \hi 2 k \lo n)  \Big) + P(\Omega \hi {\complement}).
\end{split}
\end{equation*} }
Finally, combining  the result \eqref{step1} from Step 1  with the estimates for $R_1, \ldots, R_7$, we conclude that
\begin{equation}
\begin{split}
& P \Big( \max \lo {j \notin \sf N \hi i \lo n  } \|\hat \Sigma \hi n \lo { \xi \hi j  {\sf N \hi i \lo n }} \hat B \hi {{\sf N \hi i \lo n }} \lo n - \hat {\Sigma}  \hi n \lo { \xi \hi j {\xi \hi i}} \| \lo {F} \geq \lambda \lo n \tsum \lo {j \ne i} \hi {p} \| \hat B \hi {ij} \lo n \| \lo F \Big) \\
& \quad \quad  \quad \lesssim  \exp \Big( - C \lo 1 \frac{n \hi {1- \alpha (2+3\beta)} (\lambda \lo n \tsum \lo {j \in  \sf N \hi i \lo n } \|  B \hi {*ij} \lo {m \lo n k \lo n} \| \lo F) \hi 2  }{n \hi i m \lo n \hi 2  k \lo n \hi {4} }+2 \log (p  m \lo n k \lo n) \Big), 
\end{split} \nonumber 
  \end{equation}
  such that   $ \lambda \lo n (n \hi i) \hi {3/2} \lesssim (b \hi {*i} \lo n)  \hi {3} (\tsum \lo {j \in  \sf N \hi i \lo n } \|  B \hi {*ij} \lo {m \lo n k \lo n} \| \lo F) \hi {-2}$.
     \eop

  %%%%%%%%%%%%%%%%%%%%%%%%%%%%%%%%%%%%%%%%%%%
  % theorem 3
 \paragraph{Proof of Proposition \ref{sampletheorem}  }
  First, note that the event $\hat {\sf N \hi i \lo n}= \sf N \hi i \lo n$ holds if and only if
$$
  \|  \hat B \hi { ij} \lo n \| \lo F \ne 0   \quad \forall j \in \sf N \hi i \lo n  \quad  \mbox{ and }  \quad 
    \|  \hat B \hi { ij} \lo n \| \lo F = 0  \quad \forall j \notin \sf N \hi i \lo n,
$$  
  which is implied by the conditions
   \begin{align*}
      %\begin{cases}
&\|  \hat B \hi { \sf N \hi i \lo n}  \lo n - B \hi {* \sf N \hi i \lo n } \lo n \| \lo F < b \hi {*i} \lo n  =  \min \lo {j \in \sf N \hi i} \| B \hi {*ij} \lo {m_n k_n} \| _F, \\
& \max \lo {j \notin \sf N \hi i \lo n} \| \hat \Sigma  \hi n \lo { \xi \hi j \sf N \hi i} \hat B \hi {\sf N \hi i \lo n} \lo n - \hat {\Sigma}  \hi n \lo { \xi \hi j {\xi \hi i}} \| \lo F  \leq \lambda \lo n \tsum \lo {j \in  \sf N \hi i} \| \hat B \hi {ij} \lo n \| \lo F.
     %\end{cases}
    \end{align*}  
     Thus,  the Proposition \ref{propconsistency} (with $\delta=b^{*i}_n$, an obvious estimate of the probability using condition \eqref{cond:propconsistency}) and Proposition \ref{bhatoptim} we can  conclude that  
     \begin{equation*}
\begin{split}
 P( \hat {\sf N \hi i \lo n} \ne {\sf N \hi i \lo n} \text{ and } \mathcal{N} )   \lesssim   \exp \Big( - C \lo 1 \frac{n \hi {1- \alpha (2+3\beta)} (\lambda \lo n \tsum \lo {j \in  \sf N \hi i \lo n } \|  B \hi {*ij} \lo {m \lo n k \lo n} \| \lo F) \hi 2  }{n \hi i m \lo n \hi 2  k \lo n \hi {4} }+2 \log (p m \lo n k \lo n) \Big)  
 \end{split}
  \end{equation*}     
 and this completes the proof of Proposition \ref{sampletheorem}.
 \eop

\end{document}